\documentclass[11pt]{article}

\usepackage[a4paper,margin=1in]{geometry}
\usepackage{amsmath}
\usepackage{amssymb}
\usepackage{amsthm}
\usepackage{mathtools}
\usepackage{microtype}
\usepackage{graphicx}
\graphicspath{{fig/}}
\usepackage{xcolor}
\usepackage{tikz}
\usetikzlibrary{arrows.meta}
\usepackage{enumitem}
\setlist[enumerate]{leftmargin=.5in}
\setlist[itemize]{leftmargin=.5in}
\usepackage[labelsep=period,font=small,labelfont=bf]{caption}
\definecolor{linkblue}{rgb}{0.1,0.2,0.55}
\usepackage[hyphens]{url}
\usepackage[colorlinks,allcolors=linkblue,hypertexnames=false,breaklinks,hyperfootnotes=false]{hyperref}
\usepackage[capitalize,nameinlink]{cleveref}
\usepackage{etoolbox}
\apptocmd{\thebibliography}{\raggedright}{}{\PackageWarning{GQS}{thebibliography not patched}}

\DeclareMathOperator{\sign}{sgn}
\long\def\tmpcomment#1\endtmpcomment{#1}

\theoremstyle{plain}
\newtheorem{theorem}{Theorem}[section]
\newtheorem{lemma}[theorem]{Lemma}

\newtheorem{proposition}[theorem]{Proposition}
\theoremstyle{definition}
\newtheorem{definition}[theorem]{Definition}
\theoremstyle{remark}
\newtheorem{remark}[theorem]{Remark}

\numberwithin{equation}{section}

\newcommand{\email}[1]{\href{mailto:#1}{\nolinkurl{#1}}}
\newcommand{\TheTitle}{Generalized Quasi-Symmetric Nets: A Constructive Approach to the Equimodular Elliptic Type and to Flexible $m \times n$ Nets}
\newcommand{\TheAuthors}{A. Nurmatov, F. Rist, J. Klein, and D. L. Michels}

\title{\TheTitle}
\author{A.~Nurmatov\textsuperscript{1}, F.~Rist\textsuperscript{1}, J.~Klein\textsuperscript{1}, and D.~L.~Michels\textsuperscript{1}}
\date{}

\hypersetup{pdftitle={Generalized Quasi-Symmetric Nets}, pdfauthor={\TheAuthors}}

\begin{document}

\maketitle
% The common affiliation and the e-mail addresses: an unmarked footnote on the first page, keyed by the superscripts 1 of the author line.
\begingroup
\makeatletter
\renewcommand\@makefntext[1]{\noindent\raggedright\settowidth{\hangindent}{\textsuperscript{1}}\hangafter=1 #1}%
\footnotetext{\textsuperscript{1}KAUST Computational Sciences Group, Campus Building 1, Thuwal 23955-6900, Saudi Arabia\newline
  \emph{Email addresses:} \begin{tabular}[t]{@{}l@{\hspace{1.5em}}l@{}}
    \email{abdukhomid.nurmatov@kaust.edu.sa}, & \email{florian.rist@kaust.edu.sa},\\
    \email{jonathan.klein@kaust.edu.sa}, & \email{dominik.michels@kaust.edu.sa}
  \end{tabular}}
\makeatother
\endgroup

\begin{abstract}
Quadrilateral meshes that flex with rigid planar faces are of interest to
transformable architecture, rigid-foldable origami, and mechanical
metamaterials. Their building blocks are flexible $3 \times 3$ nets, that
is, Kokotsakis polyhedra with a quadrangular base, classified by
Izmestiev. We study one of his types, the equimodular elliptic type,
defined by three conditions on the flat angles, the third of them through
Jacobi elliptic functions. We turn these conditions into equivalent
algebraic ones, free of elliptic functions, which lead to the generalized
quasi-symmetric nets (GQS-nets), polyhedra of this type that include the
known elliptic quasi-symmetric nets. For every polyhedron of this type,
we give a flexible polyhedron with the same flat angles and its flexion
in closed form. From the simplest GQS-nets, we assemble arbitrarily large
flexible $m \times n$ nets. Examples include nets that, at one position,
form letters; further nets in the supplementary materials form digits or,
seen along their rows, approximately follow the side silhouettes of a bunny
and a fish. Unlike numerically computed
mechanisms, these nets are exactly flexible. 
As an application, we design
three such nets as transportable
arch halls of $80$~mm sandwich panels, which fold to fit through the
door of a standard container or into the well of a low-bed
semi-trailer, and we check them numerically under snow load and build a
$1{:}20$ model of one. Supplementary
Mathematica notebooks verify the symbolic computations, a Python
tool constructs the nets and checks them numerically, including for
self-intersections, and further scripts recheck the halls.

\end{abstract}

\section{Introduction and prior work}
\label{sec:intro}
A quadrilateral mesh with planar faces is a mechanism with the shape of
a surface if it can deform while every face stays rigid; such a
deformation is called a flexion. These mechanisms are of interest to transformable
architecture~\cite{Jiang2024}, to rigid-foldable
origami~\cite{DemaineORourke2007, Tachi2009, HePhD2021}, and to
mechanical metamaterials~\cite{Silverberg2014, Dudte2016, Callens2018,
Dieleman2020}; in architecture, for example, a facade or a roof may pass
from a flat to a curved state~\cite{Tachi2010freeform}. 
The mathematical theory of these mechanisms is that of flexible
polyhedra, which goes back to Bricard's octahedra~\cite{Bricard1897} and, for
quadrilateral meshes, to Sauer and Graf~\cite{SauerGraf1931} and
Kokotsakis~\cite{Kokotsakis1933}.

\begin{figure}[!t]
    \centering
    \includegraphics[width=0.48\textwidth]{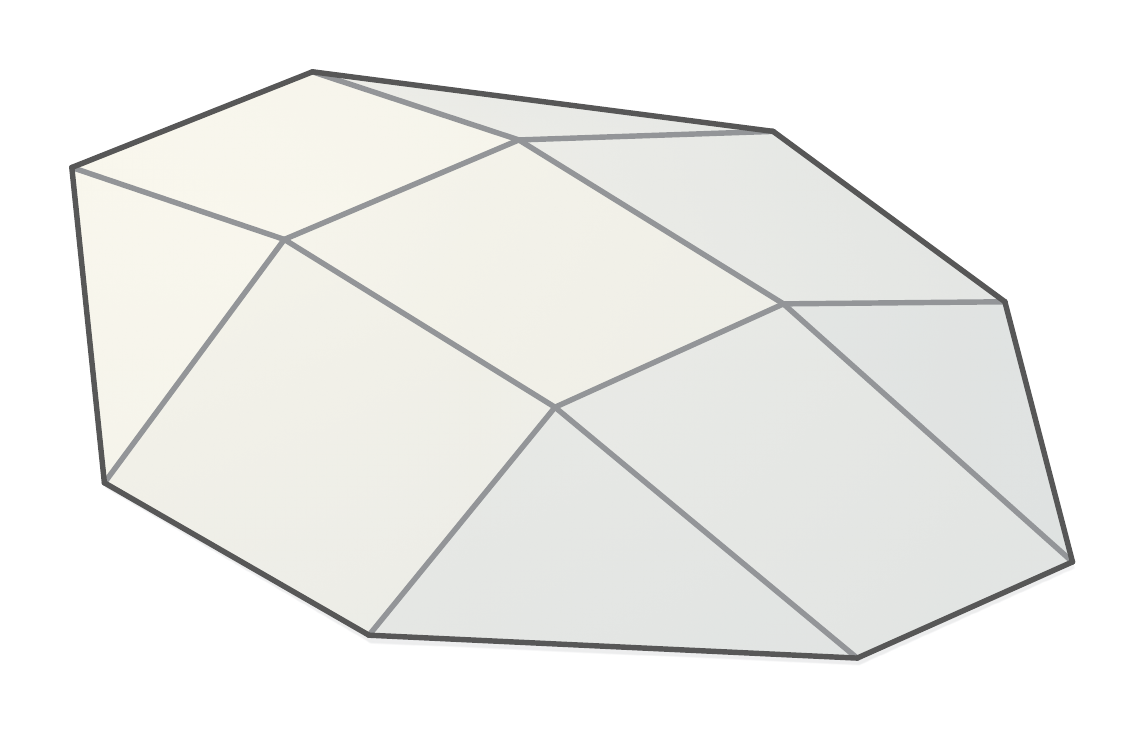}
    \hfill
    \includegraphics[width=0.48\textwidth]{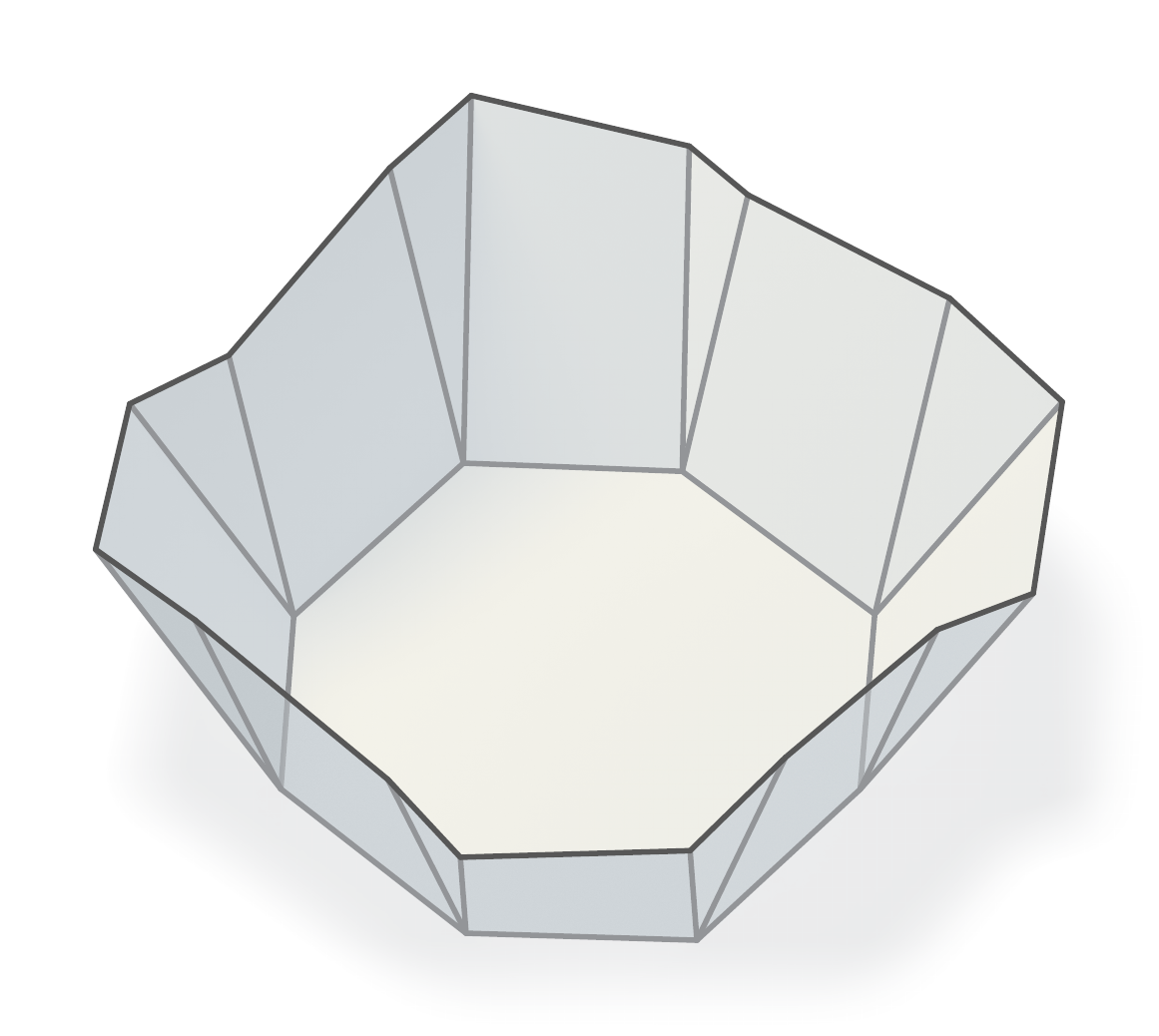}%
    \caption{Kokotsakis polyhedra with a quadrangular base (left) and an
octagonal base (right).}
    \label{fig:01ab}
\end{figure}

Kokotsakis polyhedra consist of a central polygon (the base), one
quadrilateral attached to each of its sides, and one triangle between
consecutive quadrilaterals, see~\cite[Section~1.1]{izmestiev-2017} and
Figure~\ref{fig:01ab}. Those with a quadrangular base and convex faces
(Figure~\ref{fig:01ab}, left) are the $3 \times 3$ nets of
Section~\ref{sec:notation} (grids of nine convex quadrilaterals) with
their four corner points omitted. These nets are the building blocks of
larger flexible meshes: an $m \times n$ net with no two coplanar
adjacent faces is flexible if all of its $3 \times 3$ subnets are
flexible~\cite[Theorem~5]{IzmestievNawratilSamarehRadSharifmoghaddam2026}.
This criterion corrects Theorem~3.2 of~\cite{SchiefBobenkoHoffmann2008},
whose hypothesis is insufficient, as shown
by~\cite[Theorems~6 and~7]{IzmestievNawratilSamarehRadSharifmoghaddam2026}.

The flexible $3 \times 3$ nets themselves were classified by
Izmestiev~\cite{izmestiev-2017}, after earlier work~\cite{SauerGraf1931,
SchiefBobenkoHoffmann2008, Karpenkov2010, Stachel2010, Nawratil2011,
Nawratil2012}. Up to switching boundary strips, an operation that
preserves flexibility~\cite[Section~2.5]{izmestiev-2017}, they fall into
eight classes: orthodiagonal, isogonal, equimodular, conjugate-modular, linear
compound, linearly conjugate, chimera, and trivial. 
The classes split into subclasses given by equations on the flat angles (the face angles at
the base vertices), derived from the complexified configuration space.
As Izmestiev points out~\cite[Section~3]{izmestiev-2017}, this leaves
open whether every subclass contains a flexible polyhedron, whether the
real part of the configuration space is nonempty, and whether it
contains polyhedra without self-intersections.

Flexible polyhedra are known explicitly in some of these classes. The
orthodiagonal and isogonal classes contain the T-nets and V-nets of
Sauer and Graf~\cite{SauerGraf1931}, whose flexions have been
studied~\cite{SharifmoghaddamMaleczekNawratil2023,
IzmestievRasoulzadehTervooren2024,
KilianNawratilRaffaelliRasoulzadehSharifmoghaddam2024}. The
orthodiagonal class also contains the anti-involutive polyhedra, which
Erofeev and Ivanov showed to exist and to be flexible, with explicit
flexions~\cite{ErofeevIvanov2020}. Some of these classes have been
extended to meshes with skew faces~\cite{AikynEtAl2024, Nawratil2022,
Nawratil2023}, and Kokotsakis polyhedra with skew faces have been
studied algebraically~\cite{LiuOuyangMichels2024, Liu2026}. 
Recently, the P-nets were introduced: flexible discrete and semi-discrete surfaces
whose spatial shapes are controlled by three
polylines~\cite{Nawratil2024Pnets}.

In origami, He and Guest~\cite[Section~5]{HeGuest2020} built large
rigid-foldable quadrilateral meshes from $3 \times 3$ nets of
equimodular conic and linear compound type,
see~\cite[Section~4.3]{HePhD2021} for the correspondence with these
types. The pluripotent
tiling patterns of~\cite{Dieleman2020} are a special case of the
stitching of nets of equimodular conic
type~\cite[Section~4.6.7]{HePhD2021}, and large rigid-foldable meshes
have also been obtained by repetitive stitching of units of linear and
equimodular couplings~\cite[Section~3.1]{HeHayakawaOhsaki2026}. Recent
work gives explicit real kinematics for equimodular couplings of
developable vertices and for developable triangle-based Kokotsakis
polyhedra~\cite{HuChen2026} and studies rigid-foldable developable
polygons of equimodular type~\cite{HuLinZhengBi2026}. At a developable
vertex, however, the flat angles sum to~$2\pi$, which excludes the
equimodular elliptic type, the subject of this paper.

Other approaches are approximate or numerical. Approximately flexible
polyhedra have been constructed~\cite{GorkavyiMilka2018}, and large
quadrilateral-mesh mechanisms are computed numerically by
optimization~\cite{Jiang2024}, for instance starting from flexible nets
of isotropic geometry~\cite{PirahmadPottmannSkopenkov2025, Yorov2026}.
For two mechanisms obtained in this way, Yorov et al.\ observed that all
$3 \times 3$ subnets are of equimodular or linear compound type, within
a numerical precision that does not distinguish the equimodular elliptic
type from its limiting cases~\cite[Section~3]{Yorov2026}.

The equimodular elliptic type itself is defined, for polyhedra of
elliptic type, by three conditions on the flat angles, the third of them
through Jacobi elliptic functions~\cite[Section~3.3.1]{izmestiev-2017},
see Section~\ref{sec:notation}. The quasi-symmetric nets (QS-nets)
of~\cite{NURMATOV2026104102}, defined by equalities and complements
to~$\pi$ between their flat angles, have this type whenever they are
elliptic~\cite[Theorem~1]{NURMATOV2026104102}, and the flexible QS-net
of~\cite[Example~4.1 and Proposition~2]{NURMATOV2026104102} shows that
the type is not empty, which answers Izmestiev's first
question~\cite[Section~3]{izmestiev-2017} for this type. 
The paper~\cite{NURMATOV2026104102} also characterized the polyhedra of
the type with dihedral angles at the base in~$(0, \pi)$ by algebraic
relations between their flat and dihedral angles, but kept the third
condition in its transcendental
form~\cite[Proposition~1]{NURMATOV2026104102}.

Assembling larger nets from polyhedra of this type involves three tasks:
recognizing the type, knowing how these polyhedra move, and choosing
neighbouring $3 \times 3$ subnets that share flat angles, fit together, and
move compatibly. This paper advances all three beyond~\cite{izmestiev-2017}
and~\cite{NURMATOV2026104102}. First, we rewrite the conditions of
the type in equivalent forms: the first two become simpler equations
(Proposition~\ref{prop:eqMrs}), and the third, which cannot be tested by
elementary algebra, becomes an algebraic condition, free of elliptic functions
(Propositions~\ref{prop:shifts} and~\ref{prop:signs}), which
leads to the generalized quasi-symmetric
nets (GQS-nets) of Definition~\ref{def:gqs}, polyhedra of this type
(Theorem~\ref{th:gqs}) that include the elliptic QS-nets (Section~\ref{sec:results}).
Second, closed-form flexions were known for the flat angles of the elliptic
QS-nets~\cite[Theorem~1]{NURMATOV2026104102} and of one further
polyhedron~\cite[Example~4.3 and Proposition~4]{NURMATOV2026104102}; we give
them for the flat angles of every polyhedron of the type
(Theorem~\ref{th:main}), which answers
Izmestiev's second question~\cite[Section~3]{izmestiev-2017} for this type.
Third, both earlier works concern single $3 \times 3$ nets; we assemble
GQS-nets into flexible $m \times n$ nets for all $m, n \geq 3$, with flexions
in closed form (Proposition~\ref{prop:ex:flex}).

As examples, the three nets of Figure~\ref{fig:ex:nets:GQS} each form, at
one position, one letter of ``GQS''; Section~\ref{sm:ex} of the
supplementary materials adds nets that form the digits of ``2026'' and nets
that, seen along their rows of blocks,
approximately follow the side silhouettes of the Stanford
bunny~\cite{StanfordRepo} and of Blub the fish~\cite{CraneRepo}. The
motions between the positions shown are checked numerically to be free
of self-intersections. Unlike the mechanisms computed numerically
in~\cite{Jiang2024, Yorov2026}, these nets are exactly flexible, with
flexions in closed form. In the repetitive stitching
of~\cite[Section~3.1]{HeHayakawaOhsaki2026}, the flat angles are
periodic and solved numerically, and flexibility is inferred from
Theorem~3.2 of~\cite{SchiefBobenkoHoffmann2008}; here, the flat angles are
explicit and need not be periodic, and flexibility follows from these flexions.

As an application, we design three further such nets as transportable
arch halls that fold to fit through the door of a standard container or
into the well of a low-bed semi-trailer, a use that needs a motion of rigid
faces but no flat position, and we check their clear height and snow load
numerically and build a $1{:}20$ physical model of one of them
(Section~\ref{sec:application}).

Section~\ref{sec:notation} opens with an informal overview of these results
and then fixes the notation; Sections~\ref{sec:results}--\ref{sec:concl}
contain the main results, their proofs, the flexible $m \times n$ nets with
examples, the application, the discussion, and the conclusion. The
supplement~\cite{Nurmatov2026Repo}, an online repository, contains Mathematica notebooks
(\texttt{helper.nb}, \texttt{ehelper.nb}) that verify the symbolic
computations of Sections~\ref{sec:proofs} and~\ref{sec:examples}; a Python
tool (\texttt{gqs.py}) that constructs the nets, shows their motion, and
checks them numerically, including for self-intersections; the data and
checks of the halls (\texttt{Application}); the print files and the
assembly guide of the physical model (\texttt{Prototype}); and a short
video that summarizes the paper and shows the nets and the physical
model in motion (\texttt{Video}).
The supplementary materials contain
the data of the nets of Figure~\ref{fig:ex:nets:GQS}, the further examples
(Section~\ref{sm:ex}), and the details of the application
(Section~\ref{sm:app}).

\section{Notation}
\label{sec:notation}
In this section, we fix the notation, recall the equimodular elliptic
type of~\cite{izmestiev-2017}, define the GQS-nets, and state Bricard's
system. We first explain informally what these notions describe and how
the results of this paper address the three tasks of
Section~\ref{sec:intro}; the precise definitions follow.

The first task, recognizing the equimodular elliptic type, concerns its
definition, which consists of three conditions on the sixteen flat angles
of a polyhedron of elliptic type. The first two are equations
that mix sines of the flat angles and of their combinations;
Proposition~\ref{prop:eqMrs} rewrites them as eight simpler, unmixed equations between
ratios of products of two sines, each linking two neighbouring vertices.
Section~\ref{sec:examples} uses this form to assemble flexible $m \times n$ nets.
The third condition is defined through Jacobi elliptic functions. If the
first two conditions hold, Proposition~\ref{prop:shifts} replaces it by
an equivalent algebraic condition with three radicals.

This algebraic condition also shows where the GQS-nets come from. 
When the flat angles at two pairs of vertices satisfy one of four relations,
its three radicals cancel each other or reduce to expressions without
radicals, and the condition holds. These relations, together
with the first two conditions, define the GQS-nets of
Definition~\ref{def:gqs}, and
Theorem~\ref{th:gqs} shows that every GQS-net has equimodular elliptic
type. In the simplest case, the four flat angles at a single vertex
determine all sixteen through equalities and complements to~$\pi$, see
Lemma~\ref{lem:ex:patterns}.

For the second task, knowing how the polyhedra of this type move,
Theorem~\ref{th:main} shows that the flat angles of every such
polyhedron, in particular of every GQS-net, are the flat angles of a
flexible polyhedron. It gives a
flexion of this polyhedron by explicit formulas, in which one dihedral
angle along an edge of the central face serves as the parameter and the
other three such angles are explicit functions of it. These formulas
give solutions of Bricard's system of Definition~\ref{def:bricard},
which relates these four dihedral angles to the flat angles. The theorem
also requires four signs, each equal to $1$ or $-1$, for which the third
condition holds; Proposition~\ref{prop:signs} determines these signs by
elementary algebra, and Section~\ref{sec:results} gives them a geometric
meaning.

The third task, choosing neighbouring $3 \times 3$ subnets of an
$m \times n$ net, involves three requirements: neighbouring subnets must
share flat angles, fit together, and move compatibly. For the simplest
GQS-nets in Section~\ref{sec:examples}, all three requirements can be
met: by Proposition~\ref{prop:ex:flex}, it suffices to choose four
suitable numbers in~$(0, \pi)$ as the flat angles at one vertex and, for
each row of blocks, one of the four patterns of
Lemma~\ref{lem:ex:patterns}; these choices determine the flat
angles of a flexible $m \times n$ net, with a flexion given on each
$3 \times 3$ subnet by the formulas of Theorem~\ref{th:main}.
To build intuition before the technical part,
Figure~\ref{fig:ex:nets:GQS} shows three such nets.

We now turn to the precise definitions.

\begin{figure}[htbp!]
        \centering
        \input{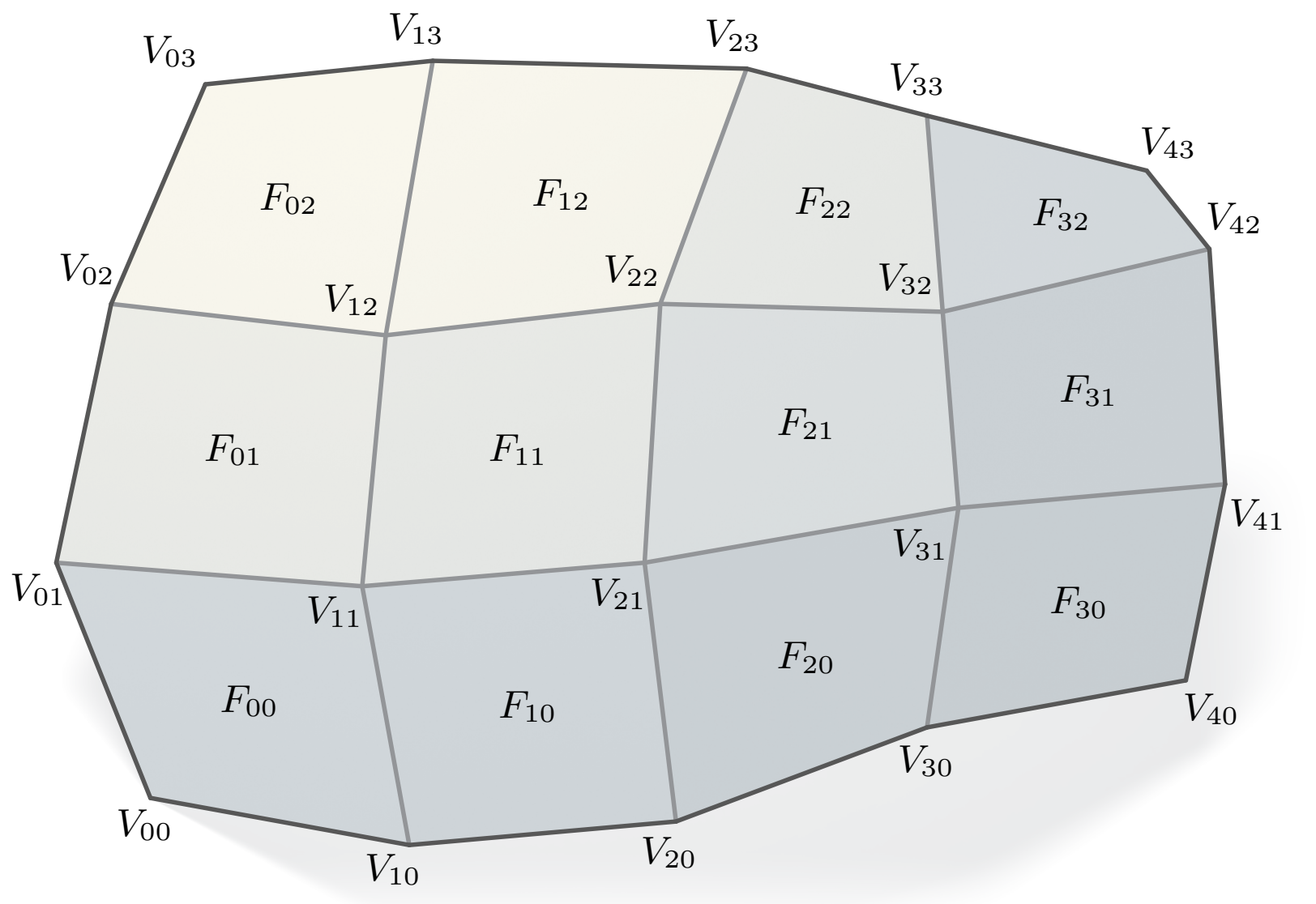}%
    \caption{Notation for vertices and faces of a $4 \times 3$ net.}
    \label{fig:02}
\end{figure}

An $m \times n$ \emph{net} (\emph{with convex faces}) 
is a collection of $(m+1)(n+1)$ points $V_{ij} \in \mathbb{R}^3$, 
indexed by integers $0 \leq i \leq m$ and $0 \leq j \leq n$, 
such that, for all $0 \leq i < m$ and $0 \leq j < n$, 
the points $V_{ij}, V_{i+1,j}, V_{i+1,j+1}, V_{i,j+1}$ 
form the consecutive vertices of a convex quadrilateral, 
denoted $F_{ij}$. See Figure~\ref{fig:02}. 
The vertices $V_{ij}$ with $0 < i < m$ and $0 < j < n$ are the
\emph{inner vertices} of the net.
The \emph{flat angles} of the net are the angles of its faces at the inner
vertices, four at each of them.

An $m \times n$ net $V_{ij}$ is \emph{flexible} 
if it belongs to a continuous family 
(called a \emph{flexion}) of $m \times n$ nets $V_{ij}(t)$, 
where $t \in [0, \epsilon)$ for some $\epsilon > 0$, 
such that all corresponding quadrilaterals $F_{ij}(t)$ 
are congruent to $F_{ij}(0)$ for all $t \in [0, \epsilon)$, 
$V_{ij}(0)$ is the original net $V_{ij}$, 
and the nets $V_{ij}(t)$ are pairwise non-congruent.

In what follows, unless stated otherwise, 
we consider only $3\times 3$ nets and call 
them simply \emph{polyhedra}. 
The corner points $V_{00}$, $V_{30}$, $V_{03}$, and $V_{33}$ 
are omitted as irrelevant to flexibility.

\begin{figure}[!htbp]
    \centering
    \input{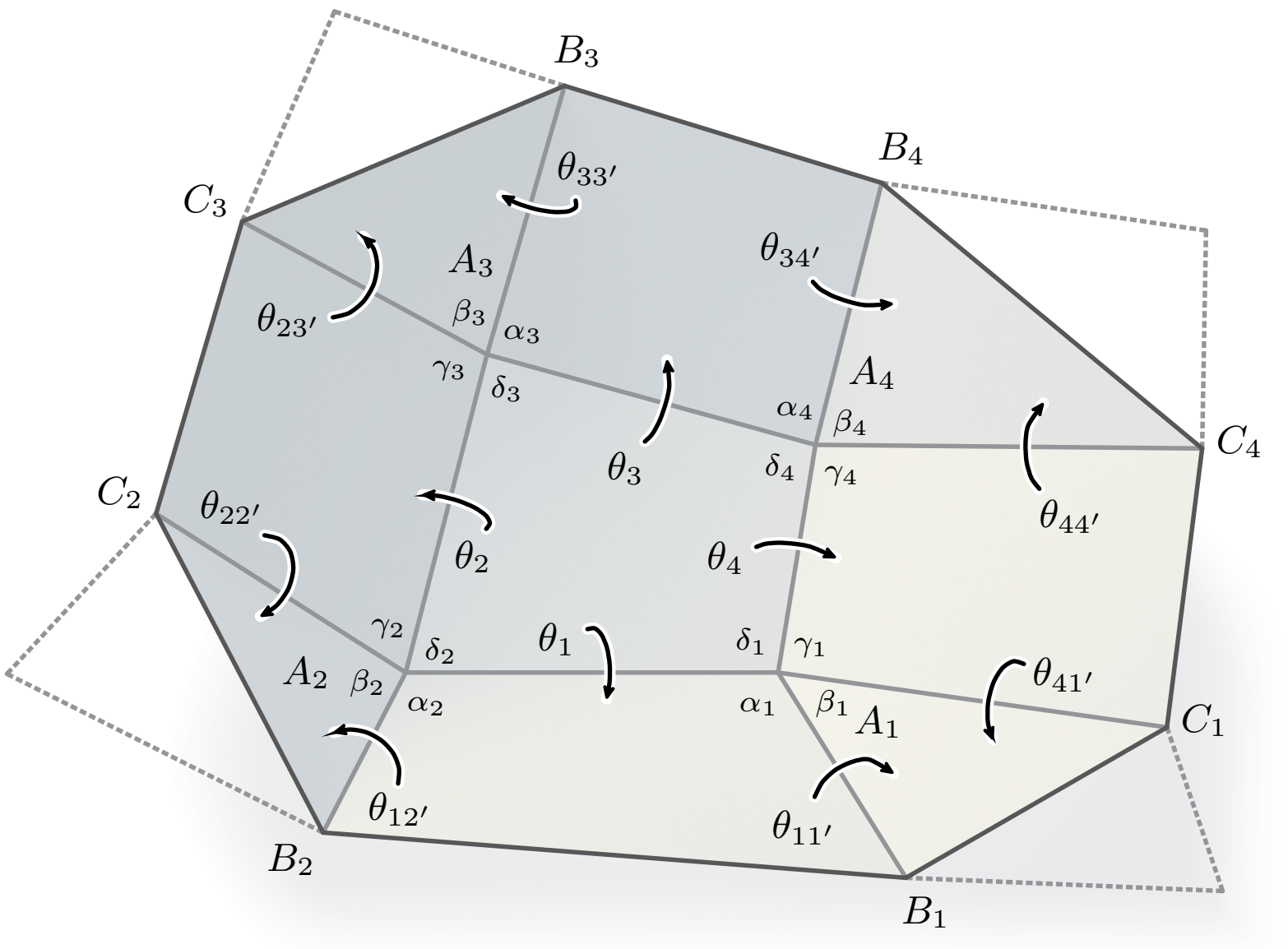}%
    \caption{Notation for vertices, flat and dihedral angles of a
$3 \times 3$ net; dashed segments indicate the omitted corners.}
    \label{fig:03}
\end{figure}

We label the vertices as in Figure~\ref{fig:03}:
$A_2\coloneqq V_{11}$, 
$B_2\coloneqq V_{10}$, $C_2\coloneqq V_{01}$, etc.,
and set $A_0 \coloneqq A_4$ and $A_5 \coloneqq A_1$. 
The flat angles at $A_i$, $i = 1, \ldots, 4$, are 
$\alpha_i \coloneqq \angle A_{i-(-1)^i}A_iB_i$, 
$\beta_i \coloneqq \angle B_iA_iC_i$, 
$\gamma_i \coloneqq \angle C_iA_iA_{5-i}$, and 
$\delta_i \coloneqq \angle A_{i-1}A_iA_{i+1}$.
By definition, $\alpha_i$, $\beta_i$, $\gamma_i$, $\delta_i \in (0, \pi)$.

The oriented normals of the faces are
$\vec{n}_0 \coloneqq \overrightarrow{A_2A_1}\times\overrightarrow{A_2A_3}$, 
$\vec{n}_1 \coloneqq \overrightarrow{A_2A_1}\times\overrightarrow{A_2B_2}$,
$\vec{n}_2 \coloneqq  \overrightarrow{A_2C_2}\times\overrightarrow{A_2A_3}$,
$\vec{n}_3 \coloneqq  \overrightarrow{A_3B_3}\times\overrightarrow{A_3A_4}$,
$\vec{n}_4 \coloneqq  \overrightarrow{A_1A_4}\times\overrightarrow{A_1C_1}$,
$\vec{n}_{1'} \coloneqq  \overrightarrow{A_1B_1}\times\overrightarrow{A_1C_1}$,
$\vec{n}_{2'} \coloneqq  \overrightarrow{A_2C_2}\times\overrightarrow{A_2B_2}$,
$\vec{n}_{3'} \coloneqq  \overrightarrow{A_3B_3}\times\overrightarrow{A_3C_3}$,
$\vec{n}_{4'} \coloneqq  \overrightarrow{A_4C_4}\times\overrightarrow{A_4B_4}$.
For nonzero vectors $\vec{w}_1$, $\vec{w}_2$, $\vec{w}_3 \in \mathbb{R}^3$,
we define the oriented angle $\theta(\vec{w}_1,\vec{w}_2,\vec{w}_3) \in [-\pi,\pi)$ by
\begin{equation*}
    \theta(\vec{w}_1,\vec{w}_2,\vec{w}_3) \coloneqq
  \begin{cases}
    \phantom{-}\angle(\vec{w}_1,\vec{w}_2)\in(0,\pi),
      & \text{if } \det(\vec{w}_1,\vec{w}_2,\vec{w}_3)>0,\\[2pt]
    -\angle(\vec{w}_1,\vec{w}_2)\in[-\pi,0],
      & \text{if } \det(\vec{w}_1,\vec{w}_2,\vec{w}_3)\leq 0,
  \end{cases}
\end{equation*}
where $\angle(\vec{w}_1,\vec{w}_2) \in [0,\pi]$ denotes the unoriented angle.

The (\emph{oriented}) \emph{dihedral angles}
of the $3 \times 3$ net are the four \emph{central} ones,
$\theta_i \coloneqq \theta(\vec{n}_0,\vec{n}_i,\overrightarrow{A_{i+1}A_i})$,
$i = 1, \ldots, 4$, along the edges of the central face, and the eight
\emph{corner} ones, along the edges $A_iB_i$ and $A_iC_i$ of the corner
faces:
$\theta_{11'} \coloneqq \theta(\vec{n}_1,\vec{n}_{1'},\overrightarrow{B_1A_1})$,
$\theta_{12'} \coloneqq \theta(\vec{n}_1,\vec{n}_{2'},\overrightarrow{A_2B_2})$,
$\theta_{22'} \coloneqq \theta(\vec{n}_2,\vec{n}_{2'},\overrightarrow{C_2A_2})$,
$\theta_{23'} \coloneqq \theta(\vec{n}_2,\vec{n}_{3'},\overrightarrow{A_3C_3})$,
$\theta_{33'} \coloneqq \theta(\vec{n}_3,\vec{n}_{3'},\overrightarrow{B_3A_3})$,
$\theta_{34'} \coloneqq \theta(\vec{n}_3,\vec{n}_{4'},\overrightarrow{A_4B_4})$,
$\theta_{41'} \coloneqq \theta(\vec{n}_4,\vec{n}_{1'},\overrightarrow{A_1C_1})$,
$\theta_{44'} \coloneqq \theta(\vec{n}_4,\vec{n}_{4'},\overrightarrow{C_4A_4})$;
see Figure~\ref{fig:03}.

For $j = 0, 1, 2, 3$, the \emph{cyclic permutation of the indices by $j$}
is the relabelling of $A_{i+j}$ as $A_i$ and of $B_{i+j}$, $C_{i+j}$ as
$B_i$, $C_i$ for even $j$ and as $C_i$, $B_i$ for odd $j$, where $i + j$
is taken modulo~$4$, with values in $\{1, \ldots, 4\}$.
It relabels $\theta_{i+j}$ as $\theta_i$, and, for odd
$j$, it exchanges the roles of $\alpha_i$ and $\gamma_i$.

A polyhedron has \emph{elliptic type} (or is an \emph{elliptic polyhedron}) if
\begin{equation}\label{eq:N1}
\alpha_i \pm \beta_i \pm \gamma_i \pm \delta_i \neq 0 \pmod{2\pi} 
\end{equation}
for all choices of the signs $\pm$ and each $i=1, \ldots, 4$. 
In what follows, we consider elliptic polyhedra only, unless stated
otherwise. We set
\begin{equation}\label{eq:N2}
\begin{aligned}
&\sigma_i \coloneqq \tfrac{1}{2}(\alpha_i + \beta_i + \gamma_i + \delta_i),
\qquad
&&\overline{\alpha}_i \coloneqq \sigma_i - \alpha_i,
\qquad
&&\overline{\beta}_i \coloneqq \sigma_i - \beta_i,
\\
&\varepsilon_i\coloneqq \sign{(\sin{\sigma_i})},
\qquad
&&\overline{\gamma}_i \coloneqq \sigma_i - \gamma_i,
\qquad
&&\overline{\delta}_i \coloneqq \sigma_i - \delta_i,
\end{aligned}
\end{equation}
and
\begin{equation}\label{eq:N3}
\begin{aligned}
&a_i \coloneqq \frac{\sin \alpha_i}{\sin \overline{\alpha}_i},
\qquad
&&b_i \coloneqq \frac{\sin \beta_i}{\sin \overline{\beta}_i},
\qquad
&&c_i \coloneqq \frac{\sin \gamma_i}{\sin \overline{\gamma}_i},
\qquad
&&d_i \coloneqq \frac{\sin \delta_i}{\sin \overline{\delta}_i},
\\
&M_i \coloneqq a_i b_i  c_i  d_i,
\qquad
&&r_i \coloneqq a_i d_i,
\qquad
&&s_i \coloneqq c_i d_i,
\qquad
&&f_i \coloneqq a_i c_i,
\\
&u_i \coloneqq 1 - M_i,
\qquad
&&x_i \coloneqq \frac{1}{r_i - 1},
\qquad
&&y_i \coloneqq \frac{1}{s_i - 1},
\qquad
&&z_i \coloneqq \frac{1}{f_i - 1}.
    \end{aligned}
\end{equation}
These quantities are well-defined, and $M_i$, $r_i$, $s_i$, $f_i\neq 0, 1$, see Lemma~2 of~\cite{NURMATOV2026104102}.

A polyhedron has \emph{equimodular elliptic type} 
if it has elliptic type and the following three conditions hold:
\begin{itemize}
    \item[(1)] The vertices have equal \emph{moduli}:
    \begin{equation}\label{eq:N4} 
        M_1 = M_2 = M_3 = M_4 \eqqcolon M. 
    \end{equation}
   \item[(2)] The amplitudes at common edges match:
\begin{equation}\label{eq:N5}
   r_1 = r_2,
   \quad
   s_1 = s_4,
   \quad
   r_3 = r_4,
   \quad
   s_2 = s_3.
\end{equation}
    \item[(3)] A signed sum of phase shifts is a \emph{period}: 
    for some $e_1$, $e_2$, $e_3$, $e_4 \in \{\pm 1\}$,
    \begin{equation}\label{eq:N6} 
        \sum_{i=1}^4 e_it_i \in \Lambda.
    \end{equation}
\end{itemize}

In the latter equation, we use the following notation. 
First, we introduce the \emph{elliptic modulus} $k$ 
and the \emph{complementary modulus} $k'$
\begin{equation}\label{eq:N7}
        k\coloneqq
    \begin{cases}
        \sqrt{1-M}, & \text{if } M < 1, \\
        \sqrt{1-M^{-1}}, & \text{if } M > 1,
    \end{cases}
    \qquad
    k' \coloneqq \sqrt{1 - k^2}.
    \end{equation}
For elliptic polyhedra, $M > 0$ by Lemma~5 of~\cite{NURMATOV2026104102}, see Remark~\ref{rem:existence}, 
and $k$, $k' \in (0, 1)$. The \emph{basis quarter periods} are
\begin{equation}\label{eq:N8}
\begin{aligned}
    &K\coloneqq K(k) = \int_{0}^{1} \frac{dx}{\sqrt{(1-x^2)(1-k^2x^2)}}, 
    \\
    &K'\coloneqq K(k') = \int_{0}^{1} \frac{dx}{\sqrt{(1-x^2)(1-k'^2x^2)}}.
\end{aligned}
\end{equation}
Then we introduce the lattice
\begin{equation}\label{eq:N9}
    \Lambda \coloneqq
    \begin{cases}
        4K\mathbb{Z} + 2\mathbf{i}K'\mathbb{Z}, & \text{if } M < 1, 
        \\
        4K\mathbb{Z} + \left(2K + 2\mathbf{i}K'\right)\mathbb{Z}, & \text{if } M > 1,
    \end{cases}
\end{equation}
where $\mathbf{i}$ (in bold font) is the imaginary unit. 

The \emph{amplitudes} $p_i$ and $q_i$ are given by
\begin{equation}\label{eq:N10}
    p_i\coloneqq \sqrt{r_i - 1} \in \mathbb{R}_{>0} \cup \mathbf{i} \mathbb{R}_{>0},
\qquad
q_i\coloneqq \sqrt{s_i - 1} \in \mathbb{R}_{>0} \cup \mathbf{i} \mathbb{R}_{>0}.
\end{equation}

Finally, the \emph{phase shift} $t_i$ satisfies
\begin{equation}\label{eq:N11} 
\operatorname{dn}{t_i} \coloneqq\operatorname{dn}{(t_i, k)} =
    \begin{cases}
        \sqrt{f_i}, & \text{if } M < 1, \\
        \dfrac{1}{\sqrt{f_i}}, & \text{if } M > 1.
    \end{cases}
\end{equation}
Here, \eqref{eq:N11} alone does not determine $t_i$, because
$\operatorname{dn}{(t_i + 2K)} = \operatorname{dn}{t_i}$;
Table~\ref{tab:tpqsigma} gives the interval where $t_i$ lies, depending
on the values of $p_iq_i$ and $\sigma_i$, and on this interval
\eqref{eq:N11} determines $t_i$ uniquely: 
$\operatorname{dn}$ maps the intervals for $p_iq_i \in \mathbb{R}$
bijectively onto $(1, +\infty)$ and those for
$p_iq_i \in \mathbf{i}\mathbb{R}_{>0}$ onto $(0, k')$,
see~\cite[Chapter~22]{WhittakerWatson2021},
and the value in~\eqref{eq:N11} lies in the corresponding range by Lemma~4
of~\cite{NURMATOV2026104102}, \eqref{eq:N7}, and~\eqref{eq:N10}.

\begin{table}[htbp!]
\centering
\footnotesize
\caption{Phase shift intervals for different $p_iq_i$ and $\sigma_i$.}
\label{tab:tpqsigma}
\renewcommand{\arraystretch}{1.2}
\begin{tabular}{c|c|c|c}
& $p_iq_i \in \mathbb{R}_{>0}$ & $p_iq_i \in \mathbf{i} \mathbb{R}_{>0}$ & $p_iq_i \in \mathbb{R}_{<0}$ \\
\hline
$\sigma_i < \pi$ & $(0, \mathbf{i}K')$ & $(K, K + \mathbf{i}K')$ & $(2K, 2K + \mathbf{i}K')$ \\
\hline
$\sigma_i > \pi$ & $(2K, 2K + \mathbf{i}K')$ & $(3K, 3K + \mathbf{i}K')$ & $(0, \mathbf{i}K')$
\end{tabular}
\end{table}

Table~\ref{tab:sncn} expresses
$\operatorname{sn}{t_i}\coloneqq\operatorname{sn}{(t_i, k)}$ 
and $\operatorname{cn}{t_i}\coloneqq\operatorname{cn}{(t_i, k)}$ 
as functions of $\operatorname{dn}{t_i}$ on each of the four
intervals of Table~\ref{tab:tpqsigma}, where $k \in (0, 1)$, by the
standard evaluations of the Jacobi elliptic functions,
see~\cite[Chapter~22]{WhittakerWatson2021};
all radicands in Table~\ref{tab:sncn} are strictly positive, because
$\operatorname{dn}{t_i} > 1$ if $p_iq_i \in \mathbb{R}$ and
$0 < \operatorname{dn}{t_i} < k'$ if $p_iq_i \in \mathbf{i}\mathbb{R}_{>0}$.

\begin{table}[htbp!]
\centering
\footnotesize
\caption{Expressions for $\operatorname{sn}{t_i}$ and $\operatorname{cn}{t_i}$ across phase shift intervals for $k\in(0,1)$.}
\label{tab:sncn}
\renewcommand{\arraystretch}{2.0}
\begin{tabular}{c|c|c|c|c}
${t_i}$ & $(0, \mathbf{i}K')$ & $(K, K + \mathbf{i}K') $ & $(2K, 2K + \mathbf{i}K')$ & $ (3K, 3K + \mathbf{i}K') $\\
\hline
$\operatorname{sn}{t_i}$  & $\mathbf{i} \, \frac{\sqrt{\operatorname{dn}^2{t_i} - 1}}{k}$ & $\frac{\sqrt{1- \operatorname{dn}^2{t_i}}}{k}$ & $-\mathbf{i} \, \frac{\sqrt{\operatorname{dn}^2{t_i} - 1}}{k}$ & $-\frac{\sqrt{1- \operatorname{dn}^2{t_i}}}{k}$ \\
\hline
$\operatorname{cn}{t_i}$   & $\frac{\sqrt{\operatorname{dn}^2{t_i} - 1 + k^2}}{k}$ &  $- \mathbf{i} \, \frac{\sqrt{1 - k^2 -\operatorname{dn}^2{t_i}}}{k}$ & $-\frac{\sqrt{\operatorname{dn}^2{t_i} - 1 + k^2}}{k}$ & $\mathbf{i} \, \frac{\sqrt{1 - k^2 -\operatorname{dn}^2{t_i}}}{k}$
\end{tabular}
\end{table}

Next, we set
\begin{equation}\label{eq:N12}
     \chi \coloneqq \varepsilon_1 \varepsilon_2 \varepsilon_3 \varepsilon_4\, \sign{(z_1z_3)}. 
\end{equation}

\begin{definition}\label{def:gqs}
A polyhedron is \emph{generalized quasi-symmetric}, or a
\emph{GQS-net}, if its flat angles satisfy
conditions~\eqref{eq:N1}, \eqref{eq:N4}, and~\eqref{eq:N5}, and
one of the following four conditions holds for some distinct
$i$, $j$, $l$, $m \in \{1, 2, 3, 4\}$:
\begin{itemize}
    \item[(i)] $f_i = f_j$, $f_l = f_m$, and
    $\chi\,(f_i -1) (f_l - 1) > 0$;
    \item[(ii)] $f_i f_j = f_l f_m = M$;
    \item[(iii)] $f_i + f_j - f_i f_j
    = f_l + f_m - f_l f_m = M$ and $\chi\,(1 - M) > 0$;
    \item[(iv)] $\dfrac{f_i + f_j - 1}{f_i f_j}
    = \dfrac{f_l + f_m - 1}{f_l f_m} = \dfrac{1}{M}$.
\end{itemize}
\end{definition}

\begin{definition}[see Section~2.2 and Lemma~4.3 of~\cite{izmestiev-2017}]\label{def:bricard}
\emph{Bricard's system of equations} is the following system of
polynomial equations in the variables $Z$, $W_2$, $U$, $W_1$:
\begin{equation}\label{eq:N13}
    \begin{aligned}
        &P_1(Z, W_1) = 0, \qquad &&P_2(Z, W_2) = 0, \\
        &P_3(U, W_2) = 0, \qquad &&P_4(U, W_1) = 0,
    \end{aligned}
\end{equation}
where
$P_i(X, Y)
= c_{22i}\,X^2 Y^2 + c_{20i}\,X^2 + c_{02i}\,Y^2
+ 2 c_{11i}\,X Y + c_{00i}$, with the coefficients 
$c_{00i} = \sin{\sigma_i} \sin{\overline{\beta}_i}$ and
\begin{equation*}
\begin{aligned}
    &c_{22i} = \sin{\overline{\delta}_i}
    \sin{(\overline{\delta}_i - \beta_i)},
    \qquad
    &&c_{20i} = \sin{\overline{\alpha}_i}
    \sin{(\overline{\alpha}_i - \beta_i)},
    \\
    &c_{02i} = \sin{\overline{\gamma}_i}
    \sin{(\overline{\gamma}_i - \beta_i)}, 
    \qquad
    &&c_{11i} = -\sin{\alpha_i} \sin{\gamma_i}.
\end{aligned}
\end{equation*}
\end{definition}
For a polyhedron with central dihedral angles $\theta_1, \ldots,
\theta_4 \neq 0$, the values $Z = \cot{\frac{\theta_1}{2}}$,
$W_2 = \cot{\frac{\theta_2}{2}}$, $U = \cot{\frac{\theta_3}{2}}$,
and $W_1 = \cot{\frac{\theta_4}{2}}$ satisfy Bricard's system, see
Figure~\ref{fig:03} and Lemma~\ref{lem:bricard}. 

Note that our dihedral angles are $\pi$ minus the dihedral angles
introduced in~\cite{izmestiev-2017}, see~\cite{ErofeevIvanov2020}.
Condition~\eqref{eq:N5} corrects a misprint in Section~3.3.1
of~\cite{izmestiev-2017}, see~\cite[Section~4.3]{HePhD2021}
and~\cite[Section~2]{NURMATOV2026104102}.
System~\eqref{eq:N13} also corrects a misprint in
\cite[Eq.~(4)]{izmestiev-2017}, where, in the notation there, $w_1$ and
$w_2$ are interchanged in the equations $P_3 = 0$ and $P_4 = 0$;
see~\cite[Section~6.2]{izmestiev-2017}.

Finally, we introduce the notation used in the parametrization
of Theorem~\ref{th:main}. We define
\begin{equation}\label{eq:N14}
    u \coloneqq  1 - M
\end{equation}
and, for $x$, $t \in \mathbb{R}$,
\begin{equation}\label{eq:N15}
    D(x, t) \coloneqq \bigl(x t^2 - 1\bigr)
    \bigl(1 - u x t^2\bigr),
\end{equation}
together with the \emph{admissible set}
\begin{equation}\label{eq:N16}
    I(x) \coloneqq \bigl\{\, t \in \mathbb{R} :
    x D(x, t) > 0 \,\bigr\}.
\end{equation}
In particular, if $M > 0$ and $M \neq 1$, that is, $u < 1$ and $u \neq 0$,
then for every $x \in \mathbb{R} \setminus \{0\}$ the admissible set is
\begin{equation}\label{eq:N17}
    I(x) =
    \begin{cases}
    \Bigl(-\tfrac{1}{\sqrt{u x}},\,
    -\tfrac{1}{\sqrt{x}}\Bigr)
    \cup
    \Bigl(\tfrac{1}{\sqrt{x}},\,
    \tfrac{1}{\sqrt{u x}}\Bigr),
    & \text{if } x > 0 \text{ and } u > 0,
    \\[2mm]
    (-\infty,\, +\infty),
    & \text{if } x < 0 \text{ and } u > 0,
    \\[2mm]
    \Bigl(-\infty,\, -\tfrac{1}{\sqrt{x}}\Bigr)
    \cup
    \Bigl(\tfrac{1}{\sqrt{x}},\, +\infty\Bigr),
    & \text{if } x > 0 \text{ and } u < 0,
    \\[2mm]
    \Bigl(-\tfrac{1}{\sqrt{u x}},\,
    \tfrac{1}{\sqrt{u x}}\Bigr),
    & \text{if } x < 0 \text{ and } u < 0.
    \end{cases}
\end{equation}
In all four cases, $I(x)$ is a union of at most two nonempty open
intervals; in the first case,
$\frac{1}{\sqrt{x}} < \frac{1}{\sqrt{u x}}$ because $0 < u < 1$ there.

\begin{remark}\label{rem:existence}
The ``if'' part of Lemma~5 of~\cite{NURMATOV2026104102} does not hold in
general. For example, let
$\alpha_1 = \beta_1 = \gamma_1 = \pi/3$, $\delta_1 = 2\pi/3$,
$\alpha_2 = \beta_2 = \gamma_2 = \pi/6$, $\delta_2 = \pi/3$,
$\alpha_3 = \gamma_3 = \pi/6$, $\beta_3 = 2\pi/3$, $\delta_3 = \pi/2$,
$\alpha_4 = \beta_4 = \gamma_4 = \pi/3$, and $\delta_4 = \pi/2$. Then
$\delta_1 + \delta_2 + \delta_3 + \delta_4 = 2\pi$ and, for each
$i = 1, \ldots, 4$, $\alpha_i$, $\beta_i$, $\gamma_i$,
$\delta_i \in (0, \pi)$ satisfy condition~\eqref{eq:N1} and $M_i > 0$
by~\eqref{eq:N3}, because $\overline{\alpha}_i$, $\overline{\beta}_i$,
$\overline{\gamma}_i$, $\overline{\delta}_i \in (0, \pi)$. Nevertheless, no
polyhedron with these flat angles exists. Indeed, otherwise,
by the spherical law of cosines,
$\cos{\angle A_1A_2C_2} = \cos{\delta_2}\cos{\gamma_2}
+ \sin{\delta_2}\sin{\gamma_2}\cos{\theta_2}$ and
$\cos{\angle A_4A_3C_3} = \cos{\delta_3}\cos{\gamma_3}
+ \sin{\delta_3}\sin{\gamma_3}\cos{\theta_2}$, and the triangle
inequality for the angles between rays gives
$\angle A_1A_2C_2 \leq \alpha_2 + \beta_2 = \pi/3$ and
$\angle A_4A_3C_3 \geq \beta_3 - \alpha_3 = \pi/2$. Hence
$\cos{\theta_2} \geq 2/\sqrt{3} - 1$ and $\cos{\theta_2} \leq 0$,
a contradiction. Only the ``only if'' part of Lemma~5
of~\cite{NURMATOV2026104102} is used in this paper, and the other
results of~\cite{NURMATOV2026104102} depend only on that part as well.
\end{remark}

\section{Main results}
\label{sec:results}
Henceforth, we use the notation established
in~\eqref{eq:N2}--\eqref{eq:N17}. The
proofs of the main results below are given in Section~\ref{sec:proofs}.
The main results also
include Proposition~\ref{prop:signs}, which determines the signs
$e_1, \ldots, e_4$ satisfying condition~\eqref{eq:N6}; it is stated in
Section~\ref{sec:proofs}, because it uses notation introduced there.
Finally, Proposition~\ref{prop:ex:flex} assembles GQS-nets into flexible
$m \times n$ nets for all $m, n \geq 3$, with flexions in closed form; it is
stated and proved in Section~\ref{sec:examples}.

\begin{proposition}\label{prop:eqMrs}
Assume that $\alpha_i$, $\beta_i$, $\gamma_i$,
$\delta_i \in (0, \pi)$ satisfy condition~\eqref{eq:N1} for each
$i = 1, \ldots, 4$. Then
conditions~\eqref{eq:N4} and~\eqref{eq:N5} hold if and only if
\begin{equation*}
\begin{aligned}
\frac{\sin{\alpha_1}\sin{\delta_1}}{\sin{\beta_1}\sin{\gamma_1}} &= \frac{\sin{\alpha_2}\sin{\delta_2}}{\sin{\beta_2}\sin{\gamma_2}},
&\qquad
\frac{\sin{\overline{\alpha}_1}\sin{\overline{\delta}_1}}{\sin{\overline{\beta}_1}\sin{\overline{\gamma}_1}} &= \frac{\sin{\overline{\alpha}_2}\sin{\overline{\delta}_2}}{\sin{\overline{\beta}_2}\sin{\overline{\gamma}_2}}, 
\\
\frac{\sin{\gamma_2}\sin{\delta_2}}{\sin{\alpha_2}\sin{\beta_2}} &= \frac{\sin{\gamma_3}\sin{\delta_3}}{\sin{\alpha_3}\sin{\beta_3}},
&\qquad
\frac{\sin{\overline{\gamma}_2}\sin{\overline{\delta}_2}}{\sin{\overline{\alpha}_2}\sin{\overline{\beta}_2}} &= \frac{\sin{\overline{\gamma}_3}\sin{\overline{\delta}_3}}{\sin{\overline{\alpha}_3}\sin{\overline{\beta}_3}},
\\
\frac{\sin{\alpha_3}\sin{\delta_3}}{\sin{\beta_3}\sin{\gamma_3}} &= \frac{\sin{\alpha_4}\sin{\delta_4}}{\sin{\beta_4}\sin{\gamma_4}},
&\qquad
\frac{\sin{\overline{\alpha}_3}\sin{\overline{\delta}_3}}{\sin{\overline{\beta}_3}\sin{\overline{\gamma}_3}} &= \frac{\sin{\overline{\alpha}_4}\sin{\overline{\delta}_4}}{\sin{\overline{\beta}_4}\sin{\overline{\gamma}_4}}, 
\\
\frac{\sin{\gamma_4}\sin{\delta_4}}{\sin{\alpha_4}\sin{\beta_4}} &= \frac{\sin{\gamma_1}\sin{\delta_1}}{\sin{\alpha_1}\sin{\beta_1}},
&\qquad
\frac{\sin{\overline{\gamma}_4}\sin{\overline{\delta}_4}}{\sin{\overline{\alpha}_4}\sin{\overline{\beta}_4}} &= \frac{\sin{\overline{\gamma}_1}\sin{\overline{\delta}_1}}{\sin{\overline{\alpha}_1}\sin{\overline{\beta}_1}}.
\end{aligned}
\end{equation*}
\end{proposition}

\begin{proposition}\label{prop:shifts}
Assume that $\alpha_i$, $\beta_i$, $\gamma_i$, $\delta_i \in (0, \pi)$
satisfy conditions~\eqref{eq:N1}, \eqref{eq:N4}, and~\eqref{eq:N5} for
each $i = 1, \ldots, 4$, and $M > 0$.
Then condition~\eqref{eq:N6} holds if and only if, for some
$\zeta_1$, $\zeta_2 \in \{\pm 1\}$,
\begin{equation*}
\zeta_1 \, M \sqrt{\prod_{i = 1}^{4} (f_i - 1)} 
+ \zeta_2 \, (1 - M) \sqrt{\prod_{i = 1}^{4} f_i}  
-  \zeta_2 \, \chi \sqrt{\prod_{i = 1}^{4} (f_i - M)}
= M(1 - M).
\end{equation*}
\end{proposition}

\begin{theorem}\label{th:gqs}
Every GQS-net has equimodular elliptic type.
\end{theorem}

\begin{theorem}\label{th:main}
Assume that $\alpha_i$, $\beta_i$, $\gamma_i$,
$\delta_i \in (0, \pi)$ satisfy conditions~\eqref{eq:N1},
\eqref{eq:N4}, and~\eqref{eq:N5} for each
$i = 1, \ldots, 4$, $\delta_1 + \delta_2 + \delta_3 + \delta_4
= 2\pi$, and $M > 0$.
Assume also that $e_1$, $e_2$,
$e_3$, $e_4 \in \{\pm 1\}$ satisfy condition~\eqref{eq:N6}.
Then $\alpha_i$, $\beta_i$, $\gamma_i$, $\delta_i$ are 
flat angles of some flexible polyhedron with a flexion given by
\begin{equation*}
\begin{aligned}
    &\cot{\frac{\theta_1(t)}{2}} = t,
    \\[1mm]
    &\cot{\frac{\theta_2(t)}{2}}
    = \frac{\varepsilon_2 \Bigl(
    t\sqrt{u x_1 y_2(1 + z_2)(1 + u z_2)}
    + e_0\,e_2 \sqrt{u y_2 z_2 D(x_1, t)}
    \Bigr)}
    {y_2 u (z_2 + x_1 t^2)},
    \\[1mm]
    &\cot{\frac{\theta_4(t)}{2}}
    = \frac{\varepsilon_1 \Bigl(
    t\sqrt{u x_1 y_1(1 + z_1)(1 + u z_1)}
    - e_0\,e_1 \sqrt{u y_1 z_1 D(x_1, t)}
    \Bigr)}
    {y_1 u (z_1 + x_1 t^2)},
    \\[1mm]
   &\cot{\frac{\theta_3(t)}{2}}
    = 
    \varepsilon_2 \varepsilon_3 \sign{(z_2)}\sqrt{\frac{z_2 z_3}{x_1 x_3}}
    \cdot
    \frac{|x_1| H_{2} H_{3}\,t
    + e_0\, e_2 H_{1} \sqrt{x_1D(x_1, t)}}
    {z_2z_3 + u x_1 H_{1}^2\,t^2},
\end{aligned}
\end{equation*}
where each of the two signs $e_0 \in \{\pm 1\}$ gives a flexion, 
$t$ ranges over the admissible set
$I(x_1)$, excluding the finitely many $t$ at which any of the
denominators of the applicable formulas vanishes (the resulting
set has nonempty interior), 
and $H_{1}$, $H_{2}$, $H_{3}$ are determined as follows:
\begin{itemize}
    \item[(a)] if $uz_2z_3 \neq 1$, then they
    are given by 
\begin{equation*}
\begin{aligned}
    H_{1} &= \frac{|z_2| \sqrt{z_3 (1 + z_3) (1 + uz_3)} + e_2\, e_3 |z_3| \sqrt{z_2 (1 + z_2) (1 + uz_2)}}{uz_2z_3 - 1},
    \\
    H_{2} &= \frac{|z_2z_3| \sqrt{(1+uz_2)(1 + uz_3)} + e_2\,e_3\sqrt{z_2z_3(1+z_2)(1+z_3)}}{uz_2z_3-1},
    \\
    H_{3} &= \frac{u\sqrt{z_2z_3(1+z_2)(1+z_3)}\,\sign{(z_2z_3)} + e_2\,e_3\sqrt{(1+uz_2)(1+uz_3)}}{uz_2z_3-1};
\end{aligned}
\end{equation*}
    \item[(b)] if $uz_2z_3=1$ and
    $e_2 = -e_3$, then they are given by 
\begin{equation*}
\begin{aligned}
    H_{1} &= \frac{z_2z_3(z_3-z_2)}{2 \sqrt{z_2z_3(1+z_2)(1+z_3)}},
    \\
    H_{2} &= \frac{z_2z_3(z_2+z_3+2)}{2 \sqrt{z_2z_3(1+z_2)(1+z_3)}},
    \\
    H_{3} &= \frac{(z_2+z_3+ 2z_2z_3)\sign{(z_2z_3)}}{2 \sqrt{z_2z_3(1+z_2)(1+z_3)}};
\end{aligned}
\end{equation*}
    \item[(c)] if $uz_2z_3=1$ and
    $e_2 = e_3$, then the formula for $\cot{\frac{\theta_3(t)}{2}}$
    is replaced by
\begin{equation*}
    \cot{\frac{\theta_3(t)}{2}}
    = 
    \varepsilon_2 \varepsilon_3 \sqrt{\frac{z_2 z_3}{x_1 x_3}}
    \cdot
    \frac{\sign{(y_2)}}
    {t}.
\end{equation*}
\end{itemize}
Moreover, for every such $t$,
\begin{equation*}
e_1\cot\frac{\theta_{41'}(t)}{2}
    = -e_2\cot\frac{\theta_{22'}(t)}{2}
    = e_0\sqrt{\frac{(1 + ux_1)(x_1t^2 - 1)}{(1 + x_1)(1 - ux_1t^2)}}.
\end{equation*}
\end{theorem}

Under the hypotheses of Proposition~\ref{prop:eqMrs}, no
denominator in the proposition vanishes; under the hypotheses of
Proposition~\ref{prop:shifts}, all the radicands there are
strictly positive; under the hypotheses of Theorem~\ref{th:main},
all the radicands in the theorem are strictly positive for every
$t \in I(x_1)$, see~\eqref{eq:N16}. All three
facts are shown in Section~\ref{sec:proofs}.

Every elliptic QS-net, that is, every QS-net
of~\cite{NURMATOV2026104102} that has elliptic type, is a
GQS-net. Indeed, as established in the proof of Theorem~1
of~\cite{NURMATOV2026104102}, the flat angles of such a
polyhedron satisfy conditions~\eqref{eq:N4} and~\eqref{eq:N5},
together with $f_1 = f_4$ and $f_1 f_2 = M$. Moreover, by
Definition~1 of~\cite{NURMATOV2026104102}, we have
$(\alpha_3, \beta_3, \gamma_3, \delta_3)
= (\pi - \alpha_2, \pi - \beta_2, \pi - \gamma_2,
\pi - \delta_2)$; since $f_i$ depends only on the flat angles at
the $i$-th vertex, it follows by Lemma~3
of~\cite{NURMATOV2026104102} that $f_2 = f_3$. Hence
$f_1 f_2 = f_4 f_3 = M$, which is condition~(ii) of
Definition~\ref{def:gqs} with $(i, j, l, m) = (1, 2, 4, 3)$. 

The dihedral angles in Theorem~\ref{th:main} are considered in the full
range $[-\pi, \pi)$ excluding $0$. 
The two flexions of Theorem~\ref{th:main} are mirror images of each other:
its formulas give, with the parameter~$-t$ and the sign~$e_0$, the
negatives of the values of $\cot\frac{\theta_1}{2}, \ldots,
\cot\frac{\theta_4}{2}$ with the parameter~$t$ and the sign~$-e_0$,
because $D(x_1, t)$ is even in~$t$ and $H_1$, $H_2$, $H_3$ depend on
neither~$t$ nor~$e_0$.
Applied after the cyclic permutation
of the indices by $j = 1, 2, 3$, the theorem takes
$\cot\frac{\theta_{j+1}}{2}$ as the parameter, ranging over $I(y_2)$,
$I(x_3)$, and $I(y_4)$, respectively, because for odd $j$ the roles of
$x_i$ and $y_i$ are interchanged, see~\eqref{eq:N3}. In this form, it
expresses the other six corner dihedral angles through $\theta_2$,
$\theta_3$, and $\theta_4$, respectively. In particular, we get
$e_2\cot\frac{\theta_{12'}}{2} = -e_3\cot\frac{\theta_{33'}}{2}$,
$e_3\cot\frac{\theta_{23'}}{2} = -e_4\cot\frac{\theta_{44'}}{2}$, and
$e_4\cot\frac{\theta_{34'}}{2} = -e_1\cot\frac{\theta_{11'}}{2}$
for the polyhedra of the flexions of Theorem~\ref{th:main} with all but
finitely many values of the parameter $t$, see
Remark~\ref{rem:th:main:2}; together with the relation
$e_1\cot\frac{\theta_{41'}}{2} = -e_2\cot\frac{\theta_{22'}}{2}$ of
Theorem~\ref{th:main}, we have
$\theta_{41'} \equiv -e_1e_2\,\theta_{22'}$,
$\theta_{12'} \equiv -e_2e_3\,\theta_{33'}$,
$\theta_{23'} \equiv -e_3e_4\,\theta_{44'}$, and
$\theta_{34'} \equiv -e_4e_1\,\theta_{11'} \pmod{2\pi}$,
which gives the signs $e_1, \ldots, e_4$ a geometric meaning.

Theorem~\ref{th:main} asserts the existence of a flexible
polyhedron with the flat angles $\alpha_i$, $\beta_i$,
$\gamma_i$, $\delta_i$, not the flexibility of every polyhedron
with these flat angles: by Proposition~3 of~\cite{NURMATOV2026104102},
two polyhedra may have the same flat angles, which satisfy the
hypotheses of the theorem, while only one of them is flexible. 
By Remark~\ref{rem:lem:main:2}, all solutions of Bricard's system
except finitely many lie on the one-parameter families of solutions
given by Lemma~\ref{lem:main:2}, and every flexion of a polyhedron with
these flat angles runs within these families while its central
dihedral angles $\theta_1, \ldots, \theta_4$ are nonzero. In particular,
a polyhedron for which
$\bigl(\cot{\frac{\theta_1}{2}}, \ldots, \cot{\frac{\theta_4}{2}}\bigr)$
lies on none of these families is rigid; an example is the polyhedron of
Proposition~3 of~\cite{NURMATOV2026104102} that is not flexible.

\section{Proofs}
\label{sec:proofs}
The proof of Proposition~\ref{prop:eqMrs} relies on
Lemma~\ref{lem:eqMrs}, and the proof of Proposition~\ref{prop:shifts}
on Lemma~\ref{lem:shifts:3}, which in turn relies on
Lemmas~\ref{lem:shifts:1} and~\ref{lem:shifts:2}. The proof of
Theorem~\ref{th:gqs} relies on Proposition~\ref{prop:shifts}, the proof
of Theorem~\ref{th:main} on Lemmas~\ref{lem:shifts:2}
and~\ref{lem:main:1}--\ref{lem:bricard} and on
Remarks~\ref{rem:lem:gqs} and~\ref{rem:lem:main:6:2}, and the proof of
Proposition~\ref{prop:signs} on Lemmas~\ref{lem:main:3},
\ref{lem:main:6}, and~\ref{lem:signs} and on
Remark~\ref{rem:lem:main:6:1}.

\begin{lemma}\label{lem:eqMrs}
Assume that $\alpha_i$, $\beta_i$, $\gamma_i$,
$\delta_i \in (0, \pi)$ satisfy condition~\eqref{eq:N1} for each
$i = 1, \ldots, 4$, and let $l$, $m \in \{1, \ldots, 4\}$ be
distinct. Then:
\begin{itemize}
    \item[(i)] $r_l = r_m$ and $M_l = M_m$ if and only if
    \begin{equation*}
    \frac{\sin{\alpha_l}\sin{\delta_l}}{\sin{\beta_l}\sin{\gamma_l}} = \frac{\sin{\alpha_m}\sin{\delta_m}}{\sin{\beta_m}\sin{\gamma_m}},
    \qquad
    \frac{\sin{\overline{\alpha}_l}\sin{\overline{\delta}_l}}{\sin{\overline{\beta}_l}\sin{\overline{\gamma}_l}} = \frac{\sin{\overline{\alpha}_m}\sin{\overline{\delta}_m}}{\sin{\overline{\beta}_m}\sin{\overline{\gamma}_m}};
    \end{equation*}
    \item[(ii)] $s_l = s_m$ and $M_l = M_m$ if and only if
    \begin{equation*}
    \frac{\sin{\gamma_l}\sin{\delta_l}}{\sin{\alpha_l}\sin{\beta_l}} = \frac{\sin{\gamma_m}\sin{\delta_m}}{\sin{\alpha_m}\sin{\beta_m}},
    \qquad
    \frac{\sin{\overline{\gamma}_l}\sin{\overline{\delta}_l}}{\sin{\overline{\alpha}_l}\sin{\overline{\beta}_l}} = \frac{\sin{\overline{\gamma}_m}\sin{\overline{\delta}_m}}{\sin{\overline{\alpha}_m}\sin{\overline{\beta}_m}}.
    \end{equation*}
\end{itemize}
\end{lemma}

\begin{proof}
By~\eqref{eq:N1}, $\overline{\alpha}_i$,
$\overline{\beta}_i$, $\overline{\gamma}_i$,
$\overline{\delta}_i \notin \pi\mathbb{Z}$, see~\eqref{eq:N2},
and, by Lemma~2 of~\cite{NURMATOV2026104102}, $r_i$,
$s_i$, $M_i$, $M_i/r_i$, $M_i/s_i \neq 0, 1$ are well-defined reals.
Since, in addition, $\alpha_i$, $\beta_i$, $\gamma_i$,
$\delta_i \in (0, \pi)$, it follows that both sides of the
identities below are well-defined. For each $i = 1, \ldots, 4$,
we have
\begin{equation*}
\begin{aligned}
    &\Phi_i \coloneqq \frac{\sin{\alpha_i}\sin{\delta_i}}{\sin{\beta_i}\sin{\gamma_i}}
    = \frac{r_i(r_i - M_i)}{M_i(r_i - 1)},
    \qquad
    &&\overline{\Phi}_i \coloneqq \frac{\sin{\overline{\alpha}_i}\sin{\overline{\delta}_i}}{\sin{\overline{\beta}_i}\sin{\overline{\gamma}_i}}
    = \frac{r_i - M_i}{r_i(r_i - 1)}, \\
    &\Psi_i \coloneqq \frac{\sin{\gamma_i}\sin{\delta_i}}{\sin{\alpha_i}\sin{\beta_i}}
    = \frac{s_i(s_i - M_i)}{M_i(s_i - 1)},
    \qquad
    &&\overline{\Psi}_i \coloneqq \frac{\sin{\overline{\gamma}_i}\sin{\overline{\delta}_i}}{\sin{\overline{\alpha}_i}\sin{\overline{\beta}_i}}
    = \frac{s_i - M_i}{s_i(s_i - 1)},
\end{aligned}
\end{equation*}
as verified in Section~1 of the supplement
(\texttt{helper.nb})~\cite{Nurmatov2026Repo}.
Observe that $\Phi_i$, $\Psi_i > 0$ 
because $\alpha_i$, $\beta_i$, $\gamma_i$, $\delta_i \in (0, \pi)$. 
Moreover, $r_i^2$, $s_i^2\neq M_i$; otherwise $\Phi_i$, $\Psi_i = -1$,
contradicting $\Phi_i$, $\Psi_i > 0$.
Since $r_i$,
$s_i$, $r_i^2$, $s_i^2\neq M_i$, it follows that $\overline{\Phi}_i$,
$\overline{\Psi}_i \neq -1, 0$. In this notation, the equations in
(i) and (ii) are $\Phi_l = \Phi_m$, $\overline{\Phi}_l = \overline{\Phi}_m$ 
and $\Psi_l = \Psi_m$, $\overline{\Psi}_l = \overline{\Psi}_m$, respectively.
We prove (i); statement (ii) follows by replacing 
$(r_i, \Phi_i, \overline{\Phi}_i)$ with $(s_i, \Psi_i, \overline{\Psi}_i)$.

($\Rightarrow$) If $r_l = r_m$ and $M_l = M_m$, then $\Phi_l = \Phi_m$ and
$\overline{\Phi}_l = \overline{\Phi}_m$ by
the identities for $\Phi_i$ and $\overline{\Phi}_i$.

($\Leftarrow$) Because $\Phi_i > 0$ and $\overline{\Phi}_i \neq -1, 0$,
the identities for $\Phi_i$ and $\overline{\Phi}_i$ give
$r_i = \frac{\Phi_i\,(1 + \overline{\Phi}_i)}{\overline{\Phi}_i\,(1 + \Phi_i)}$ 
and $M_i = \frac{r_i^2\,\overline{\Phi}_i}{\Phi_i}$.
Thus $\Phi_l = \Phi_m$ and
$\overline{\Phi}_l = \overline{\Phi}_m$ imply $r_l = r_m$ and
$M_l = M_m$.
\end{proof}

\begin{proof}[Proof of Proposition~\ref{prop:eqMrs}]
By~(i) and~(ii) of Lemma~\ref{lem:eqMrs}, applied with
$(l, m)$ equal to $(1, 2)$ and $(3, 4)$, and to $(2, 3)$ and
$(4, 1)$, respectively, the eight stated equations hold if and
only if $(r_1, M_1) = (r_2, M_2)$, $(s_2, M_2) = (s_3, M_3)$, 
$(r_3, M_3) = (r_4, M_4)$, and $(s_4, M_4) = (s_1, M_1)$, 
which are conditions~\eqref{eq:N4} and~\eqref{eq:N5}.
\end{proof}

The equalities appearing in conditions (i)--(iv) of
Definition~\ref{def:gqs} are relations between the quantities
$f_1, \ldots, f_4$ and $M$. Lemma~\ref{lem:gqs} below
characterizes each of them by two equations between sine ratios,
one involving only the angles $\alpha_i, \ldots, \delta_i$ and
the other only the angles
$\overline{\alpha}_i, \ldots, \overline{\delta}_i$.

\begin{lemma}\label{lem:gqs}
Assume that $\alpha_i$, $\beta_i$, $\gamma_i$, $\delta_i \in (0, \pi)$
satisfy conditions~\eqref{eq:N1} and~\eqref{eq:N4} for each
$i = 1, \ldots, 4$, and let $l$, $m \in \{1, \ldots, 4\}$ be distinct. Consider the equations
\begin{equation*}
\begin{aligned}
    &\text{(a)}\; \frac{\sin{\alpha_l}\sin{\gamma_l}}{\sin{\beta_l}\sin{\delta_l}} = \frac{\sin{\alpha_m}\sin{\gamma_m}}{\sin{\beta_m}\sin{\delta_m}},
    \qquad
    &&\text{(b)}\; \frac{\sin{\overline{\alpha}_l}\sin{\overline{\gamma}_l}}{\sin{\overline{\beta}_l}\sin{\overline{\delta}_l}} = \frac{\sin{\overline{\alpha}_m}\sin{\overline{\gamma}_m}}{\sin{\overline{\beta}_m}\sin{\overline{\delta}_m}}, \\
    &\text{(c)}\; \frac{\sin{\alpha_l}\sin{\gamma_l}}{\sin{\beta_l}\sin{\delta_l}} = \frac{\sin{\beta_m}\sin{\delta_m}}{\sin{\alpha_m}\sin{\gamma_m}},
    \qquad
    &&\text{(d)}\; \frac{\sin{\overline{\alpha}_l}\sin{\overline{\gamma}_l}}{\sin{\overline{\beta}_l}\sin{\overline{\delta}_l}} = \frac{\sin{\overline{\beta}_m}\sin{\overline{\delta}_m}}{\sin{\overline{\alpha}_m}\sin{\overline{\gamma}_m}}.
\end{aligned}
\end{equation*}
Then:
\begin{itemize}
    \item[(i)] $f_l = f_m$ if and only if (a) and (b) hold;
    \item[(ii)] $f_l f_m = M$ if and only if (c) and (d) hold;
    \item[(iii)] $f_l + f_m - f_l f_m = M$ if and only if (a) and (d) hold;
    \item[(iv)] $f_l f_m = M(f_l + f_m - 1)$ if and only if (b) and (c) hold.
\end{itemize}
\end{lemma}

\begin{proof}
First, by~\eqref{eq:N1}, $\overline{\alpha}_i$,
$\overline{\beta}_i$, $\overline{\gamma}_i$,
$\overline{\delta}_i \notin \pi\mathbb{Z}$, see~\eqref{eq:N2},
and, by~\eqref{eq:N4} and Lemma~2 of~\cite{NURMATOV2026104102},
$f_i$, $M$, $M/f_i\neq 0, 1$ are well-defined reals.
Since, in addition, $\alpha_i$, $\beta_i$,
$\gamma_i$, $\delta_i \in (0, \pi)$, it follows that both sides of
the identities below are well-defined. 
For each $i = 1, \ldots, 4$, we have
\begin{equation*}
    \Omega_i \coloneqq \frac{\sin{\alpha_i}\sin{\gamma_i}}{\sin{\beta_i}\sin{\delta_i}} = \frac{f_i(f_i - M)}{M(f_i - 1)},
    \qquad
    \overline{\Omega}_i \coloneqq \frac{\sin{\overline{\alpha}_i}\sin{\overline{\gamma}_i}}{\sin{\overline{\beta}_i}\sin{\overline{\delta}_i}} = \frac{f_i - M}{f_i(f_i - 1)},
\end{equation*}
as verified in Section~1 of the supplement (\texttt{helper.nb})~\cite{Nurmatov2026Repo}. 
Observe that $\Omega_i > 0$ 
because $\alpha_i$, $\beta_i$, $\gamma_i$, $\delta_i \in (0, \pi)$.
Moreover, $f_i^2\neq M$; otherwise $\Omega_i = -1$,
contradicting $\Omega_i > 0$.
In this notation, equations (a)--(d) are $\Omega_l = \Omega_m$,
$\overline{\Omega}_l = \overline{\Omega}_m$,
$\Omega_l\,\Omega_m = 1$, and
$\overline{\Omega}_l\,\overline{\Omega}_m = 1$, respectively. 

Next, using the expressions for $\Omega_i$ and $\overline{\Omega}_i$ with $i \in \{l, m\}$, 
we get the following \emph{$\Omega$-identities}:
\begin{equation*}
\begin{aligned}
    \Omega_l - \Omega_m &= \frac{(f_l - f_m)\,\bigl(f_l f_m - (f_l + f_m) + M\bigr)}{M(f_l - 1)(f_m - 1)},
    \\
    \Omega_l\,\Omega_m - 1 &= \frac{(f_l f_m - M)\,\bigl(f_l f_m - M(f_l + f_m - 1)\bigr)}{M^2(f_l - 1)(f_m - 1)}, 
    \\
    \overline{\Omega}_l - \overline{\Omega}_m &= \frac{(f_m - f_l)\,\bigl(f_l f_m - M(f_l + f_m - 1)\bigr)}{f_l f_m (f_l - 1)(f_m - 1)},
    \\
    \overline{\Omega}_l\,\overline{\Omega}_m - 1 &= \frac{(f_l f_m - M)\,(f_l + f_m - f_l f_m - M)}{f_l f_m (f_l - 1)(f_m - 1)},
\end{aligned}
\end{equation*}
as verified in Section~2 of the supplement (\texttt{helper.nb})~\cite{Nurmatov2026Repo}. 
Since $M \neq 0$ and $f_l$, $f_m \neq 0, 1$, it
follows that each of the four equations (a)--(d) 
holds if and only if the
numerator of the right-hand side of the corresponding
$\Omega$-identity vanishes.
Hence (a) holds if and
only if $f_l = f_m$ or $f_l + f_m - f_l f_m = M$; (b), if and only
if $f_l = f_m$ or $f_l f_m = M(f_l + f_m - 1)$; (c), if and only
if $f_l f_m = M$ or $f_l f_m = M(f_l + f_m - 1)$; and (d), if and
only if $f_l f_m = M$ or $f_l + f_m - f_l f_m = M$.

Finally, observe that the equalities $f_l + f_m - f_l f_m = M$ and
$f_l f_m = M(f_l + f_m - 1)$ cannot hold simultaneously.
Indeed, otherwise, since $M\neq 1$, it follows that $f_lf_m=-M$ and 
$f_l=-f_m$, which implies $f_l^2 = M$, contradicting $f_l^2 \neq M$. 
Analogously, $f_l = f_m$ and $f_l f_m = M$ cannot hold simultaneously.
Thus, the statement of the lemma follows.
\end{proof}

\begin{remark}\label{rem:lem:gqs}
For distinct
$l$, $m \in \{1, \ldots, 4\}$, the equalities
$f_l + f_m - f_l f_m = M$ and $f_l f_m = M(f_l + f_m - 1)$
cannot hold simultaneously under the hypotheses of Lemma~\ref{lem:gqs}, 
as shown in its proof.
\end{remark}

In the remainder of this section,
$\operatorname{sn}$, $\operatorname{cn}$, and $\operatorname{dn}$ 
denote the Jacobi elliptic functions with modulus $k$, see~\eqref{eq:N7}.
We set $\Lambda' \coloneqq 2K\mathbb{Z} + 2\mathbf{i}K'\mathbb{Z}$,
so that, for $k \in (0, 1)$, the set $\mathbf{i}K' + \Lambda'$ is
the pole set of $\operatorname{sn}$, $\operatorname{cn}$,
and $\operatorname{dn}$~\cite[\S 22.35]{WhittakerWatson2021};
accordingly, we call a point of $\mathbb{C}$ a \emph{pole} 
if it belongs to $\mathbf{i}K' + \Lambda'$. 
Throughout, we use without further mention that $k \in (0, 1)$ implies
$k' \in (0, 1)$ and $K$, $K' > 0$, see~\eqref{eq:N7} and~\eqref{eq:N8}.

\begin{lemma}\label{lem:shifts:1}
Let $k \in (0, 1)$ and $v_1$, $v_2$, $v_3$,
$v_4 \in \mathbb{C}\setminus(\mathbf{i}K' + \Lambda')$, and define
\begin{equation*}
    G \coloneqq k'^2
    - k^2 k'^2 \prod_{i = 1}^{4} \operatorname{sn}{v_i}
    + k^2 \prod_{i = 1}^{4} \operatorname{cn}{v_i}
    - \prod_{i = 1}^{4} \operatorname{dn}{v_i}.
\end{equation*}
Then $G = 0$ if and only if
\begin{equation*}
    \eta_1 v_1 + \eta_2 v_2 + \eta_3 v_3 + \eta_4 v_4
    \in 4K\mathbb{Z} + 4\mathbf{i}K'\mathbb{Z}
\end{equation*}
for some $\eta_1$, $\eta_2$, $\eta_3$, $\eta_4 \in \{\pm 1\}$ with
$\eta_1 \eta_2 \eta_3 \eta_4 = 1$.
\end{lemma}

\begin{proof}
First, set $\tilde{v}_i \coloneqq \pi v_i/(2K)$ for $i = 1, \ldots, 4$, and let
$\vartheta_1$, $\vartheta_2$, $\vartheta_3$, $\vartheta_4$ be the
Jacobi theta functions with lattice parameter
$\tau \coloneqq \mathbf{i}K'/K$ and nome $e^{-\pi K'/K}$,
see~\cite[\S\S 21.1, 21.11]{WhittakerWatson2021}. Then,
by~\cite[\S\S 21.61, 21.7, 22.11, 22.3, and 22.32]{WhittakerWatson2021},
\begin{equation*}
    \operatorname{sn}{v_i}
    = \frac{\vartheta_3(0)\, \vartheta_1(\tilde{v}_i)}
    {\vartheta_2(0)\, \vartheta_4(\tilde{v}_i)},
    \qquad
    \operatorname{cn}{v_i}
    = \frac{\vartheta_4(0)\, \vartheta_2(\tilde{v}_i)}
    {\vartheta_2(0)\, \vartheta_4(\tilde{v}_i)},
    \qquad
    \operatorname{dn}{v_i}
    = \frac{\vartheta_4(0)\, \vartheta_3(\tilde{v}_i)}
    {\vartheta_3(0)\, \vartheta_4(\tilde{v}_i)},
\end{equation*}
with $k = \vartheta_2^2(0)/\vartheta_3^2(0)$ and
$k' = \vartheta_4^2(0)/\vartheta_3^2(0)$, where $\vartheta_2(0)$,
$\vartheta_3(0) \neq 0$, see~\cite[\S 21.12]{WhittakerWatson2021}.
Since $\tilde{v}_i = \pi v_i/(2K)$, $K\tau = \mathbf{i}K'$, and
$v_i \notin \mathbf{i}K' + \Lambda'$, it follows that
$\vartheta_4(\tilde{v}_i) \neq 0$ for each $i \in \{1, 2, 3, 4\}$, 
because the zeros of $\vartheta_4$ are the points of
$\pi\tau/2 + \pi\mathbb{Z} + \pi\tau\mathbb{Z}$,
see~\cite[\S 21.12]{WhittakerWatson2021}.
Thus, taking the product over $i = 1, \ldots, 4$ in each of the 
three displayed identities and substituting $k$ and $k'$, we get
\begin{equation*}
    G= \frac{k'^2}{\prod_{i = 1}^{4} \vartheta_4(\tilde{v}_i)}
    \biggl(
    \prod_{i = 1}^{4} \vartheta_4(\tilde{v}_i)
    - \prod_{i = 1}^{4} \vartheta_1(\tilde{v}_i)
    + \prod_{i = 1}^{4} \vartheta_2(\tilde{v}_i)
    - \prod_{i = 1}^{4} \vartheta_3(\tilde{v}_i)
    \biggr).
\end{equation*}

Next, by~\cite[\S 21.22, formula~(iv)]{WhittakerWatson2021} and the parity
of the theta functions, see~\cite[\S 21.11]{WhittakerWatson2021}, we have
\begin{equation*}
    \prod_{i = 1}^{4} \vartheta_4(\tilde{v}_i)
    - \prod_{i = 1}^{4} \vartheta_1(\tilde{v}_i)
    + \prod_{i = 1}^{4} \vartheta_2(\tilde{v}_i)
    - \prod_{i = 1}^{4} \vartheta_3(\tilde{v}_i)
    = -2\prod_{i = 1}^{4} \vartheta_1(w_i),
\end{equation*}
where $2w_1 \coloneqq \tilde{v}_1 + \tilde{v}_2 + \tilde{v}_3 + \tilde{v}_4$,
$2w_2 \coloneqq \tilde{v}_1 + \tilde{v}_2 - \tilde{v}_3 - \tilde{v}_4$,
$2w_3 \coloneqq \tilde{v}_1 - \tilde{v}_2 + \tilde{v}_3 - \tilde{v}_4$, and
$2w_4 \coloneqq \tilde{v}_1 - \tilde{v}_2 - \tilde{v}_3 + \tilde{v}_4$;
here formula~(iv) is applied, in the notation
of~\cite[\S 21.22]{WhittakerWatson2021}, with
$(w, x, y, z) \coloneqq (-w_1, w_2, w_3, w_4)$, for which
$(w', x', y', z') = (\tilde{v}_1, -\tilde{v}_2, -\tilde{v}_3,
-\tilde{v}_4)$. 

Finally, since $k' \neq 0$ and
$\vartheta_4(\tilde{v}_i) \neq 0$ for each $i \in \{1, 2, 3, 4\}$, it follows that 
$G=0$ if and only if
$\vartheta_1(w_i) = 0$ for some $i \in \{1, 2, 3, 4\}$. The zeros of $\vartheta_1$ are
the points of $\pi\mathbb{Z} + \pi\tau\mathbb{Z}$,
see~\cite[\S 21.12]{WhittakerWatson2021}. Hence
$\vartheta_1(w_i) = 0$ if and only if
$2w_i \in 2\pi\mathbb{Z} + 2\pi\tau\mathbb{Z}$, that is, if and only
if $4Kw_i/\pi \in 4K\mathbb{Z} + 4\mathbf{i}K'\mathbb{Z}$. 
Since $\tilde{v}_i = \pi v_i/(2K)$, it follows that
$4Kw_1/\pi = v_1 + v_2 + v_3 + v_4$,
$4Kw_2/\pi = v_1 + v_2 - v_3 - v_4$,
$4Kw_3/\pi = v_1 - v_2 + v_3 - v_4$, and
$4Kw_4/\pi = v_1 - v_2 - v_3 + v_4$. 
Thus, since
$-(4K\mathbb{Z} + 4\mathbf{i}K'\mathbb{Z})
= 4K\mathbb{Z} + 4\mathbf{i}K'\mathbb{Z}$, it follows that
$4Kw_i/\pi \in 4K\mathbb{Z} + 4\mathbf{i}K'\mathbb{Z}$ for some
$i \in \{1, 2, 3, 4\}$ if and only if
$\eta_1 v_1 + \eta_2 v_2 + \eta_3 v_3 + \eta_4 v_4
\in 4K\mathbb{Z} + 4\mathbf{i}K'\mathbb{Z}$ for some $\eta_1$,
$\eta_2$, $\eta_3$, $\eta_4 \in \{\pm 1\}$ with
$\eta_1 \eta_2 \eta_3 \eta_4 = 1$.
\end{proof}

\begin{lemma}\label{lem:shifts:2}
Assume that $\alpha_i$, $\beta_i$, $\gamma_i$,
$\delta_i \in (0, \pi)$ satisfy conditions~\eqref{eq:N1}
and~\eqref{eq:N4} for each $i = 1, \ldots, 4$, and $M > 0$.
Then, for each $i = 1, \ldots, 4$, 
\begin{equation*}
\sign{(\operatorname{dn}^2{t_i} - 1)}
= \sign{(\operatorname{dn}^2{t_i} - k'^2)}
= \sign{(p_i^2 q_i^2)} \neq 0,
\end{equation*}
and
\begin{equation*}
    \operatorname{sn}{t_i}
    = \frac{\mathbf{i}\,\varepsilon_i}{k\, p_i q_i}
    \sqrt{(\operatorname{dn}^2{t_i} - 1)\, p_i^2 q_i^2},
    \qquad
    \operatorname{cn}{t_i}
    = \frac{\varepsilon_i}{k\, p_i q_i}
    \sqrt{(\operatorname{dn}^2{t_i} - k'^2)\, p_i^2 q_i^2}.
\end{equation*}
\end{lemma}

\begin{proof}
First, by~\eqref{eq:N4}
and Lemmas~2 and~4 of~\cite{NURMATOV2026104102}, 
$r_i$, $s_i$, $f_i$, $M$, $M/f_i\neq 0, 1$ are well-defined reals, 
$f_i > 0$, $\overline{\alpha}_i$, $\overline{\beta}_i$,
$\overline{\gamma}_i$, $\overline{\delta}_i \in (0, \pi)$, 
and, by~\eqref{eq:N1} and~\eqref{eq:N2},
$\varepsilon_i \in \{\pm 1\}$. 
Hence $k$, $k' \in (0, 1)$, because $M > 0$ and $M \neq 1$, see~\eqref{eq:N7}, 
and the phase shift $t_i$ is
well-defined, see~\eqref{eq:N11} and
Table~\ref{tab:tpqsigma}. 
Then $p_i$, $q_i\in \mathbb{R}_{>0} \cup \mathbf{i} \mathbb{R}_{>0}$, and
$p_i^2$, $q_i^2$, $\operatorname{dn}^2{t_i} - k'^2$, $\operatorname{dn}^2{t_i} - 1$ 
are well-defined nonzero reals, and 
$(\operatorname{dn}^2{t_i} - 1)\, p_i^2 q_i^2$, 
$(\operatorname{dn}^2{t_i} - k'^2)\, p_i^2 q_i^2$ are strictly positive
by~\eqref{eq:N7}, \eqref{eq:N10}, \eqref{eq:N11}, and Lemmas~2 and~4
of~\cite{NURMATOV2026104102}. 
Thus $\sign{(\operatorname{dn}^2{t_i} - 1)} = \sign{(\operatorname{dn}^2{t_i} - k'^2)} = \sign{(p_i^2 q_i^2)}\neq 0$.

Second, recall that $t_i$ is determined by~\eqref{eq:N11} and
Table~\ref{tab:tpqsigma}; the table places $t_i$ in accordance with
the value of the elliptic cosine, where
$\operatorname{cn}{t_i}
= \frac{\operatorname{dn}{t_i}\sin{\sigma_i}}
{p_i q_i \sin{\overline{\delta}_i}}$ when $M < 1$, and
$\operatorname{cn}{t_i}
= \frac{\sin{\sigma_i}}{p_i q_i \sin{\overline{\delta}_i}}$ when
$M > 1$ (see~\cite[Eq.~(31)]{izmestiev-2017} and
Remark~\ref{rem:lem:shifts:2} below). 

Finally, the formulas stated in the lemma follow
from these values of $\operatorname{cn}{t_i}$
by~\eqref{eq:N2}, \eqref{eq:N3}, \eqref{eq:N7}, \eqref{eq:N10}, \eqref{eq:N11},
and a case-by-case comparison of Table~\ref{tab:tpqsigma} with
Table~\ref{tab:sncn}, over the three classes of $p_i q_i$ and
the two ranges of $\sigma_i$, as verified in Section~3 of the
supplement
(\texttt{helper.nb})~\cite{Nurmatov2026Repo}.
\end{proof}

\begin{remark}\label{rem:lem:shifts:2}
Equation~(31) of~\cite{izmestiev-2017} contains a typographical
error for $M < 1$; the case $M > 1$ is stated correctly.
Specifically, for $M < 1$, equating the coefficients of $xy$ in
equations~(28) and~(30) of~\cite{izmestiev-2017} gives
\begin{equation*}
    \operatorname{cn}{t_0}\operatorname{dn}{t_0}
    = \frac{pq \sin{\alpha}\sin{\gamma}\sin{\overline{\delta}}}
    {\sin{\sigma}\sin{(\sigma - \alpha - \beta)}
    \sin{(\sigma - \beta - \gamma)}}
    = \frac{f \sin{\sigma}}{pq \sin{\overline{\delta}}}
    = \frac{\operatorname{dn}^2{t_0}\,\sin{\sigma}}
    {pq \sin{\overline{\delta}}},
\end{equation*}
where $t_0$ denotes the phase shift in the notation
of~\cite{izmestiev-2017} and the vertex index is suppressed; the
second equality holds by Lemma~1 of~\cite{NURMATOV2026104102},
using~\eqref{eq:N2}, \eqref{eq:N3}, and~\eqref{eq:N10}, and the
third by~\eqref{eq:N11}. Hence
$\operatorname{cn}{t_0}
= \frac{\operatorname{dn}{t_0}\sin{\sigma}}
{pq \sin{\overline{\delta}}}$, because
$\operatorname{dn}^2{t_0} = f \neq 0$, see~\eqref{eq:N11} and
Lemma~2 of~\cite{NURMATOV2026104102}. 
This error does not affect the table in Theorem~2.7
of~\cite{izmestiev-2017} (Table~\ref{tab:tpqsigma}), because the table
is determined by $\operatorname{cn}{t_0}/|\operatorname{cn}{t_0}|$, see
Table~\ref{tab:sncn}, and the value given by equation~(31) equals the
corrected value divided by $\operatorname{dn}{t_0} = \sqrt{f} > 0$,
see~\eqref{eq:N11} and Lemma~4 of~\cite{NURMATOV2026104102}.
Throughout this paper, we adopt the corrected formula, as applied
in the proof of Lemma~\ref{lem:shifts:2}.
\end{remark}

In the proofs below, we use the parity and half-period 
shifts~\cite[\S\S 22.12, 22.31, and 22.34]{WhittakerWatson2021}:
for all $m, n \in \mathbb{Z}$ and $v \in \mathbb{C}$,
\begin{equation}\label{eq:P1}
    (\operatorname{sn},\operatorname{cn},\operatorname{dn})(2mK + 2n\mathbf{i}K' \pm v)
    = \bigl(\pm(-1)^{m}\operatorname{sn}{v}, \,(-1)^{m+n}\operatorname{cn}{v}, \,(-1)^{n}\operatorname{dn}{v}\bigr).
\end{equation}

\begin{lemma}\label{lem:shifts:3}
Assume that $\alpha_i$, $\beta_i$, $\gamma_i$, $\delta_i \in (0, \pi)$
satisfy conditions~\eqref{eq:N1}, \eqref{eq:N4}, and~\eqref{eq:N5} for
each $i = 1, \ldots, 4$, and $M > 0$. 
Then condition~\eqref{eq:N6} holds if and only if 
$\tilde{G} = 0$ for some
$e \in \{\pm 1\}$ and $n \in \{0, 1\}$, where
$\varrho \coloneqq 0$ if $M < 1$ and $\varrho \coloneqq n$ if $M > 1$, and
\begin{equation*}
\begin{aligned}
    \tilde{G}\coloneqq k^2 k'^2
    &- (-1)^{\varrho}\, e\, k'^2 \chi
    \sqrt{\textstyle\prod_{i = 1}^{4}(\operatorname{dn}^2{t_i} - 1)}
    \\
    &+ (-1)^{n + \varrho}\, \chi
    \sqrt{\textstyle\prod_{i = 1}^{4}(\operatorname{dn}^2{t_i} - k'^2)}
    - (-1)^{n}\, k^2 \prod_{i = 1}^{4} \operatorname{dn}{t_i}.
\end{aligned}
\end{equation*}
Moreover, both radicands in the definition of $\tilde{G}$ are
strictly positive reals.
\end{lemma}

\begin{proof}
First, by~\eqref{eq:N4} and Lemmas~2 and~4 of~\cite{NURMATOV2026104102},
$r_i$, $s_i$, $f_i$, $M\neq 0,1$ are well-defined reals and $f_i > 0$.
Hence $k$, $k' \in (0, 1)$, because $M > 0$ and $M \neq 1$, see~\eqref{eq:N7}, 
and the phase shift $t_i$ is
well-defined, see~\eqref{eq:N11} and
Table~\ref{tab:tpqsigma}.
Also, by~\eqref{eq:N3} and~\eqref{eq:N14},
$u$, $x_i$, $y_i$, $z_i \neq 0$ are well-defined reals,
and $\chi \in \{\pm 1\}$, because
$\varepsilon_i \in \{\pm 1\}$
by~\eqref{eq:N1} and~\eqref{eq:N2}, and $z_i \neq 0$, see~\eqref{eq:N12}.
Now, fix $n \in \{0, 1\}$ and $e_1$, $e_2$, $e_3$,
$e_4 \in \{\pm 1\}$, set $e \coloneqq e_1 e_2 e_3 e_4$, and denote
$v_i \coloneqq e_i t_i$ for $i = 1, 2, 3$ and
$v_4 \coloneqq e_4 t_4 - \lambda$, where $\lambda \coloneqq 2n\mathbf{i}K'$ if
$M < 1$ and $\lambda \coloneqq 2nK + 2n\mathbf{i}K'$ if $M > 1$. 

Next, observe that $v_1$, $v_2$, $v_3$, $v_4 \in \mathbb{C}\setminus(\mathbf{i}K' + \Lambda')$. 
Indeed, $t_i \notin \mathbf{i}K' + \Lambda'$ 
because $\operatorname{Im}{t_i} \in (0, K')$, see~\eqref{eq:N11} and Table~\ref{tab:tpqsigma}.
Moreover, $\mathbf{i}K' + \Lambda'$ is invariant under translation by $\lambda$ (because $\lambda \in \Lambda'$) 
and under reflection (since $\Lambda' = -\Lambda'$ and $2\mathbf{i}K' \in \Lambda'$, it follows that 
$-(\mathbf{i}K' + \Lambda') = -\mathbf{i}K' + \Lambda' = \mathbf{i}K' + \Lambda'$). 
Hence Lemma~\ref{lem:shifts:1} applies to $v_1$, $v_2$, $v_3$, $v_4$. Let $G$ be as in Lemma~\ref{lem:shifts:1}.

We now show that $\tilde{G} = k^2 G$. 
Observe first that $\prod_{i = 1}^{4} p_i q_i = p_1^2\, p_3^2\, q_1^2\, q_2^2
= 1/(x_1 x_3 y_1 y_2)$ by~\eqref{eq:N3}, \eqref{eq:N5},
and~\eqref{eq:N10}. By~\eqref{eq:N3}, \eqref{eq:N4}, \eqref{eq:N14}, and Lemma~4
of~\cite{NURMATOV2026104102}, we have $u x_i y_i z_i > 0$ for
$i \in \{1, 2, 3, 4\}$.
Hence
$\sign{(z_1 z_3)}
= \sign{(u x_1 y_1 u x_3 y_3)}
= \sign{(x_1 x_3 y_1 y_2)}$, because $y_3 = y_2$ by~\eqref{eq:N3}
and~\eqref{eq:N5}. 
Thus,
$\sign{\bigl(\prod_{i = 1}^{4} p_i q_i\bigr)}
= \sign{(x_1 x_3 y_1 y_2)}
= \sign{(z_1 z_3)}$.
Since $\bigl(\prod_{i = 1}^{4} p_i q_i\bigr)^2 > 0$,
it follows, by Lemma~\ref{lem:shifts:2}, that
$\prod_{i = 1}^{4}(\operatorname{dn}^2{t_i} - 1)$ and
$\prod_{i = 1}^{4}(\operatorname{dn}^2{t_i} - k'^2)$ are strictly
positive reals, and 
\begin{equation*}
    k^4 \prod_{i = 1}^{4} \operatorname{sn}{t_i}
    = \chi \sqrt{\textstyle\prod_{i = 1}^{4}
    (\operatorname{dn}^2{t_i} - 1)},
    \qquad
    k^4 \prod_{i = 1}^{4} \operatorname{cn}{t_i}
    = \chi \sqrt{\textstyle\prod_{i = 1}^{4}
    (\operatorname{dn}^2{t_i} - k'^2)},
\end{equation*}
where $\varepsilon_1 \varepsilon_2 \varepsilon_3 \varepsilon_4\,
\sign{(z_1 z_3)} = \chi$ by~\eqref{eq:N12}.
Thus, the expression for $\tilde{G}$ becomes
\begin{equation*}
    \tilde{G} = k^2 \Bigl(k'^2
    - (-1)^{\varrho} \, e\, k'^2 k^2
    \prod_{i = 1}^{4} \operatorname{sn}{t_i}
    + (-1)^{n + \varrho} \, k^2
    \prod_{i = 1}^{4} \operatorname{cn}{t_i}
    - (-1)^{n} \prod_{i = 1}^{4} \operatorname{dn}{t_i}\Bigr).
\end{equation*}
Next, by~\eqref{eq:P1} and the definition of~$\varrho$, we have
$(\operatorname{sn}, \operatorname{cn}, \operatorname{dn}){(v_i)}
= (e_i \operatorname{sn}{t_i}, \,\operatorname{cn}{t_i}, \,\operatorname{dn}{t_i})$
for $i = 1, 2, 3$, and
$(\operatorname{sn}, \operatorname{cn}, \operatorname{dn}){(v_4)}
= \bigl((-1)^{\varrho} e_4 \operatorname{sn}{t_4}, \,(-1)^{n + \varrho}\operatorname{cn}{t_4}, \,(-1)^{n}\operatorname{dn}{t_4}\bigr)$.
Multiplying these identities over $i = 1, \ldots, 4$, we get
$\tilde{G} = k^2 G$. Since $k \neq 0$, it follows that
$\tilde{G} = 0$ if and only if $G = 0$.

($\Rightarrow$) Assume that condition~\eqref{eq:N6} holds, that is,
$e_1 t_1 + e_2 t_2 + e_3 t_3 + e_4 t_4 \in \Lambda$ for some
$e_1, \ldots, e_4 \in \{\pm 1\}$. By~\eqref{eq:N9}, there are
$m$, $n' \in \mathbb{Z}$ such that
$\sum_{i = 1}^{4} e_i t_i = 4mK + 2n'\mathbf{i}K'$ if $M < 1$, and
$\sum_{i = 1}^{4} e_i t_i = 4mK + (2K + 2\mathbf{i}K')n'$ if
$M > 1$. Choose $n \in \{0, 1\}$ with $n \equiv n' \pmod{2}$. Then
$\sum_{i = 1}^{4} e_i t_i - \lambda
\in 4K\mathbb{Z} + 4\mathbf{i}K'\mathbb{Z}$, because $n' - n$ is
even. Hence
$\eta_1 e_1 t_1 + \eta_2 e_2 t_2 + \eta_3 e_3 t_3
+ \eta_4 (e_4 t_4 - \lambda)
\in 4K\mathbb{Z} + 4\mathbf{i}K'\mathbb{Z}$ with
$\eta_1 = \eta_2 = \eta_3 = \eta_4 = 1$. Thus, by
Lemma~\ref{lem:shifts:1}, $G = 0$, and hence $\tilde{G} = 0$
for $e = e_1 e_2 e_3 e_4$ and this $n$.

($\Leftarrow$) Assume that $\tilde{G} = 0$ for some $e$ and $n$, and choose
$e_1$, $e_2$, $e_3$, $e_4 \in \{\pm 1\}$ with
$e_1 e_2 e_3 e_4 = e$. Then $G = 0$, and hence, by
Lemma~\ref{lem:shifts:1}, there are
$\eta_1, \ldots, \eta_4 \in \{\pm 1\}$ with
$\eta_1\eta_2\eta_3\eta_4 = 1$ such that
$\eta_1 e_1 t_1 + \eta_2 e_2 t_2 + \eta_3 e_3 t_3
+ \eta_4 (e_4 t_4 - \lambda)
\in 4K\mathbb{Z} + 4\mathbf{i}K'\mathbb{Z}$. Since
$\eta_4 \lambda - \lambda \in \{0, -2\lambda\}$ and
$2\lambda \in 4K\mathbb{Z} + 4\mathbf{i}K'\mathbb{Z}$, it follows
that $\sum_{i = 1}^{4} (\eta_i e_i)\, t_i
    \in \lambda + 4K\mathbb{Z} + 4\mathbf{i}K'\mathbb{Z}
    \subseteq \Lambda$,
because $\lambda \in \Lambda$ and
$4K\mathbb{Z} + 4\mathbf{i}K'\mathbb{Z} \subseteq \Lambda$,
see~\eqref{eq:N9}; for $M < 1$ this is immediate, and for $M > 1$
it follows from $4\mathbf{i}K' = 2(2K + 2\mathbf{i}K') - 4K$.
Hence condition~\eqref{eq:N6} holds with the signs
$\eta_1 e_1, \ldots, \eta_4 e_4$.
\end{proof}

\begin{proof}[Proof of Proposition~\ref{prop:shifts}]
First, by~\eqref{eq:N4} and Lemmas~2 and~4
of~\cite{NURMATOV2026104102}, $f_1, \ldots, f_4 \neq 1$ 
are well-defined strictly positive reals, $M \neq 1$, 
and $\chi \in \{\pm 1\}$, because
$\varepsilon_1, \ldots, \varepsilon_4 \in \{\pm 1\}$
by~\eqref{eq:N1} and~\eqref{eq:N2}, and $z_1$, $z_3 \neq 0$
by~\eqref{eq:N3}, see~\eqref{eq:N12}.
Next, fix $\zeta_1$, $\zeta_2 \in \{\pm 1\}$, set $e \coloneqq \zeta_1\,\chi$,
and choose $n \in \{0, 1\}$ such that $(-1)^n = \zeta_2$; 
this is a bijection between the pairs
$\zeta_1$, $\zeta_2$ and the pairs $e$, $n$, because $\chi \in \{\pm 1\}$.
By~\eqref{eq:N7} and~\eqref{eq:N11}, for each $i \in \{1, \ldots, 4\}$, 
$\operatorname{dn}{t_i} = \sqrt{f_i}$, 
$k^2 = 1 - M$, $k'^2 = M$ if $M < 1$, 
and $\operatorname{dn}{t_i} = 1/\sqrt{f_i}$,
$k^2 = 1 - 1/M$, $k'^2 = 1/M$ if $M > 1$. 
Moreover, by Lemma~\ref{lem:shifts:3},
$\prod_{i = 1}^{4}(\operatorname{dn}^2{t_i} - 1)$, $\prod_{i = 1}^{4}(\operatorname{dn}^2{t_i} - k'^2)>0$. 
Hence, in both cases,
$\prod_{i = 1}^{4}(f_i - 1)$,
$\prod_{i = 1}^{4}(f_i - M) > 0$, because $f_i$, $M > 0$. 
Thus, all three radicands in the proposition are strictly positive reals.
Further, the displayed equation of the proposition 
is obtained from the equation $\tilde{G} = 0$ 
of Lemma~\ref{lem:shifts:3}, evaluated at
these $e$ and $n$, by multiplying by $-1$ if $M < 1$ and by
$-\zeta_2\, M^2 \sqrt{f_1 f_2 f_3 f_4}$ if $M > 1$, and moving the
term $-M(1-M)$ to the right-hand side. 
Since $M^2 \sqrt{f_1 f_2 f_3 f_4} \neq 0$, it follows that, for each
$\zeta_1$ and $\zeta_2$, the two equations are equivalent.
Thus, by Lemma~\ref{lem:shifts:3}, condition~\eqref{eq:N6} holds if and
only if the displayed equation of the proposition holds for some
$\zeta_1$, $\zeta_2 \in \{\pm 1\}$.
\end{proof}

\begin{proof}[Proof of Theorem~\ref{th:gqs}]
Conditions~\eqref{eq:N1}, \eqref{eq:N4}, and~\eqref{eq:N5} hold by
Definition~\ref{def:gqs}, so it suffices to verify condition~\eqref{eq:N6}. 
First, by~\eqref{eq:N4} and Lemmas~2, 4,
and~5 of~\cite{NURMATOV2026104102}, $M$, $f_i$, $M/f_i \neq 0, 1$
are well-defined positive reals and $(f_i - 1)(f_i - M) > 0$ for
each $i \in \{1, 2, 3, 4\}$.
Also, $\chi \in \{\pm 1\}$, because
$\varepsilon_1, \ldots, \varepsilon_4 \in \{\pm 1\}$
by~\eqref{eq:N1} and~\eqref{eq:N2}, and $z_1$, $z_3 \neq 0$
by~\eqref{eq:N3}, see~\eqref{eq:N12}. 
Now fix indices
$\{i, j, l, m\} = \{1, 2, 3, 4\}$, as in Definition~\ref{def:gqs}.
In each of the cases (i)--(iv) of Definition~\ref{def:gqs}, the
displayed equation of Proposition~\ref{prop:shifts} holds for
$\zeta_1 = \chi$ and $\zeta_2 = 1$ in cases (i)--(iii) and for
$\zeta_1 = \chi$ and $\zeta_2 = -\chi\,\sign{(1 - M)}$ in case (iv), by
elementary algebra from the conditions of each case, as verified in
Section~4 of the supplement (\texttt{helper.nb})~\cite{Nurmatov2026Repo}.
Thus, by Proposition~\ref{prop:shifts}, condition~\eqref{eq:N6} holds, and
the polyhedron has equimodular elliptic type.
\end{proof}

We now present the lemmas needed to prove Theorem~\ref{th:main}. First, we
describe all solutions of Bricard's system~\eqref{eq:N13} over $\mathbb{C}$
for flat angles satisfying conditions~\eqref{eq:N1}, \eqref{eq:N4},
and~\eqref{eq:N5} with $M > 0$, see Lemma~\ref{lem:main:2}, which relies
on Lemmas~\ref{lem:main:1} and~\ref{lem:main:1-2} below: 
the solutions form one-parameter families,
one for each choice of signs satisfying condition~\eqref{eq:N6}, together
with finitely many further solutions, see Remark~\ref{rem:lem:main:2}. 
Hereafter, we set $F \coloneqq \operatorname{sn}$ and
$\Gamma \coloneqq 2K + \Lambda$ if $M < 1$, and $F \coloneqq \operatorname{cn}$
and $\Gamma \coloneqq \Lambda$ if $M > 1$.
We also use without further mention that an elliptic function of order $N$
attains every value of $\mathbb{C} \cup \{\infty\}$ exactly $N$ times
modulo its period lattice, counting multiplicity,
see~\cite[\S 20.13]{WhittakerWatson2021}. Recall that the order of an
elliptic function is the number of its poles modulo its period
lattice, counted with multiplicity.

\begin{lemma}\label{lem:main:1}
Assume that $\alpha_i$, $\beta_i$, $\gamma_i$, $\delta_i \in (0, \pi)$
satisfy conditions~\eqref{eq:N1} and~\eqref{eq:N4} for each
$i = 1, \ldots, 4$, and $M > 0$. Then, for each $i \in \{1, 2, 3, 4\}$ and
each $e \in \{\pm 1\}$, a pair $(X, Y) \in \mathbb{C}^2$ solves the equation
$P_i(X, Y) = 0$ of Bricard's system~\eqref{eq:N13} if and only if
\begin{equation*}
    (X, Y) = \bigl(p_i\,F(v),\ q_i\,F(v + e\,t_i)\bigr)
\end{equation*}
for some $v \in \mathbb{C}$ at which both coordinates are finite.
\end{lemma}

\begin{proof}
First, by~\eqref{eq:N4} and Lemmas~2 and~4
of~\cite{NURMATOV2026104102}, $r_i$, $s_i$, $f_i$,
$M \neq 0, 1$ are well-defined reals and $f_i > 0$. Hence
$p_i$, $q_i \in \mathbb{R}_{>0} \cup \mathbf{i}\mathbb{R}_{>0}$,
see~\eqref{eq:N10}, $k \in (0, 1)$, because $M > 0$ and
$M \neq 1$, see~\eqref{eq:N7}, and the phase shift $t_i$ is
well-defined, see~\eqref{eq:N11} and
Table~\ref{tab:tpqsigma}.

Next, it suffices to treat the case $e = +1$. Indeed, by~\eqref{eq:P1},
the replacement of $v$ by $-v$ maps the set of pairs
$(p_i F(v), q_i F(v + t_i))$ onto
the set of pairs $(p_i F(v), q_i F(v - t_i))$ if $F = \operatorname{cn}$,
because $\operatorname{cn}$ is even, and onto the set of pairs
$(-p_i F(v), -q_i F(v - t_i))$ if $F = \operatorname{sn}$, because
$\operatorname{sn}$ is odd. Since $P_i(-X, -Y) = P_i(X, Y)$, see
Definition~\ref{def:bricard}, it follows that 
the solution set of $P_i(X, Y) = 0$ is
invariant under $(X, Y) \mapsto (-X, -Y)$. Hence, in both cases, the
pairs with $e = -1$ form the solution set of $P_i(X, Y) = 0$ if and
only if the pairs with $e = +1$ do.

Finally, let $e = +1$.
By Lemma~4 of~\cite{NURMATOV2026104102}, $\overline{\alpha}_i$,
$\overline{\beta}_i$, $\overline{\gamma}_i$,
$\overline{\delta}_i \in (0, \pi)$,
so, by Lemma~4.1 of~\cite{izmestiev-2017},
there is a spherical quadrilateral with the side lengths $\alpha_i$,
$\beta_i$, $\gamma_i$, $\delta_i$.
By Definition~\ref{def:bricard} together with Section~2.2,
Lemma~4.3, and Definition~4.4 of~\cite{izmestiev-2017}, the solution
set of $P_i(X, Y) = 0$ in $\mathbb{C}^2$ is the affine part of the
complexified configuration space of this quadrilateral,
written in the variables of Definition~\ref{def:bricard}.
Since, by condition~\eqref{eq:N1}, this quadrilateral is of elliptic type
in the sense of Definition~2.3 of~\cite{izmestiev-2017}, it follows, by
Lemma~4.9 and Theorem~2.7 of~\cite{izmestiev-2017}, the former proved
in~\cite[Lemma~3.11]{Izmestiev2015}, that the map
$v \mapsto (p_i F(v), q_i F(v + t_i))$, with $k$, $p_i$, $q_i$, $t_i$ as
in~\eqref{eq:N7}, \eqref{eq:N10}, and~\eqref{eq:N11}, parametrizes this
configuration space, and this parametrization is a bijection of
$\mathbb{C}/\Lambda$ onto this configuration space, see the proof of Lemma~4.10
of~\cite{izmestiev-2017}. Hence the pair
$(p_i F(v), q_i F(v + t_i))$ satisfies $P_i(X, Y) = 0$ for every
$v \in \mathbb{C}$ at which both coordinates are finite, and every
solution $(X, Y) \in \mathbb{C}^2$ of $P_i(X, Y) = 0$ is such a pair for
some $v \in \mathbb{C}$ at which both coordinates are finite.
\end{proof}

\begin{lemma}\label{lem:main:1-2}
Assume that $M > 0$ and $M \neq 1$. Then, for all $v$, $v' \in \mathbb{C}$,
we have $F(v') = F(v)$ in $\mathbb{C} \cup \{\infty\}$ if and only if
$v' - v \in \Lambda$ or $v' + v \in \Gamma$.
\end{lemma}

\begin{proof}
Observe that $k \in (0, 1)$, because $M > 0$ and $M \neq 1$, see~\eqref{eq:N7}.

($\Rightarrow$) The function $F$ has a single simple pole
modulo $\Lambda'$, located at
$\mathbf{i}K'$~\cite[\S 22.35]{WhittakerWatson2021}, and
$\Lambda' = \Lambda \cup (2K + \Lambda)$, see~\eqref{eq:N9}; for $M < 1$
this is immediate, and for $M > 1$ it follows from
$2\mathbf{i}K' - 2K = (2K + 2\mathbf{i}K') - 4K \in \Lambda$. 
Hence $F$ has exactly two simple poles modulo $\Lambda$, located at
$\mathbf{i}K'$ and $2K + \mathbf{i}K'$, which are distinct modulo $\Lambda$
because $2K \notin \Lambda$, see~\eqref{eq:N9}; therefore $F$ is elliptic
of order two with respect to $\Lambda$, so that $F$ attains every value of
$\mathbb{C} \cup \{\infty\}$ exactly twice modulo $\Lambda$, counting
multiplicity.
Now let $F(v') = F(v)$ and $c \in \Gamma$. Since, by~\eqref{eq:P1}, $F(c - v) = F(v)$, it follows that
the function $F$ attains the value $F(v)$ at $v$ and at $c - v$. If these
two points are distinct modulo $\Lambda$, then they are the two solutions
of the equation $F(x) = F(v)$ modulo $\Lambda$, so that $v' - v \in \Lambda$ or
$v' + v \in c + \Lambda = \Gamma$. 
If they coincide modulo $\Lambda$, that is,
$c - v \equiv v \pmod{\Lambda}$, then $v$ is not a pole; otherwise
$v \equiv \mathbf{i}K'$ or $v \equiv 2K + \mathbf{i}K' \pmod{\Lambda}$,
and, since $c \equiv 2K + 2\mathbf{i}K' \pmod{\Lambda}$, see~\eqref{eq:N9}, 
it follows that $c - v \equiv 2K + \mathbf{i}K'$ or
$c - v \equiv \mathbf{i}K' \pmod{\Lambda}$, respectively, so that
$2K \in \Lambda$, which is ruled out by~\eqref{eq:N9}.
Thus $F(v)$ is finite when $c - v \equiv v \pmod{\Lambda}$, and, using
the identities $F(x) = F(c - x)$ and $F(x + \lambda) = F(x)$ for all
$x \in \mathbb{C}$ and $\lambda \in \Lambda$, which hold
by~\eqref{eq:P1}, we get $F(v + h) = F(c - v - h) = F(v - h)$ for every
$h \in \mathbb{C}$; thus the function $h \mapsto F(v + h) - F(v)$ is even,
hence vanishes at $h = 0$ to even order, so that $v$ is a solution of
the equation $F(x) = F(v)$ of multiplicity two, hence the only one modulo $\Lambda$,
and $v' - v \in \Lambda$.

($\Leftarrow$) By~\eqref{eq:P1}, $F(x + \lambda) = F(x)$ and
$F(c - x) = F(x)$ for all $x \in \mathbb{C}$, $\lambda \in \Lambda$,
$c \in \Gamma$; applying the first with $x = v$, $\lambda = v' - v$ and the
second with $x = v$, $c = v' + v$, respectively, gives $F(v') = F(v)$.
\end{proof}

\begin{remark}\label{rem:lem:main:1-2}
By Lemma~\ref{lem:main:1-2}, if $F(v') = F(v)$, then there is
$e \in \{\pm 1\}$ such that $F(v' + s) = F(v + e s)$ for every
$s \in \mathbb{C}$, because $(v' + s) - (v + s) = v' - v$ and
$(v' + s) + (v - s) = v' + v$.
\end{remark}

\begin{lemma}\label{lem:main:2}
Assume that $\alpha_i$, $\beta_i$, $\gamma_i$, $\delta_i \in (0, \pi)$
satisfy conditions~\eqref{eq:N1}, \eqref{eq:N4}, and~\eqref{eq:N5} for each
$i = 1, \ldots, 4$, and $M > 0$. Then $(Z, W_2, U, W_1) \in \mathbb{C}^4$ is
a solution of Bricard's system~\eqref{eq:N13} if and only if there exist
$e_1$, $e_2$, $e_3$, $e_4 \in \{\pm 1\}$ and $w \in \mathbb{C}$ such that
\begin{equation*}
\begin{aligned}
        &Z
        = p_1\,F(w),
        \qquad
        &&W_2
        = q_2\,F(w - e_2 t_2),
        \\
        &U
        = p_3\,F(w - e_2 t_2 - e_3 t_3),
        \qquad
        &&W_1
        = q_1\,F(w + e_1 t_1),
\end{aligned}
\end{equation*}
with all right-hand sides finite, and either $e_1 t_1 + e_2 t_2 + e_3 t_3 + e_4 t_4 \in \Lambda$, or
$2w + e_1 t_1 - e_2 t_2 - e_3 t_3 - e_4 t_4 \in \Gamma$.
\end{lemma}

\begin{proof}
First, by~\eqref{eq:N4} and Lemmas~2 and~4 of~\cite{NURMATOV2026104102},
$r_i$, $s_i$, $f_i$, $M \neq 0, 1$ are well-defined reals and $f_i > 0$.
Hence $p_i$, $q_i \in \mathbb{R}_{>0} \cup \mathbf{i}\mathbb{R}_{>0}$,
see~\eqref{eq:N10}, $k \in (0, 1)$, because $M > 0$ and $M \neq 1$,
see~\eqref{eq:N7}, and the phase shift $t_i$ is well-defined,
see~\eqref{eq:N11} and Table~\ref{tab:tpqsigma}. Moreover, by~\eqref{eq:N5}
and~\eqref{eq:N10}, we have $p_2 = p_1$, $p_4 = p_3$, $q_3 = q_2$, and
$q_4 = q_1$.

($\Rightarrow$) Let $(Z, W_2, U, W_1) \in \mathbb{C}^4$
solve~\eqref{eq:N13}. By Lemma~\ref{lem:main:1} with $e = 1$, there are
$v_1$, $v_2$, $v_3$, $v_4 \in \mathbb{C}$ such that
$(Z, W_1) = \bigl(p_1 F(v_1), q_1 F(v_1 + t_1)\bigr)$,
$(Z, W_2) = \bigl(p_2 F(v_2), q_2 F(v_2 + t_2)\bigr)$,
$(U, W_2) = \bigl(p_3 F(v_3), q_3 F(v_3 + t_3)\bigr)$,
and $(U, W_1) = \bigl(p_4 F(v_4), q_4 F(v_4 + t_4)\bigr)$.
Set $w \coloneqq v_1$ and $e_1 \coloneqq 1$. 
Then $Z/p_1 = F(w) = F(v_2)$, $W_2/q_2 = F(v_2 + t_2) = F(v_3 + t_3)$, 
$U/p_3 = F(v_3) = F(v_4)$, and $W_1/q_1 = F(w + e_1 t_1) = F(v_4 + t_4)$
because $p_2 = p_1$, $p_4 = p_3$, $q_3 = q_2$, $q_4 = q_1$, 
and these quantities are nonzero.
By Remark~\ref{rem:lem:main:1-2}, using $F(v_2) = F(w)$, we get 
$F(v_2 + t_2) = F(w - e_2 t_2)$ for some $e_2 \in \{\pm 1\}$; 
using $F(v_3 + t_3) = F(v_2 + t_2)= F(w - e_2 t_2)$, we get 
$F(v_3) = F(w - e_2 t_2 - e_3 t_3)$ for some $e_3 \in \{\pm 1\}$;
using $F(v_4) =F(v_3) = F(w - e_2 t_2 - e_3 t_3)$, we get 
$F(v_4 + t_4)  = F(w - e_2 t_2 - e_3 t_3 - e_4 t_4)$ for some
$e_4 \in \{\pm 1\}$. 
This proves the formulas of the lemma, all right-hand sides being finite
because $Z$, $W_2$, $U$, $W_1 \in \mathbb{C}$.
Finally, since
$q_1 F(w + e_1 t_1) = W_1 = q_1 F(w - e_2 t_2 - e_3 t_3 - e_4 t_4)$ and
$q_1 \neq 0$, it follows, by Lemma~\ref{lem:main:1-2}, that
$e_1 t_1 + e_2 t_2 + e_3 t_3 + e_4 t_4 \in \Lambda$, or
$2w + e_1 t_1 - e_2 t_2 - e_3 t_3 - e_4 t_4 \in \Gamma$.

($\Leftarrow$) Assume that $(Z, W_2, U, W_1)$ is given by the formulas of
the lemma. The first three equations of~\eqref{eq:N13} hold by
Lemma~\ref{lem:main:1}: $P_1(Z, W_1) = 0$ with $i = 1$, $v = w$, $e = e_1$;
$P_2(Z, W_2) = 0$ with $i = 2$, $v = w$, $e = -e_2$, using $p_2 = p_1$;
$P_3(U, W_2) = 0$ with $i = 3$, $v = w - e_2 t_2 - e_3 t_3$, $e = e_3$,
using $q_3 = q_2$. Moreover, since either $e_1 t_1 + e_2 t_2 + e_3 t_3 + e_4 t_4 \in \Lambda$, or
$2w + e_1 t_1 - e_2 t_2 - e_3 t_3 - e_4 t_4 \in \Gamma$, it follows,
by Lemma~\ref{lem:main:1-2} applied with $v' = w + e_1 t_1$ and
$v = w - e_2 t_2 - e_3 t_3 - e_4 t_4$,
that $F(w + e_1 t_1) = F(w - e_2 t_2 - e_3 t_3 - e_4 t_4)$, so that
$(U, W_1) = \bigl(p_4 F(w - e_2 t_2 - e_3 t_3),
q_4 F(w - e_2 t_2 - e_3 t_3 - e_4 t_4)\bigr)$, because $p_4 = p_3$ and
$q_4 = q_1$, and $P_4(U, W_1) = 0$ holds by Lemma~\ref{lem:main:1} with
$i = 4$, $v = w - e_2 t_2 - e_3 t_3$, $e = -e_4$.
\end{proof}

\begin{remark}\label{rem:lem:main:2}
Under the hypotheses of Lemma~\ref{lem:main:2}, 
fix $e_1$, $e_2$, $e_3$, $e_4 \in \{\pm 1\}$. If they satisfy
condition~\eqref{eq:N6}, which is the first alternative of
Lemma~\ref{lem:main:2}, then the formulas of Lemma~\ref{lem:main:2} give a
solution of~\eqref{eq:N13} for every $w \in \mathbb{C}$ at which their
right-hand sides are finite; we call these solutions the family of $e_1$,
$e_2$, $e_3$, $e_4$. 
Otherwise, only the second alternative can hold, and it determines $2w$
modulo $\Lambda$, because $\Gamma$ is a coset of $\Lambda$. Hence $w$
takes at most four values modulo $\Lambda$, because
$\tfrac{1}{2}\Lambda/\Lambda$ has order four, and these values give at
most four solutions by~\eqref{eq:P1}.
Thus all but finitely many
solutions of~\eqref{eq:N13} lie on the families. 
A solution obtained from the second alternative may still lie on a
family, because $e_4$ does not occur in the formulas of
Lemma~\ref{lem:main:2}; for instance, it lies on the family of $e_1$,
$e_2$, $e_3$, $-e_4$ if these satisfy condition~\eqref{eq:N6}.
Each family is closed in $\mathbb{C}^4$, because the right-hand sides of
the formulas of Lemma~\ref{lem:main:2} depend on $w$ only modulo
$\Lambda$ and each of them tends to $\infty$ as $w$ tends to one of its
poles. Hence each solution lying on none of the families has a
neighbourhood in $\mathbb{C}^4$ containing no other solution
of~\eqref{eq:N13}.
Now consider a flexion of a polyhedron with flat angles $\alpha_i$,
$\beta_i$, $\gamma_i$, $\delta_i$, with parameter $t \in [0, \epsilon)$.
While the central dihedral angles $\theta_1, \ldots, \theta_4$ of its
polyhedra are nonzero, the solution $(Z, W_2, U, W_1) = \bigl(\cot{\frac{\theta_1}{2}}, \,
\cot{\frac{\theta_2}{2}}, \, \cot{\frac{\theta_3}{2}}, \,
\cot{\frac{\theta_4}{2}}\bigr)$ of~\eqref{eq:N13}, see Lemma~\ref{lem:bricard},
depends continuously on $t$ and is not constant on any subinterval of $[0, \epsilon)$,
because the polyhedra of a flexion have congruent faces and are pairwise
non-congruent, and a polyhedron is determined up to congruence by its
faces and its central dihedral angles. Hence this solution always lies on
one of the families, because a continuous solution reaching a solution
lying on none of them would stay at it on a subinterval. In particular, a
polyhedron with flat angles $\alpha_i$, $\beta_i$, $\gamma_i$,
$\delta_i$ whose solution $(Z, W_2, U, W_1)$ lies on none of the families
is rigid; an example is the second polyhedron in Proposition~3
of~\cite{NURMATOV2026104102}, see its proof.
\end{remark}

Setting $Z = t$ in the formulas of Lemma~\ref{lem:main:2},
with $e_1$, $e_2$, $e_3$, $e_4$ satisfying condition~\eqref{eq:N6},
we show next that, for the flat angles of every polyhedron of equimodular
elliptic type, there is a nonempty real interval of values of
$t$ on which $W_2$, $U$, and $W_1$ are real, so that, together
with $Z = t$, they form a
one-parameter family of real solutions of Bricard's
system~\eqref{eq:N13}, see Lemma~\ref{lem:main:7}.
To this end, we expand the right-hand
sides of the formulas of Lemma~\ref{lem:main:2}, which requires
the addition theorems~\cite[\S\S 22.2 and
22.21]{WhittakerWatson2021}:
\begin{equation}\label{eq:P2}
\begin{aligned}
    \operatorname{sn}(v + w) &= \frac{\operatorname{sn}{v}\operatorname{cn}{w}\operatorname{dn}{w} + \operatorname{sn}{w}\operatorname{cn}{v}\operatorname{dn}{v}}{1 - k^2\operatorname{sn}^2{v}\operatorname{sn}^2{w}},
    \\
    \operatorname{cn}(v + w) &= \frac{\operatorname{cn}{v}\operatorname{cn}{w} - \operatorname{sn}{v}\operatorname{sn}{w}\operatorname{dn}{v}\operatorname{dn}{w}}{1 - k^2\operatorname{sn}^2{v}\operatorname{sn}^2{w}},
    \\
    \operatorname{dn}(v + w) &= \frac{\operatorname{dn}{v}\operatorname{dn}{w} - k^2\operatorname{sn}{v}\operatorname{sn}{w}\operatorname{cn}{v}\operatorname{cn}{w}}{1 - k^2\operatorname{sn}^2{v}\operatorname{sn}^2{w}},
\end{aligned}
\end{equation}
identities of meromorphic functions of $v$ and $w$. In
particular, they hold pointwise whenever $v$ and $w$ are not
poles and $k^2\operatorname{sn}^2{v}\operatorname{sn}^2{w} \neq
1$. We also repeatedly use the \emph{Pythagorean relations}
$\operatorname{sn}^2{v} + \operatorname{cn}^2{v} = 1$ and
$\operatorname{dn}^2{v} + k^2\operatorname{sn}^2{v} = 1$, valid
for every $v \in \mathbb{C}$ that is not a pole,
see~\cite[\S 22.11]{WhittakerWatson2021}.
The right-hand sides of the formulas of Lemma~\ref{lem:main:2}
for $W_2$ and $W_1$ can be expanded directly, and we present
their expansions together with that of $U$. Expanding the
right-hand side of the formula for $U$ requires extra work:
Lemmas~\ref{lem:main:3}--\ref{lem:main:6} below serve this
purpose.

In Lemmas~\ref{lem:main:3}--\ref{lem:main:5}, we assume
that $k \in (0, 1)$ and
$t_i \in (\mu_i K, \mu_i K + \mathbf{i}K')$ with
$\mu_i \in \{0, 1, 2, 3\}$ for $i \in \{1, 2, 3, 4\}$; the $\mu_i$ are
uniquely determined by the $t_i$, because the intervals
$(0, \mathbf{i}K')$, $(K, K + \mathbf{i}K')$,
$(2K, 2K + \mathbf{i}K')$, and $(3K, 3K + \mathbf{i}K')$ are pairwise
disjoint. Note that the phase shifts of flat angles satisfying~\eqref{eq:N1}
and~\eqref{eq:N4} with $M > 0$ fulfill this assumption:
$k \in (0, 1)$, because $M \neq 1$ by~\eqref{eq:N4} and Lemma~2
of~\cite{NURMATOV2026104102}, see~\eqref{eq:N7}, and
$t_i \in (\mu_i K, \mu_i K + \mathbf{i}K')$ by~\eqref{eq:N11}
and Table~\ref{tab:tpqsigma}.
Lemmas~\ref{lem:main:3}--\ref{lem:main:5} concern an arbitrary pair of the
points $t_1$, $t_2$, $t_3$, $t_4$, whose indices need not be distinct.

\begin{lemma}\label{lem:main:3}
Let $i$, $j \in \{1, 2, 3, 4\}$, and let $e$, $e' \in \{\pm 1\}$. Then
$e t_i + e' t_j$ is a pole if and only if $e = e'$ and $t_i + t_j$ is a
pole.
\end{lemma}

\begin{proof}
First, by assumption, $k \in (0, 1)$ and  $\operatorname{Im}{t_i}$, $\operatorname{Im}{t_j} \in (0, K')$.
Next, by the definition of
$\Lambda'$, each element of
$\mathbf{i}K' + \Lambda'$ has imaginary part an odd integer multiple of $K'$. 
Since $-(\mathbf{i}K' + \Lambda') = \mathbf{i}K' + \Lambda'$ (because 
$\Lambda' = -\Lambda'$ and $2\mathbf{i}K' \in \Lambda'$), 
and since
$t_i + e e' t_j = e(e t_i + e' t_j)$ and $e \in \{\pm 1\}$, it follows
that $e t_i + e' t_j \in \mathbf{i}K' + \Lambda'$ if and only if
$t_i + e e' t_j \in \mathbf{i}K' + \Lambda'$. For $e \neq e'$ the latter
membership reads $t_i - t_j \in \mathbf{i}K' + \Lambda'$ and fails,
because $\operatorname{Im}(t_i - t_j) = \operatorname{Im}{t_i} - \operatorname{Im}{t_j} \in (-K', K')$, an
interval containing no odd integer multiple of $K'$. Thus,
$e t_i + e' t_j \in \mathbf{i}K' + \Lambda'$ if and only if $e = e'$ and
$t_i + t_j \in \mathbf{i}K' + \Lambda'$.
\end{proof}

We further use the quarter-period
shifts (see~\cite[\S 22.34]{WhittakerWatson2021})
\begin{equation}\label{eq:P3}
    (\operatorname{sn}, \operatorname{cn}, \operatorname{dn}){(v \pm \mathbf{i}K')} 
    = \biggl(\frac{1}{k \operatorname{sn}{v}}, \, \mp\frac{\mathbf{i}\operatorname{dn}{v}}{k \operatorname{sn}{v}},
    \, \mp \frac{\mathbf{i}\operatorname{cn}{v}}{\operatorname{sn}{v}}\biggr),
\end{equation}
identities of meromorphic functions of $v$; 
and the zero sets (see~\cite[\S 22.35]{WhittakerWatson2021})
\begin{equation}\label{eq:P4}
    (\operatorname{sn}^{-1}, \operatorname{cn}^{-1}, \operatorname{dn}^{-1}){(0)} 
    = \bigl(\Lambda', \, K + \Lambda', \, K + \mathbf{i}K' + \Lambda'\bigr),
\end{equation}
which are pairwise disjoint.

\begin{lemma}\label{lem:main:4}
Let $i$, $j \in \{1, 2, 3, 4\}$. Then $t_i + t_j$ is a
pole if and only if $k^2 \operatorname{sn}^2{t_i} \operatorname{sn}^2{t_j} = 1$.
\end{lemma}

\begin{proof}
First, by assumption, $k \in (0, 1)$ and $\operatorname{Im}{t_i}$,
$\operatorname{Im}{t_j} \in (0, K')$.
Next, observe that every zero or pole of
$\operatorname{sn}$ has imaginary part an 
integer multiple of $K'$, see~\eqref{eq:P4} and the definition of a pole.
Since every point at which $\operatorname{sn}$ is evaluated below has imaginary part lying
in $(-K', 0) \cup (0, K')$, it follows that none of these points is a
zero or a pole of $\operatorname{sn}$.

($\Rightarrow$) Suppose $t_i + t_j \in \mathbf{i}K' + \Lambda'$, that is,
$t_i + t_j = 2mK + (2n+1)\mathbf{i}K'$ for some $m$, $n \in \mathbb{Z}$.
Then $t_j = 2mK + 2n\mathbf{i}K' + (\mathbf{i}K' - t_i)$, 
so~\eqref{eq:P1} and~\eqref{eq:P3} give
\begin{equation*}
    \operatorname{sn}{t_j}
    = (-1)^{m}\operatorname{sn}(\mathbf{i}K' - t_i)
    = \frac{(-1)^{m}}{k\operatorname{sn}(-t_i)}
    = \frac{(-1)^{m+1}}{k\operatorname{sn}{t_i}}.
\end{equation*}
Thus $k\operatorname{sn}{t_i}\operatorname{sn}{t_j} = (-1)^{m+1}$. Squaring gives
$k^2\operatorname{sn}^2{t_i}\operatorname{sn}^2{t_j} = 1$.

($\Leftarrow$) Assume $k^2\operatorname{sn}^2{t_i}\operatorname{sn}^2{t_j} = 1$. Since
$\operatorname{sn}{t_i} \neq 0$, it follows, using~\eqref{eq:P3}, that
$\operatorname{sn}^2{t_j} = 1/(k^2\operatorname{sn}^2{t_i})
= \operatorname{sn}^2(t_i - \mathbf{i}K')$.
The function $\operatorname{sn}^2$ is, by~\eqref{eq:P1}, even and elliptic
with respect to the lattice $\Lambda'$. Since $\operatorname{sn}$ has a
single simple pole modulo $\Lambda'$, located at
$\mathbf{i}K'$~\cite[\S 22.35]{WhittakerWatson2021}, it follows that
$\operatorname{sn}^2$ has one double pole modulo $\Lambda'$ and is
therefore elliptic of order two, so that $\operatorname{sn}^2$ attains every
value exactly twice modulo $\Lambda'$, counting multiplicity.
By evenness, $\operatorname{sn}^2$
attains the value $\operatorname{sn}^2(t_i - \mathbf{i}K')$ at both
$t_i - \mathbf{i}K'$ and $\mathbf{i}K' - t_i$, and these two points are
distinct modulo $\Lambda'$. Indeed, otherwise
$2(t_i - \mathbf{i}K') \in \Lambda'$, so that
$\operatorname{Im}{(t_i - \mathbf{i}K')} = \operatorname{Im}{t_i} - K'$
would lie in $K'\mathbb{Z}$, contradicting
$\operatorname{Im}{t_i} \in (0, K')$. Hence, modulo $\Lambda'$, the
value $\operatorname{sn}^2(t_i - \mathbf{i}K')$ is attained only at
$t_i - \mathbf{i}K'$ and $\mathbf{i}K' - t_i$, and since
$\operatorname{sn}^2{t_j} = \operatorname{sn}^2(t_i - \mathbf{i}K')$, it follows that
$t_j \equiv \pm(t_i - \mathbf{i}K') \pmod{\Lambda'}$.
The branch $t_j \equiv t_i - \mathbf{i}K' \pmod{\Lambda'}$ reads
$t_i - t_j \in \mathbf{i}K' + \Lambda'$, that is, $t_i - t_j$ is a pole,
which is ruled out by Lemma~\ref{lem:main:3}, applied with
$e = 1$ and $e' = -1$. 
Thus $t_j \equiv \mathbf{i}K' - t_i \pmod{\Lambda'}$, that is,
$t_i + t_j \in \mathbf{i}K' + \Lambda'$.
\end{proof}

\begin{remark}\label{rem:lem:main:4}
It is shown in the proof of Lemma~\ref{lem:main:4} that, for $i$, $j \in \{1, 2, 3, 4\}$, if
$t_i + t_j = 2mK + (2n + 1)\mathbf{i}K'$ with $m$, $n \in \mathbb{Z}$,
then $k \operatorname{sn}{t_i} \operatorname{sn}{t_j} = (-1)^{m + 1}$. 
\end{remark}

Hereafter, for $i \in \{1, 2, 3, 4\}$, we denote $g_i \coloneqq \operatorname{dn}{t_i}$.

\begin{lemma}\label{lem:main:5}
Let $i$, $j \in \{1, 2, 3, 4\}$, and assume that
$t_i + t_j$ is a pole. Then $t_i + t_j = 2mK + \mathbf{i}K'$ for some
$m \in \mathbb{Z}$ with $(-1)^{m} = -k\operatorname{sn}{t_i}\operatorname{sn}{t_j}$, and
\begin{equation*}
    (\operatorname{sn}, \operatorname{cn}, \operatorname{dn}){(t_i - t_j)} =
    \biggl(\mathbf{i}\,(-1)^{m}\,\frac{g_i^2 - g_j^2}{2 k\, g_i g_j}, 
    \, (-1)^{m}\,\frac{g_i^2 g_j^2 - k'^2}{2 k\, g_i g_j},
    \, \frac{g_i^2 + g_j^2}{2\, g_i g_j}\biggr).
\end{equation*}
\end{lemma}

\begin{proof}
First, by assumption, $k \in (0, 1)$ and $\operatorname{Im}{t_i}$,
$\operatorname{Im}{t_j} \in (0, K')$.
Since $t_i + t_j$ is a pole and
$\operatorname{Im}(t_i + t_j)  \in (0, 2K')$, an interval
containing no odd integer multiple of $K'$ other than $K'$, it follows
that $t_i + t_j = 2mK + \mathbf{i}K'$ for some $m \in \mathbb{Z}$, and
$k\operatorname{sn}{t_i}\operatorname{sn}{t_j} = (-1)^{m + 1}$ by Remark~\ref{rem:lem:main:4}. 
Moreover, $t_i$ and $t_j$ are not poles and, by~\eqref{eq:P4}, none of
$\operatorname{sn}$, $\operatorname{cn}$, $\operatorname{dn}$
vanishes at $t_i$ or $t_j$, because
$\operatorname{Im}{t_i}$, $\operatorname{Im}{t_j} \in (0, K')$.
Also, $t_j = -(t_i - \mathbf{i}K') + 2mK$, so~\eqref{eq:P1}
and~\eqref{eq:P3} give
$g_j = \operatorname{dn}(t_i - \mathbf{i}K')
= \mathbf{i}\,\operatorname{cn}{t_i}/\operatorname{sn}{t_i}$, that is,
$\operatorname{cn}{t_i} = -\mathbf{i}\,g_j\operatorname{sn}{t_i}$.

Next, observe that $1 - k^2\operatorname{sn}^4{t_i} = 0$ if and only if
$t_i = \mu_i K + \tfrac{1}{2}\mathbf{i}K'$. Indeed, by Lemma~\ref{lem:main:4},
applied with $j = i$, we get $1 - k^2\operatorname{sn}^4{t_i} = 0$ if and only if
$2t_i$ is a pole. 
Since $\operatorname{Re}{(2t_i)} = 2\mu_i K$ and
$\operatorname{Im}{(2t_i)} \in (0, 2K')$, an interval containing no
odd integer multiple of $K'$ other than $K'$, it follows that $2t_i$
is a pole if and only if $2t_i = 2\mu_i K + \mathbf{i}K'$, that is,
if and only if $t_i = \mu_i K + \tfrac{1}{2}\mathbf{i}K'$.

Assume first that $t_i \neq \mu_i K + \tfrac{1}{2}\mathbf{i}K'$. Then
$2t_i$ is not a pole and $1 - k^2\operatorname{sn}^4{t_i} \neq 0$.
Moreover, $\operatorname{sn}{(2t_i)} \neq 0$ by~\eqref{eq:P4}, because no element of
$\Lambda'$ has imaginary part in $(0, 2K')$.
Since $t_i$ and $2t_i$ are not poles and
$\operatorname{sn}{(2t_i)}$, $1 - k^2\operatorname{sn}^4{t_i} \neq 0$,
the formulas stated in the lemma follow from
$t_i - t_j = 2t_i - (t_i + t_j) = -2mK + (2t_i - \mathbf{i}K')$,
\eqref{eq:P1}, \eqref{eq:P2} applied with $v = w = t_i$, \eqref{eq:P3}, 
$\operatorname{cn}{t_i} = -\mathbf{i}\, g_j\operatorname{sn}{t_i}$,
$k'^2 = 1 - k^2$, see~\eqref{eq:N7}, and the Pythagorean relations, as verified in
Section~5 of the supplement (\texttt{helper.nb})~\cite{Nurmatov2026Repo}.

Assume now that $t_i = \mu_i K + \tfrac{1}{2}\mathbf{i}K'$. 
By Remark~\ref{rem:lem:main:4}, applied
with $j = i$, we get $k\operatorname{sn}^2{t_i} = (-1)^{\mu_i + 1}$, whence,
by the Pythagorean relations,
$g_i^2 = 1 - k^2\operatorname{sn}^2{t_i}
= 1 + (-1)^{\mu_i} k$. Moreover,
$t_i + \mathbf{i}K' = 2\mu_i K + 2\mathbf{i}K' - t_i$, so~\eqref{eq:P1}
and~\eqref{eq:P3} give
$g_i = -\operatorname{dn}(t_i + \mathbf{i}K')
= \mathbf{i}\,\operatorname{cn}{t_i}/\operatorname{sn}{t_i} = g_j$.
By $g_j = g_i$,
$g_i^2 = 1 + (-1)^{\mu_i} k$, and
$k'^2 = 1 - k^2$, see~\eqref{eq:N7}, the right-hand
sides of the formulas stated in the lemma equal $0$,
$(-1)^{\mu_i - m}$, $1$, respectively.
On the other hand,
$t_j = (t_i + t_j) - t_i = (2m - \mu_i)K + \tfrac{1}{2}\mathbf{i}K'$,
whence $t_i - t_j = 2(\mu_i - m)K$. Then, by~\eqref{eq:P1} and 
$(\operatorname{sn}, \operatorname{cn}, \operatorname{dn})(0) = (0, 1, 1)$,
see~\cite[\S 22.11]{WhittakerWatson2021}, we have
$(\operatorname{sn}, \operatorname{cn}, \operatorname{dn}){(t_i - t_j)} = (0, \,(-1)^{\mu_i - m}, \,1)$,
in agreement with the right-hand sides.
\end{proof}

Lemma~\ref{lem:main:6} below is the last piece needed to expand the right-hand
side of the formula for $U$, see Lemma~\ref{lem:main:2}. 
Hereafter, for $i$, $j \in \{1, 2, 3, 4\}$, we introduce the notation
\begin{equation}\label{eq:P5}
    \xi_{ij} \coloneqq \sign{(z_i z_j)},
    \qquad
    S_{ij} \coloneqq \frac{\varepsilon_i \varepsilon_j \sqrt{p_i^2 q_i^2\, p_j^2 q_j^2\, \xi_{ij}}}
    {p_i q_i\, p_j q_j},
\end{equation}
and, for $i \neq j$, $e_i$, $e_j \in \{\pm 1\}$, in the case
$k^2 \neq (g_i^2 - 1)(g_j^2 - 1)$, we set
\begin{equation}\label{eq:P6}
\begin{aligned}
    h_{ij1} &\coloneqq \frac{e_i \, g_j \sqrt{(g_i^2 - 1)(g_j^2 - k'^2)\, \xi_{ij}}
    + e_j\, g_i \sqrt{(g_j^2 - 1)(g_i^2 - k'^2)\, \xi_{ij}}}
    {k^2 - (g_i^2 - 1)(g_j^2 - 1)},
    \\
    h_{ij2} &\coloneqq \frac{\sqrt{(g_i^2 - k'^2)(g_j^2 - k'^2)\, \xi_{ij}}
    + e_i \, e_j \, g_i g_j \sqrt{(g_i^2 - 1)(g_j^2 - 1)\, \xi_{ij}}}
    {k^2 - (g_i^2 - 1)(g_j^2 - 1)},
    \\
    h_{ij3} &\coloneqq \frac{k^2 g_i g_j
    + e_i \, e_j \, \xi_{ij}
    \sqrt{(g_i^2 - 1)(g_i^2 - k'^2)(g_j^2 - 1)(g_j^2 - k'^2)}}
    {k^2 - (g_i^2 - 1)(g_j^2 - 1)},
\end{aligned}
\end{equation}
while, in the case $k^2 = (g_i^2 - 1)(g_j^2 - 1)$ and $e_i = -e_j$, we set
\begin{equation}\label{eq:P7}
    h_{ij1} \coloneqq \frac{e_i\,(g_i^2 - g_j^2)}{2 k\, g_i g_j},
    \qquad
    h_{ij2} \coloneqq \frac{g_i^2 g_j^2 - k'^2}{2 k\, g_i g_j},
    \qquad
    h_{ij3} \coloneqq \frac{g_i^2 + g_j^2}{2\, g_i g_j}.
\end{equation}
Note that, in the latter case, the common denominator
$k^2 - (g_i^2 - 1)(g_j^2 - 1)$ of the fractions in~\eqref{eq:P6}
vanishes, which is why the quantities $h_{ij1}$, $h_{ij2}$, $h_{ij3}$
are defined separately in this case.

\begin{lemma}\label{lem:main:6}
Assume that $\alpha_i$, $\beta_i$, $\gamma_i$, $\delta_i \in (0, \pi)$
satisfy conditions~\eqref{eq:N1} and~\eqref{eq:N4} for each
$i = 1, \ldots, 4$, and $M > 0$. Let $i$, $j \in \{1, 2, 3, 4\}$
with $i \neq j$. 
Then, for all $e_i$, $e_j \in \{\pm 1\}$ such that
$e_i t_i + e_j t_j$ is not a pole,
\begin{equation*}
    (\operatorname{sn}, \operatorname{cn}, \operatorname{dn}){(e_i t_i + e_j t_j)} 
    = (\mathbf{i}\, S_{ij} h_{ij1}, \,S_{ij} h_{ij2}, \,h_{ij3}),
\end{equation*}
where $h_{ij1}$, $h_{ij2}$, $h_{ij3}$ are given by~\eqref{eq:P6} if
$k^2 \neq (g_i^2 - 1)(g_j^2 - 1)$ and by~\eqref{eq:P7} otherwise.
Moreover, all radicands in~\eqref{eq:P5} and~\eqref{eq:P6} are
strictly positive reals.
\end{lemma}

\begin{proof}
First, $r_1, \ldots, r_4$, $s_1, \ldots, s_4$,
$f_1, \ldots, f_4$, $M \neq 0, 1$ are well-defined reals and
$f_1, \ldots, f_4 > 0$, by~\eqref{eq:N4} and Lemmas~2 and~4
of~\cite{NURMATOV2026104102}, and
$\varepsilon_1, \ldots, \varepsilon_4 \in \{\pm 1\}$, by~\eqref{eq:N1}
and~\eqref{eq:N2}.
Hence $p_1, \ldots, p_4$, $q_1, \ldots, q_4 \in \mathbb{R}_{>0} \cup
\mathbf{i}\mathbb{R}_{>0}$, see~\eqref{eq:N10}, and $k$, $k' \in (0, 1)$,
because $M > 0$ and $M \neq 1$, see~\eqref{eq:N7}.
Then, by~\eqref{eq:N7}, \eqref{eq:N11}, and
Lemma~4 of~\cite{NURMATOV2026104102}, we get
$g_1, \ldots, g_4 \in (0, k') \cup (1, +\infty)$ and,
by~\eqref{eq:N11} and Table~\ref{tab:tpqsigma},
$t_1 \in (\mu_1 K, \mu_1 K + \mathbf{i}K')$, \ldots,
$t_4 \in (\mu_4 K, \mu_4 K + \mathbf{i}K')$ with
$\mu_1, \ldots, \mu_4 \in \{0, 1, 2, 3\}$.
Moreover,
$\sign{(g_i^2 - 1)} = \sign{(g_i^2 - k'^2)} = \sign{(p_i^2 q_i^2)}\neq 0$
and
$\sign{(g_j^2 - 1)} = \sign{(g_j^2 - k'^2)} = \sign{(p_j^2 q_j^2)}\neq 0$, 
by Lemma~\ref{lem:shifts:2}.

Second, observe that $\xi_{ij} = \sign{\bigl(p_i^2 q_i^2\, p_j^2 q_j^2\bigr)} \in \{\pm 1\}$
by~\eqref{eq:N3}, \eqref{eq:N4}, \eqref{eq:N10}, 
and Lemma~4 of~\cite{NURMATOV2026104102}.
Hence any product of an element of
$\{g_i^2 - 1, g_i^2 - k'^2, p_i^2 q_i^2\}$ with an element of
$\{g_j^2 - 1, g_j^2 - k'^2, p_j^2 q_j^2\}$ has sign $\xi_{ij}$. 
Thus, all radicands in~\eqref{eq:P5} and~\eqref{eq:P6} are strictly positive.

Assume first that $k^2 \neq (g_i^2 - 1)(g_j^2 - 1)$. Then, by 
the Pythagorean relations, 
$k^2\operatorname{sn}^2{t_i}\operatorname{sn}^2{t_j}
= \frac{(g_i^2 - 1)(g_j^2 - 1)}{k^2} \neq 1$. 
Next, by~\eqref{eq:P1},
$(\operatorname{sn}, \operatorname{cn}, \operatorname{dn}){(e_i t_i)}
= (e_i \operatorname{sn}{t_i}, \,\operatorname{cn}{t_i}, \,g_i)$
and
$(\operatorname{sn}, \operatorname{cn}, \operatorname{dn}){(e_j t_j)}
= (e_j \operatorname{sn}{t_j}, \,\operatorname{cn}{t_j}, \,g_j)$,
and $\pm t_i$, $\pm t_j$ are not poles, because their imaginary parts
lie in $(-K', 0) \cup (0, K')$.
Thus, by Lemma~\ref{lem:shifts:2} and~\eqref{eq:P2}, applied with
$v = e_i t_i$ and $w = e_j t_j$, the formulas stated in the lemma
follow, as verified in Section~6 of the supplement
(\texttt{helper.nb})~\cite{Nurmatov2026Repo}.

Assume now that $k^2 = (g_i^2 - 1)(g_j^2 - 1)$. Then, by the Pythagorean
relations, $k^2\operatorname{sn}^2{t_i}\operatorname{sn}^2{t_j}
= \frac{(g_i^2 - 1)(g_j^2 - 1)}{k^2} = 1$, so, by Lemma~\ref{lem:main:4},
the point $t_i + t_j$ is a pole. Hence $e_i = -e_j$ by
Lemma~\ref{lem:main:3}, because $e_i t_i + e_j t_j$ is not a pole, and,
by Lemma~\ref{lem:main:5}, $t_i + t_j = 2mK + \mathbf{i}K'$ for some
$m \in \mathbb{Z}$ with $(-1)^{m} = -k\operatorname{sn}{t_i}\operatorname{sn}{t_j}$.
Moreover,
$\xi_{ij}
= \sign{\bigl(p_i^2 q_i^2\, p_j^2 q_j^2\bigr)}
= \sign{\bigl((g_i^2 - 1)(g_j^2 - 1)\bigr)}
= \sign{(k^2)} = 1$. Hence, by the formulas for $\operatorname{sn}{t_i}$
and $\operatorname{sn}{t_j}$, see Lemma~\ref{lem:shifts:2}, we get
$(-1)^{m} = S_{ij}$, as verified in Section~6 of
the supplement (\texttt{helper.nb})~\cite{Nurmatov2026Repo}.
Finally, $e_i t_i + e_j t_j = e_i(t_i - t_j)$ because $e_j = -e_i$. 
Then, by~\eqref{eq:P1},
$(\operatorname{sn}, \operatorname{cn}, \operatorname{dn})(e_i t_i + e_j t_j)
= \bigl(e_i\operatorname{sn}(t_i - t_j),\, \operatorname{cn}(t_i - t_j),\, \operatorname{dn}(t_i - t_j)\bigr)$.
Thus, by Lemma~\ref{lem:main:5} and $(-1)^{m} = S_{ij}$, 
the formulas stated in the lemma follow.
\end{proof}

\begin{remark}\label{rem:lem:main:6:1}
Under the hypotheses of Lemma~\ref{lem:main:6}, by the Pythagorean
relations and Lemma~\ref{lem:main:4}, $k^2 = (g_i^2 - 1)(g_j^2 - 1)$
if and only if $t_i + t_j$ is a pole. 
Moreover, if $t_i + t_j$ is a pole, then,
by Lemma~\ref{lem:main:5}, $t_i + t_j = 2mK + \mathbf{i}K'$ with
$m \in \mathbb{Z}$ and $S_{ij} = (-1)^{m}$, 
as shown in the proof of
Lemma~\ref{lem:main:6}.
\end{remark}

\begin{remark}\label{rem:lem:main:6:2}
Under the hypotheses of Lemma~\ref{lem:main:6} and for $e_i$, $e_j$ as there,
$h_{ij1}$, $h_{ij2}$, $h_{ij3}$ are real, because
$e_i$, $e_j$, $g_i$, $g_j$, $k$, $k'$, $\xi_{ij}$ are real and all
the radicands in~\eqref{eq:P6} are strictly positive, see Lemma~\ref{lem:main:6}. 
Moreover, at most one of them vanishes because if two of them vanished, then, by
Lemma~\ref{lem:main:6}, two of the values
$(\operatorname{sn}, \operatorname{cn}, \operatorname{dn}){(e_i t_i + e_j t_j)}$
would vanish, contradicting the Pythagorean relations and $k \in (0, 1)$. 
Finally, since $k'^2 = 1 - k^2$, see~\eqref{eq:N7},
$g_i$, $g_j \in (0, k') \cup (1, +\infty)$, and
$\xi_{ij} = \sign{\bigl((g_i^2 - 1)(g_j^2 - 1)\bigr)}$, see the proof of
Lemma~\ref{lem:main:6}, it follows, as verified in Section~7 of the
supplement (\texttt{helper.nb})~\cite{Nurmatov2026Repo}, that
$h_{ij1} = 0$ if and only if $g_i = g_j$ and $e_i = -e_j$;
$h_{ij2} = 0$ if and only if $g_i g_j = k'$ and $e_i = -e_j$; and
$h_{ij3} = 0$ if and only if
$k^2(g_i^2 + g_j^2) = k^2 - (g_i^2 - 1)(g_j^2 - 1)$ and
$e_i = e_j$; in particular, $h_{ij3} \neq 0$ if
$k^2 = (g_i^2 - 1)(g_j^2 - 1)$. 
Note that, by~\eqref{eq:N7} and~\eqref{eq:N11}, these conditions are
precisely the relations between $f_i$, $f_j$, and $M$ from
Definition~\ref{def:gqs}: $g_i = g_j$ reads $f_i = f_j$;
$g_i g_j = k'$ reads $f_i f_j = M$;
$k^2(g_i^2 + g_j^2) = k^2 - (g_i^2 - 1)(g_j^2 - 1)$ reads
$(f_i + f_j - 1)/(f_i f_j) = 1/M$ if $M < 1$ and
$f_i + f_j - f_i f_j = M$ if $M > 1$; and
$k^2 = (g_i^2 - 1)(g_j^2 - 1)$ reads
$f_i + f_j - f_i f_j = M$ if $M < 1$ and
$(f_i + f_j - 1)/(f_i f_j) = 1/M$ if $M > 1$.
\end{remark}

\begin{lemma}\label{lem:main:7}
Assume that $\alpha_i$, $\beta_i$, $\gamma_i$,
$\delta_i \in (0, \pi)$ satisfy conditions~\eqref{eq:N1},
\eqref{eq:N4}, and~\eqref{eq:N5} for each
$i = 1, \ldots, 4$, and $M > 0$.
Assume also that $e_1$, $e_2$,
$e_3$, $e_4 \in \{\pm 1\}$ satisfy condition~\eqref{eq:N6}.
Then, for every $t \in I(x_1)$, see~\eqref{eq:N16}, excluding
the finitely many $t$ at which any of the denominators of the
applicable formulas below vanishes (the resulting set
has nonempty interior), the formulas
\begin{equation*}
\begin{aligned}
    &Z = t,
    \\[1mm]
    &W_2
    = \frac{\varepsilon_2 \Bigl(
    t\sqrt{u x_1 y_2 (1 + z_2)(1 + u z_2)}
    + e_0\,e_2 \sqrt{u y_2 z_2 D(x_1, t)}
    \Bigr)}
    {y_2 u (z_2 + x_1 t^2)},
    \\[1mm]
    &W_1
    = \frac{\varepsilon_1 \Bigl(
    t\sqrt{u x_1 y_1 (1 + z_1)(1 + u z_1)}
    - e_0\,e_1 \sqrt{u y_1 z_1 D(x_1, t)}
    \Bigr)}
    {y_1 u (z_1 + x_1 t^2)},
    \\[1mm]
    &U
    = \frac{\varepsilon_2 \varepsilon_3\sign{(u z_2)}
    }{\sqrt{|x_1 x_3|}}
    \begin{cases}
    \dfrac{|x_1| h_{232} h_{233} \,t
    + e_0 \, h_{231} \sqrt{x_1D(x_1, t)}}
    {1 + \xi_{23}\, k^2 x_1 h_{231}^2\,t^2},
    &\text{if } u > 0,
    \\[4mm]
    \dfrac{|x_1|h_{232} \,t
    + e_0 \, k'  h_{231} h_{233} \sqrt{x_1 D(x_1, t)}}
    {h_{233}^2 - \xi_{23}\, k^2 x_1 h_{231}^2\,t^2},
    &\text{if } u < 0,
    \end{cases}
\end{aligned}
\end{equation*}
define a real solution of Bricard's system~\eqref{eq:N13}, and
these solutions form a one-parameter family over $\mathbb{R}$.

Here each of the two signs $e_0 \in \{\pm 1\}$ gives such a family,
and $h_{231}$, $h_{232}$, $h_{233}\in \mathbb{R}$ are the quantities
$h_{ij1}$, $h_{ij2}$, $h_{ij3}$ for $(i, j) = (2, 3)$, as follows:
\begin{itemize}
    \item[(a)] if $k^2 \neq (g_2^2 - 1)(g_3^2 - 1)$, then they
    are given by~\eqref{eq:P6};
    \item[(b)] if $k^2 = (g_2^2 - 1)(g_3^2 - 1)$ and
    $e_2 = -e_3$, then they are given by~\eqref{eq:P7};
    \item[(c)] if $k^2 = (g_2^2 - 1)(g_3^2 - 1)$ and
    $e_2 = e_3$, then the $U$-formula
    is replaced by
    \begin{equation*}
        U
        = \frac{\varepsilon_2 \varepsilon_3\sign{(u z_2)}}
        {\sqrt{|x_1 x_3|}}
        \begin{cases}
        \dfrac{\sign{(x_1)}}{k\,t},
        &\text{if } u > 0,
        \\[4mm]
        \dfrac{e_0\, e_2\, k' \sqrt{x_1 D(x_1, t)}}
        {k(1 - x_1t^2)},
        &\text{if } u < 0.
        \end{cases}
    \end{equation*}
\end{itemize}
Moreover, all the radicands in the formulas above are strictly positive reals.
\end{lemma}

\begin{proof}
First, fix $e_0 \in \{\pm 1\}$ and let $t$ be as in the lemma.
Then, by~\eqref{eq:N4}, $M > 0$, and Lemmas~2 and~4
of~\cite{NURMATOV2026104102}, $r_i$, $s_i$, $f_i$,
$M \neq 0, 1$ are well-defined reals and $f_i > 0$.
Hence
$x_i$, $y_i$, $z_i \neq -1, 0$ are
well-defined reals, see~\eqref{eq:N3}, with 
$z_i (1 + z_i) > 0$,
$1 + u z_i > 0$,
and $u x_i y_i z_i > 0$
by~\eqref{eq:N3}, \eqref{eq:N4}, \eqref{eq:N14}, and
Lemma~4 of~\cite{NURMATOV2026104102}.
Therefore
$u x_i y_i (1 + z_i)(1 + u z_i) > 0$, 
and $x_1 D(x_1, t)$, $u y_1 z_1 D(x_1, t)$,
$u y_2 z_2 D(x_1, t) > 0$, because
$x_2 = x_1$ and $t \in I(x_1)$, see~\eqref{eq:N5}
and~\eqref{eq:N16}.
Thus all the radicands in the statement of the lemma 
are strictly positive reals, and those in~\eqref{eq:P6} are strictly positive
by Lemma~\ref{lem:main:6}.
Throughout the proof, we write
$R \coloneqq \sqrt{x_1 D(x_1, t)}$. 
Moreover, $\sign{(u z_2)} = \sign{(x_1 y_2)}$ and
$\xi_{23} = \sign{(x_1 x_3)}$, because
$\xi_{23} = \sign{(z_2 z_3)}$,
$u x_2 y_2 z_2 > 0$, $u x_3 y_3 z_3 > 0$, $x_2 = x_1$, and
$y_3 = y_2$, see~\eqref{eq:N3}, \eqref{eq:N5}, and~\eqref{eq:P5}.

Second, $\varepsilon_i \in \{\pm 1\}$ by~\eqref{eq:N1} and~\eqref{eq:N2},
$p_i$, $q_i \in \mathbb{R}_{>0} \cup \mathbf{i}\mathbb{R}_{>0}$,
see~\eqref{eq:N10}, $k$, $k' \in (0, 1)$, because $M > 0$ and
$M \neq 1$, see~\eqref{eq:N7}, and the phase shift
$t_i \in (\mu_i K, \mu_i K + \mathbf{i}K')$
with $\mu_i \in \{0, 1, 2, 3\}$ is
well-defined, see~\eqref{eq:N11} and
Table~\ref{tab:tpqsigma}.
In addition, $x_ip_i^2 =1$ and $y_iq_i^2=1$, see~\eqref{eq:N3} and~\eqref{eq:N10}.
Further, $p_2 = p_1$, $p_4 = p_3$, $q_4 = q_1$, and
$q_3 = q_2$, see~\eqref{eq:N5} and~\eqref{eq:N10}.
Also, no $\pm t_i$ is a pole,
because
$\operatorname{Im}{(\pm t_i)} \in (-K', 0) \cup (0, K')$.

Third, $F = \operatorname{sn}$ if $u > 0$ and $F = \operatorname{cn}$
if $u < 0$, because $u = 1 - M \neq 0$, see~\eqref{eq:N14}.
Recall that the function $F$
attains every value of $\mathbb{C} \cup \{\infty\}$. 
Thus, there is
$\varsigma \in \mathbb{C}$ with
$F(\varsigma) = \frac{t}{p_1}$; 
such a $\varsigma$
is not a pole, because
$t$ is real
and $p_1 \neq 0$. 
Below, we fix such a $\varsigma$ and show that the values
$q_2 F(\varsigma - e_2 t_2)$,
$p_3 F(\varsigma - e_2 t_2 - e_3 t_3)$,
$q_1 F(\varsigma + e_1 t_1)$ equal the right-hand sides of the
formulas of the lemma for $W_2$, $U$, $W_1$, respectively;
then, by Lemma~\ref{lem:main:2}, these three values,
together with $Z = t$, form a
solution of Bricard's system~\eqref{eq:N13}.

Lastly, in cases (a) and (b), $e_2 = -e_3$ whenever
$k^2 = (g_2^2 - 1)(g_3^2 - 1)$, so $\pm (e_2 t_2 + e_3 t_3)$ is not a
pole by Lemma~\ref{lem:main:3} and Remark~\ref{rem:lem:main:6:1}, and
hence $h_{231}$, $h_{232}$, $h_{233}$ are real by
Remark~\ref{rem:lem:main:6:2}. We consider two cases:
$u > 0$ (Case 1) and $u < 0$ (Case 2).

\emph{Case 1.}
Since $u > 0$, that is, $M < 1$, it follows, by \eqref{eq:N7}
and~\eqref{eq:N11}, that $k = \sqrt{u}$, $k' = \sqrt{1-u}$, and
$g_i = \operatorname{dn}{t_i} = \sqrt{f_i}$.
Here $F = \operatorname{sn}$, so
$\operatorname{sn}{\varsigma} = \frac{t}{p_1}$ and
$\operatorname{sn}^2{\varsigma} = x_1t^2$, see~\eqref{eq:N3}
and~\eqref{eq:N10}.
Hence, by the Pythagorean
relations, $k^2 = u$, and \eqref{eq:N15}, we have
\begin{equation*}
    \bigl(p_1\operatorname{cn}{\varsigma}\operatorname{dn}{\varsigma}\bigr)^2
    = p_1^2\bigl(1 - \operatorname{sn}^2{\varsigma}\bigr)
    \bigl(1 - k^2\operatorname{sn}^2{\varsigma}\bigr)
    = -\frac{x_1D(x_1, t)}{x_1^2}
    = \left(\frac{\mathbf{i}\,R}{|x_1|}\right)^{\!2}.
\end{equation*}
By \eqref{eq:P1}, replacing $\varsigma$ by $2K - \varsigma$
preserves $\operatorname{sn}{\varsigma}$ and reverses the sign of
$\operatorname{cn}{\varsigma}\operatorname{dn}{\varsigma}$. Thus, we may fix
$p_1\operatorname{cn}{\varsigma}\operatorname{dn}{\varsigma}
    = \frac{e_0\,\mathbf{i}\,R}{|x_1|}$,
so that the two values of $e_0$ correspond to the two admissible
choices of $\varsigma$.
The computations below are verified in Section~8 of the supplement
(\texttt{helper.nb})~\cite{Nurmatov2026Repo}.

\emph{$W_2$-, $W_1$-formulas.}
Let $i \in \{1, 2\}$ and $e' \in \{\pm 1\}$. 
By the Pythagorean
relations, we have
$\operatorname{sn}^2{t_i} = \frac{1 - g_i^2}{k^2} = -\frac{1}{u z_i}$, see~\eqref{eq:N3}.
Hence $1 - k^2\operatorname{sn}^2{t_i}\operatorname{sn}^2{\varsigma}
   = \frac{z_i + x_1t^2}{z_i} \neq 0$,
the nonvanishing holding by the choice of $t$. 
Since $\pm t_i$ and $\varsigma$ are not poles, it follows that the first
identity of~\eqref{eq:P2} holds pointwise at $v = \varsigma$,
$w = e' t_i$; together with the parity identities of~\eqref{eq:P1} and
Lemma~\ref{lem:shifts:2}, it gives that
$q_i\operatorname{sn}(\varsigma + e' t_i)$ equals the $W_2$-formula of
the lemma if $i = 2$, $e' = -e_2$, and the $W_1$-formula if $i = 1$,
$e' = e_1$.

\emph{$U$-formula in cases (a) and (b).}
We have $\operatorname{sn}^2{(e_2 t_2 + e_3 t_3)} = -\xi_{23}\,h_{231}^2$, 
by Lemma~\ref{lem:main:6} and~\eqref{eq:P5}.
Hence $1 - k^2\operatorname{sn}^2{\varsigma}\operatorname{sn}^2{(e_2 t_2 + e_3 t_3)}
= 1 + \xi_{23}\,k^2 x_1 h_{231}^2\,t^2 \neq 0$, the
nonvanishing holding by the choice of $t$.
Since $\pm (e_2 t_2 + e_3 t_3)$ and $\varsigma$ are not poles, it follows
that the first identity of~\eqref{eq:P2} holds pointwise at
$v = \varsigma$, $w = -e_2 t_2 - e_3 t_3$; together with the parity
identities of~\eqref{eq:P1} and Lemma~\ref{lem:main:6}, it gives that
$p_3\operatorname{sn}(\varsigma - e_2 t_2 - e_3 t_3)$ equals the
$U$-formula of the lemma in the case $u > 0$.

\emph{$U$-formula in case (c).}
Since, in this case, $k^2 = (g_2^2 - 1)(g_3^2 - 1)$, it follows,
by Remark~\ref{rem:lem:main:6:1}, that $t_2 + t_3$ is a pole, so that
$t_2 + t_3 = 2mK + \mathbf{i}K'$ for some
$m \in \mathbb{Z}$ with $(-1)^{m}= S_{23}$.
Thus, since $e_2 = e_3$ and
$\operatorname{sn}{\varsigma} = \frac{t}{p_1} \neq 0$, the
nonvanishing holding by the choice of $t$, it follows,
by \eqref{eq:P1} and the first identity of~\eqref{eq:P3}, that
$p_3 \operatorname{sn}{(\varsigma - e_2 t_2 - e_3 t_3)}$ equals the
$U$-formula of case (c) for $u > 0$.

\emph{Case 2.}
Since $u < 0$, that is, $M > 1$, it follows, by~\eqref{eq:N7}
and~\eqref{eq:N11}, that $k' = 1/\sqrt{1 - u}$, $k = k'\sqrt{-u}$, and
$g_i = \operatorname{dn}{t_i} = 1/\sqrt{f_i}$.
Here $F = \operatorname{cn}$, so
$\operatorname{cn}{\varsigma} = \frac{t}{p_1}$.
The remaining steps are as in Case~1, with the second identities
of~\eqref{eq:P2} and~\eqref{eq:P3} in place of the first; they are
carried out in Section~9 of the supplement
(\texttt{helper.nb})~\cite{Nurmatov2026Repo}, including the choice of
$\varsigma$ according to $e_0$.

\emph{The range of $t$.}
In both cases, each displayed denominator is, up to the
nonzero constant factor $\sqrt{|x_1x_3|}$ in the
$U$-formulas, a polynomial in $t$ that is not
identically zero: the denominators $y_2 u (z_2 + x_1t^2)$ and
$y_1 u (z_1 + x_1t^2)$, because $x_1$, $y_1$, $y_2$, $u \neq 0$; the denominator
of the case $u > 0$, because its constant term equals $1$; 
the denominator of the case $u < 0$, because
$(h_{233}, h_{231}) \neq (0, 0)$ by
Remark~\ref{rem:lem:main:6:2} and $k$, $x_1$, $\xi_{23} \neq 0$; 
and the
denominators $kt$ and $k(1 - x_1t^2)$ of case (c), because
$k \neq 0$. 
Thus at most finitely many $t$ are excluded from
$I(x_1)$ in the lemma.
Since $x_1 \neq 0$, it follows,
by~\eqref{eq:N17},
that $I(x_1)$ has nonempty interior; hence
so does the set obtained from it by removing finitely many
points, and, by Cases 1 and 2, the formulas define the
asserted one-parameter family of real solutions.
\end{proof}

\begin{lemma}\label{lem:bricard}
Assume that $\alpha_i$, $\beta_i$, $\gamma_i$,
$\delta_i \in (0, \pi)$ satisfy condition~\eqref{eq:N1} 
and $\theta_i \in [-\pi, \pi) \setminus \{0\}$
for each
$i = 1, \ldots, 4$.
Then a polyhedron with flat angles $\alpha_i$, $\beta_i$, $\gamma_i$,
$\delta_i$ and central dihedral angles $\theta_i$ 
exists if and only if
$\delta_1 + \delta_2 + \delta_3 + \delta_4 = 2\pi$, $M_i > 0$ for
each $i = 1, \ldots, 4$, and the values
$Z = \cot{\frac{\theta_1}{2}}$, $W_2 = \cot{\frac{\theta_2}{2}}$,
$U = \cot{\frac{\theta_3}{2}}$, $W_1 = \cot{\frac{\theta_4}{2}}$
satisfy Bricard's system~\eqref{eq:N13}.
\end{lemma}

\begin{proof}
Recall that our dihedral angles are $\pi$ minus those
of~\cite{izmestiev-2017}, see~\cite{ErofeevIvanov2020}; hence, for a
polyhedron with central dihedral angles $\theta_i \neq 0$, the
variables of Section~2.2 of~\cite{izmestiev-2017} take the values
$\tan{\frac{\pi - \theta_i}{2}} = \cot{\frac{\theta_i}{2}}$.

($\Rightarrow$) We have
$\delta_1 + \delta_2 + \delta_3 + \delta_4 = 2\pi$ and $M_i > 0$
for each $i = 1, \ldots, 4$ by Lemma~5
of~\cite{NURMATOV2026104102},
and the values
$(Z, W_2, U, W_1) = \bigl(\cot{\frac{\theta_1}{2}}, \, \cot{\frac{\theta_2}{2}}, \, \cot{\frac{\theta_3}{2}}, \, \cot{\frac{\theta_4}{2}}\bigr)$
satisfy system~\eqref{eq:N13} by Lemmas~2.1 and~4.3 and Section~2.2 of~\cite{izmestiev-2017}.

($\Leftarrow$) Assume that
$\delta_1 + \delta_2 + \delta_3 + \delta_4 = 2\pi$, 
$M_i > 0$ for
each $i = 1, \ldots, 4$, and the values
$(Z, W_2, U, W_1) = \bigl(\cot{\frac{\theta_1}{2}}, \, \cot{\frac{\theta_2}{2}}, \, \cot{\frac{\theta_3}{2}}, \, \cot{\frac{\theta_4}{2}}\bigr)$
satisfy system~\eqref{eq:N13}.
By Lemma~4 of~\cite{NURMATOV2026104102}, $\overline{\alpha}_i$,
$\overline{\beta}_i$, $\overline{\gamma}_i$,
$\overline{\delta}_i \in (0, \pi)$, so, by Lemma~4.1
of~\cite{izmestiev-2017}, a spherical quadrilateral with the side
lengths $\alpha_i$, $\beta_i$, $\gamma_i$, $\delta_i$ exists for each
$i = 1, \ldots, 4$. Hence, by Lemmas~2.1 and~4.3 and Section~2.2
of~\cite{izmestiev-2017}, the given solution
of~\eqref{eq:N13}
corresponds to a polyhedron with the flat angles $\alpha_i$,
$\beta_i$, $\gamma_i$, $\delta_i$ whose central dihedral angles
$\theta'_1, \ldots, \theta'_4$ satisfy
$\cot{\frac{\theta'_i}{2}} = \cot{\frac{\theta_i}{2}}$ for each
$i = 1, \ldots, 4$. Since $\theta \mapsto \cot{\frac{\theta}{2}}$
is a bijection of $[-\pi, \pi) \setminus \{0\}$ onto
$\mathbb{R}$, it follows that $\theta'_i = \theta_i$.
The faces of this polyhedron can be chosen convex:
moving $B_i$ and $C_i$ along their rays
from $A_i$, see Figure~\ref{fig:03}, changes neither the flat nor the
dihedral angles, and, with short enough edges $A_iB_i$ and $A_iC_i$, all
faces are convex, because the flat angles lie in $(0, \pi)$ and
$\delta_1 + \delta_2 + \delta_3 + \delta_4 = 2\pi$.
\end{proof}

\begin{proof}[Proof of Theorem~\ref{th:main}]
First, by~\eqref{eq:N4}, $M > 0$, and Lemmas~2 and~4
of~\cite{NURMATOV2026104102}, $r_i$, $s_i$, $f_i$,
$M \neq 0, 1$ are well-defined reals, $f_i > 0$, and
$\overline{\alpha}_i$, $\overline{\beta}_i$, $\overline{\gamma}_i$,
$\overline{\delta}_i \in (0, \pi)$; also
$\varepsilon_i \in \{\pm 1\}$ by~\eqref{eq:N1} and~\eqref{eq:N2}.
Hence $x_i$, $y_i$, $z_i \neq -1, 0$ are well-defined reals with
$x_i(1 + x_i)$, $z_i(1 + z_i) > 0$, $1 + u x_i$, $1 + u z_i > 0$, and
$u x_i y_i z_i > 0$ by~\eqref{eq:N3}, \eqref{eq:N4}, \eqref{eq:N14},
and Lemma~4 of~\cite{NURMATOV2026104102}; also $k$, $k' \in (0, 1)$,
$u \neq 0$, and $u < 1$, see~\eqref{eq:N7} and~\eqref{eq:N14}. Next, $x_2 = x_1$
and $y_3 = y_2$ by~\eqref{eq:N3} and~\eqref{eq:N5}; since
$u x_i y_i z_i > 0$ for $i = 2, 3$, it follows that
$\sign{(x_1 x_3)} = \sign{(z_2 z_3)} = \xi_{23}$, see~\eqref{eq:P5}, and
$\sign{(y_2)} = \sign{(u x_1 z_2)}$; in particular,
$\frac{z_2 z_3}{x_1 x_3} > 0$. The radicands in the formulas for
$H_1$, $H_2$, and $H_3$ are strictly positive because
$z_i(1 + z_i) > 0$ and $1 + u z_i > 0$ for $i = 2, 3$. Moreover,
by~\eqref{eq:N3}, \eqref{eq:N7}, \eqref{eq:N11}, and~\eqref{eq:N14},
we have $g_i = \operatorname{dn}{t_i} > 0$, and $k^2 = u$, $k'^2 = M$,
$g_i^2 = f_i = \frac{1 + z_i}{z_i}$ if $u > 0$, while
$k^2 = -\frac{u}{M}$, $k'^2 = \frac{1}{M}$,
$g_i^2 = \frac{1}{f_i} = \frac{z_i}{1 + z_i}$ if $u < 0$.

Next, by Lemma~\ref{lem:main:7}, the radicands that the formulas of the
theorem share with those of the lemma are strictly positive reals. 
The right-hand sides of the formulas of the theorem for
$\cot{\frac{\theta_1(t)}{2}}$, $\cot{\frac{\theta_2(t)}{2}}$, and
$\cot{\frac{\theta_4(t)}{2}}$ coincide with those of the formulas of
Lemma~\ref{lem:main:7} for $Z$, $W_2$, and $W_1$, respectively.
By~\eqref{eq:N3} and~\eqref{eq:N14}, the equalities $u z_2 z_3 = 1$
and $1 + u(z_2 + z_3 + z_2 z_3) = 0$ read $f_2 + f_3 - f_2 f_3 = M$ and
$f_2 f_3 = M(f_2 + f_3 - 1)$, respectively, so they cannot hold
simultaneously by Remark~\ref{rem:lem:gqs}. Using this together with
Remark~\ref{rem:lem:main:6:2}, it is verified in Section~10 of the
supplement (\texttt{helper.nb})~\cite{Nurmatov2026Repo} that, in
each of the cases (a), (b), and (c) of the theorem, the right-hand side
of the formula for $\cot{\frac{\theta_3(t)}{2}}$ equals the right-hand
side of the applicable $U$-formula of Lemma~\ref{lem:main:7}, and that
the denominators of the two formulas vanish for the same values of $t$.
Since, in addition, the denominators $u z_2 z_3 - 1$,
$2\sqrt{z_2 z_3 (1 + z_2)(1 + z_3)}$, and $x_1 x_3$ in the formulas of
the theorem do not depend on $t$ and are nonzero, it follows that the
theorem and Lemma~\ref{lem:main:7} exclude the same values of $t$.
Hence, by Lemma~\ref{lem:main:7}, for every $t \in I(x_1)$ outside the
excluded set, the right-hand sides of the four formulas of the theorem
are well-defined reals forming a real solution of Bricard's
system~\eqref{eq:N13}, and the set of such $t$ has nonempty interior.

Now fix any such $t$, and define
$\theta_1(t), \ldots, \theta_4(t) \in [-\pi, \pi) \setminus \{0\}$ by
the four formulas of the theorem; this is possible, and uniquely so,
because $\theta \mapsto \cot{\frac{\theta}{2}}$ is a bijection of
$[-\pi, \pi) \setminus \{0\}$ onto $\mathbb{R}$. 
The values
$Z = \cot{\frac{\theta_1(t)}{2}}$, $W_2 = \cot{\frac{\theta_2(t)}{2}}$,
$U = \cot{\frac{\theta_3(t)}{2}}$, $W_1 = \cot{\frac{\theta_4(t)}{2}}$
satisfy Bricard's system~\eqref{eq:N13} by Lemma~\ref{lem:main:7}.
Since, in addition,
$\delta_1 + \delta_2 + \delta_3 + \delta_4 = 2\pi$ and $M_i > 0$ by
hypothesis and~\eqref{eq:N4}, it follows, by Lemma~\ref{lem:bricard},
that a polyhedron with the flat angles $\alpha_i$, $\beta_i$,
$\gamma_i$, $\delta_i$ and the central dihedral angles
$\theta_1(t), \ldots, \theta_4(t)$ exists.
Such a polyhedron is determined by its central face $A_1A_2A_3A_4$, which,
by the proof of Lemma~2.1 of~\cite{izmestiev-2017}, can be any
quadrilateral with the angles $\delta_i$, and the lengths $|A_iB_i|$,
$|A_iC_i|$, which are free by the proof of Lemma~\ref{lem:bricard}, see
Figure~\ref{fig:03}. We choose them once for all~$t$, with the edges
$A_iB_i$ and $A_iC_i$ short enough for the faces to be convex; then the
faces do not depend on~$t$.
The set of such $t$ has nonempty
interior, so it contains a nonempty open interval $J$ with
$J \cap (-J) = \emptyset$. Since a congruence either preserves all the
oriented dihedral angles or reverses their signs, and
$\cot{\frac{\theta_1(t)}{2}} = t$, it follows that these polyhedra are
pairwise non-congruent for $t \in J$. Moreover, they depend
continuously on $t \in J$, because the values
$\cot{\frac{\theta_i(t)}{2}}$, $i = 1, \ldots, 4$, do: the radicands in
the formulas of the theorem are strictly positive reals and their
denominators do not vanish on $J$. 
Hence, for each $t^* \in J$ and $\epsilon > 0$ with
$[t^*, t^* + \epsilon) \subset J$, these polyhedra with
$t \in [t^*, t^* + \epsilon)$, reparametrized by $t - t^*$, form a flexion
of the polyhedron with $t = t^*$, see Section~\ref{sec:notation}.
Thus each of these polyhedra is flexible,
and the formulas of the theorem define its flexion.
Each of the two signs $e_0 \in \{\pm 1\}$ of
Lemma~\ref{lem:main:7} gives a flexion.

Finally, fix $e_0 \in \{\pm 1\}$ and such a $t$. The radicand in the
formula for $\theta_{41'}$ and $\theta_{22'}$ is a strictly positive
real; hence this formula excludes no further $t$. Indeed,
$x_1D(x_1, t) > 0$ because $t \in I(x_1)$, see~\eqref{eq:N16}; thus
$(x_1t^2 - 1)(1 - ux_1t^2)$ has the sign of $x_1$, which is the sign of
$1 + x_1$ because $x_1(1 + x_1) > 0$, and $1 + ux_1 > 0$. Moreover,
$\theta_{41'}$ and $\theta_{22'}$ are neither $0$ nor $\pm\pi$, because
$t \in I(x_1)$ and $x_2 = x_1$, as verified in Section~11 of the
supplement (\texttt{helper.nb})~\cite{Nurmatov2026Repo}; 
with the formulas for $\cot\frac{\theta_2(t)}{2}$ and
$\cot\frac{\theta_4(t)}{2}$, the formula for $\theta_{41'}$ and
$\theta_{22'}$ is then verified in the same section.
\end{proof}

\begin{remark}\label{rem:th:main:1}
Under the hypotheses of Theorem~\ref{th:main}, for $t$ as in
Theorem~\ref{th:main}, excluding the finitely many $t$ at which the
denominators below vanish, the corner dihedral angles $\theta_{11'}(t)$
and $\theta_{12'}(t)$ of the flexions of Theorem~\ref{th:main} are
nonzero and satisfy
\begin{equation*}
\begin{aligned}
    &\cot{\frac{\theta_{11'}(t)}{2}}
    = \frac{\varepsilon_1 \Bigl(
    t\sqrt{u x_1 \tilde{y}_1(1 + \tilde{z}_1)(1 + u \tilde{z}_1)}
    - e_0\,e_1 \sqrt{u \tilde{y}_1 \tilde{z}_1 D(x_1, t)}
    \Bigr)}
    {\tilde{y}_1 u (\tilde{z}_1 + x_1 t^2)},
    \\[1mm]
    &\cot{\frac{\theta_{12'}(t)}{2}}
    = \frac{\varepsilon_2 \Bigl(
    t\sqrt{u x_1 \tilde{y}_2(1 + \tilde{z}_2)(1 + u \tilde{z}_2)}
    + e_0\,e_2 \sqrt{u \tilde{y}_2 \tilde{z}_2 D(x_1, t)}
    \Bigr)}
    {\tilde{y}_2 u (\tilde{z}_2 + x_1 t^2)},
\end{aligned}
\end{equation*}
where $\tilde{y}_i \coloneqq -\frac{1 + y_i}{1 + u y_i}$ and 
$\tilde{z}_i \coloneqq -\frac{1 + z_i}{1 + u z_i}$, $i = 1, 2$.
Indeed, under the substitution
$(\alpha_i, \beta_i, \gamma_i, \delta_i)
\mapsto (\delta_i, \gamma_i, \beta_i, \alpha_i)$
in~\eqref{eq:N2}, \eqref{eq:N3}, and~\eqref{eq:N14},
we get $(u, x_i, y_i, z_i, \varepsilon_i)
\mapsto (u, x_i, \tilde{y}_i, \tilde{z}_i, \varepsilon_i)$,
and the substituted angles satisfy~\eqref{eq:N1} as well, because
condition~\eqref{eq:N1} is invariant under permutations of $\alpha_i$,
$\beta_i$, $\gamma_i$, $\delta_i$.
Thus, by~\eqref{eq:N1}--\eqref{eq:N4}, \eqref{eq:N14}, and Lemmas~2
and~4 of~\cite{NURMATOV2026104102}, applied to the angles and to the
substituted angles, $\overline{\alpha}_i$, $\overline{\beta}_i$,
$\overline{\gamma}_i$, $\overline{\delta}_i \in (0, \pi)$,
$\varepsilon_i \in \{\pm 1\}$, and $x_i$, $y_i$, $z_i$, $\tilde{y}_i$,
$\tilde{z}_i \neq -1, 0$ are well-defined reals with
$u x_1 y_i z_i > 0$, $\tilde{z}_i(1 + \tilde{z}_i) > 0$,
$1 + u\tilde{z}_i > 0$, and $u x_1 \tilde{y}_i \tilde{z}_i > 0$,
because $x_2 = x_1$. 
Hence all the radicands in the two formulas above are strictly positive
for $t \in I(x_1)$,
see~\eqref{eq:N16}. The substitution, the inequalities
$\theta_{11'}(t)$, $\theta_{12'}(t) \neq 0$, and the formulas are verified in
Section~12 of the supplement (\texttt{helper.nb})~\cite{Nurmatov2026Repo}.
\end{remark}

\begin{remark}\label{rem:th:main:2}
Under the hypotheses of Theorem~\ref{th:main}, for $j = 1, 2, 3$ and all but
finitely many $t$ as in Theorem~\ref{th:main}, the polyhedra of the flexions of
Theorem~\ref{th:main} with the parameter $t$ are polyhedra of the flexions of
Theorem~\ref{th:main} applied after the cyclic permutation of the indices
by~$j$, with $e_{1+j}, \ldots, e_{4+j}$ in place of $e_1, \ldots, e_4$,
where, for each $i=1, \ldots, 4$, $i+j$ is taken modulo~$4$, with values in $\{1, \ldots, 4\}$.
Indeed, $x_i$, $y_i$, $z_i \neq -1, 0$ are well-defined reals with
$z_i(1 + z_i) > 0$ and $1 + u z_i > 0$ by~\eqref{eq:N3}, \eqref{eq:N4},
\eqref{eq:N14}, and Lemmas~2 and~4 of~\cite{NURMATOV2026104102};
$F(v + \lambda) = F(v)$ for $v \in \mathbb{C}$ and $\lambda \in \Lambda$
by~\eqref{eq:P1}; the points of the family of $e_1, \ldots, e_4$ are
solutions of~\eqref{eq:N13} by Remark~\ref{rem:lem:main:2}; and, for $t$ as
in Theorem~\ref{th:main}, the points of this family with $Z = t$ are the
points $\bigl(\cot{\frac{\theta_1}{2}}, \ldots, \cot{\frac{\theta_4}{2}}\bigr)$
of the polyhedra of the flexions of Theorem~\ref{th:main} with the
parameter $t$, for $e_0 = \pm 1$, see the proofs of
Lemma~\ref{lem:main:7} and Theorem~\ref{th:main}. 
Using this, the claim is verified in Section~13 of
the supplement (\texttt{helper.nb})~\cite{Nurmatov2026Repo}.
\end{remark}

The flexion formulas of Theorem~\ref{th:main} involve the signs
$e_1, \ldots, e_4$ of condition~\eqref{eq:N6}, whereas
Proposition~\ref{prop:shifts} only decides whether such signs exist.
Proposition~\ref{prop:signs} below determines the sign patterns satisfying
condition~\eqref{eq:N6}, so that the flexion of Theorem~\ref{th:main} can be
written down explicitly; its proof relies on the following lemma.

\begin{lemma}\label{lem:signs}
Assume that $M > 0$ and $M \neq 1$, and let $v$, $w \in \mathbb{C}\setminus(\mathbf{i}K' + \Lambda')$. Then
$v + w \in \Lambda$ if and only if, for some $n \in \{0, 1\}$,
\begin{equation}\label{eq:P8}
 (\operatorname{sn}, \operatorname{cn}, \operatorname{dn})(v)
    = \bigl((-1)^{1 + \varrho}\operatorname{sn}{w},
    \,(-1)^{n + \varrho}\operatorname{cn}{w},
    \,(-1)^{n}\operatorname{dn}{w}\bigr),
\end{equation}
where $\varrho \coloneqq 0$ if $M < 1$ and $\varrho \coloneqq n$ if $M > 1$.
\end{lemma}

\begin{proof}
First, $k \in (0, 1)$, because $M > 0$ and $M \neq 1$, see~\eqref{eq:N7}.
We next set $L \coloneqq 4K\mathbb{Z} + 4\mathbf{i}K'\mathbb{Z} = 2\Lambda'$ and, for
$n \in \{0, 1\}$, $\lambda_n \coloneqq 2n\mathbf{i}K'$ if $M < 1$ and
$\lambda_n \coloneqq 2nK + 2n\mathbf{i}K'$ if $M > 1$. 
Then $\lambda_n \in \Lambda'$ and, by~\eqref{eq:N9}, $\Lambda = L \cup (\lambda_1 + L)$;
for $M < 1$ this is immediate, and for $M > 1$
it follows from $4\mathbf{i}K' = 2(2K + 2\mathbf{i}K') - 4K$.
By~\eqref{eq:P1}, $\operatorname{sn}$, $\operatorname{cn}$,
$\operatorname{dn}$ are $L$-periodic, and 
$(\operatorname{sn}, \operatorname{cn}, \operatorname{dn})(\lambda_n - w)
    = \bigl((-1)^{1 + \varrho}\operatorname{sn}{w},
    \,(-1)^{n + \varrho}\operatorname{cn}{w},
    \,(-1)^{n}\operatorname{dn}{w}\bigr)$, that is, \eqref{eq:P8} reads
$(\operatorname{sn}, \operatorname{cn}, \operatorname{dn})(v)
= (\operatorname{sn}, \operatorname{cn}, \operatorname{dn})(\lambda_n - w)$.

($\Rightarrow$) If $v + w \in \Lambda$, then $v \in \lambda_n - w + L$ for some
$n \in \{0, 1\}$, so, by $L$-periodicity, we get 
$(\operatorname{sn}, \operatorname{cn}, \operatorname{dn})(v)
= (\operatorname{sn}, \operatorname{cn}, \operatorname{dn})(\lambda_n - w)$, 
which is~\eqref{eq:P8} for that $n$.

($\Leftarrow$) Assume~\eqref{eq:P8} for some $n \in \{0, 1\}$, and set
$w' \coloneqq \lambda_n - w$. 
Then $-w' \notin \mathbf{i}K' + \Lambda'$, because
$w \notin \mathbf{i}K' + \Lambda'$, $\lambda_n \in \Lambda'$, and the set
$\mathbf{i}K' + \Lambda'$ is invariant under translations by $\Lambda'$; 
also $0 \notin \mathbf{i}K' + \Lambda'$, see~\eqref{eq:P4}.
Moreover, $(\operatorname{sn}, \operatorname{cn}, \operatorname{dn})(v)
= (\operatorname{sn}, \operatorname{cn}, \operatorname{dn})(w')$
by~\eqref{eq:P8}. Hence, by the Pythagorean relations
and~\eqref{eq:N7},
$k'^2 + k^2\operatorname{cn}{v}\operatorname{cn}{w'} - \operatorname{dn}{v}\operatorname{dn}{w'}
= k'^2 + k^2\operatorname{cn}^2{w'} - \operatorname{dn}^2{w'} = 0$,
which is $G = 0$ for
$(v_1, v_2, v_3, v_4) = (v, -w', 0, 0)$ in Lemma~\ref{lem:shifts:1}, because
$\operatorname{cn}$, $\operatorname{dn}$ are even, see~\eqref{eq:P1}, and
$(\operatorname{sn}, \operatorname{cn}, \operatorname{dn})(0) = (0, 1, 1)$, 
see~\cite[\S 22.11]{WhittakerWatson2021}. Thus, by Lemma~\ref{lem:shifts:1},
$v - w' \in L$ or $v + w' \in L$. 
In the latter case, by $L$-periodicity and~\eqref{eq:P1},
$\operatorname{sn}{v} = \operatorname{sn}(-w') = -\operatorname{sn}{w'}
= -\operatorname{sn}{v}$, so $v$, $w' \in \Lambda'$ by~\eqref{eq:P4}. 
Hence $2w' \in 2\Lambda' = L$ and
$v - w' = (v + w') - 2w' \in L$. 
Thus, $v - w' \in L$ in both cases, and
$v + w = (v - w') + \lambda_n \in \lambda_n + L \subseteq \Lambda$.
\end{proof}

\begin{proposition}\label{prop:signs}
Assume that $\alpha_i$, $\beta_i$, $\gamma_i$, $\delta_i \in (0, \pi)$ satisfy
conditions~\eqref{eq:N1}, \eqref{eq:N4}, and~\eqref{eq:N5} for each
$i = 1, \ldots, 4$, and $M > 0$. Let $e_1$, $e_2$, $e_3$, $e_4 \in \{\pm 1\}$ and
$\{i, j, l, m\} = \{1, 2, 3, 4\}$.
Then $e_1 t_1 + e_2 t_2 + e_3 t_3 + e_4 t_4 \in \Lambda$ if and only if
one of the following holds:
\begin{itemize}
\item[(i)] neither $e_i t_i + e_j t_j$ nor $e_l t_l + e_m t_m$ is a pole, and,
for some $n \in \{0, 1\}$,
\begin{equation*}
    h_{ij1} = (-1)^{1 + \varrho} \chi \xi_{ij} h_{lm1},
    \quad
    h_{ij2} = (-1)^{n + \varrho} \chi \xi_{ij} h_{lm2},
    \quad
    h_{ij3} = (-1)^{n} h_{lm3},
\end{equation*}
where $\varrho = 0$ if $M < 1$ and $\varrho = n$ if $M > 1$, and $h_{ij1}$,
$h_{ij2}$, $h_{ij3}$ are given by~\eqref{eq:P6} if
$k^2 \neq (g_i^2 - 1)(g_j^2 - 1)$ and by~\eqref{eq:P7} otherwise, and analogously
for $h_{lm1}$, $h_{lm2}$, $h_{lm3}$;
\item[(ii)] both $e_i t_i + e_j t_j$ and $e_l t_l + e_m t_m$ are poles, and
$\chi = 1$ if $M < 1$, and $e_i e_l = -\chi$ if $M > 1$.
\end{itemize}
\end{proposition}

\begin{proof}
First, by~\eqref{eq:N4} and Lemmas~2 and~4
of~\cite{NURMATOV2026104102}, $r_\nu$, $s_\nu$, $f_\nu$,
$M \neq 0, 1$ are well-defined reals and $f_\nu > 0$, for each
$\nu = 1, \ldots, 4$. Hence
$p_\nu$, $q_\nu \in \mathbb{R}_{>0} \cup \mathbf{i}\mathbb{R}_{>0}$,
see~\eqref{eq:N10}, $k \in (0, 1)$, because $M > 0$ and
$M \neq 1$, see~\eqref{eq:N7}, and the phase shift $t_\nu$ is
well-defined, with $t_\nu \in (\mu_\nu K, \mu_\nu K + \mathbf{i}K')$ for
some $\mu_\nu \in \{0, 1, 2, 3\}$, see~\eqref{eq:N11} and
Table~\ref{tab:tpqsigma}. 
Also, by~\eqref{eq:N3} and~\eqref{eq:N14},
$u$, $x_\nu$, $y_\nu$, $z_\nu \neq 0$ are well-defined reals, 
and $\chi \in \{\pm 1\}$, because
$\varepsilon_1, \ldots, \varepsilon_4 \in \{\pm 1\}$
by~\eqref{eq:N1} and~\eqref{eq:N2}, and $z_1$, $z_3 \neq 0$, see~\eqref{eq:N12}. 

Second, observe that $\prod_{\nu = 1}^{4} p_\nu q_\nu = p_1^2\, p_3^2\, q_1^2\, q_2^2
= 1/(x_1 x_3 y_1 y_2)$ is a nonzero real, by~\eqref{eq:N3}, \eqref{eq:N5},
and~\eqref{eq:N10}. By~\eqref{eq:N3}, \eqref{eq:N4}, \eqref{eq:N14}, and Lemma~4
of~\cite{NURMATOV2026104102}, we have $u x_\nu y_\nu z_\nu > 0$ for
$\nu \in \{1, 2, 3, 4\}$.
Hence
$\sign{(z_1 z_3)}
= \sign{(u x_1 y_1 u x_3 y_3)}
= \sign{(x_1 x_3 y_1 y_2)}$, because $y_3 = y_2$ by~\eqref{eq:N3}
and~\eqref{eq:N5}. 
Thus,
$\sign{\bigl(\prod_{\nu = 1}^{4} p_\nu q_\nu\bigr)}
= \sign{(x_1 x_3 y_1 y_2)}
= \sign{(z_1 z_3)}$.
Moreover, by Lemma~\ref{lem:main:6}, the radicands
in~\eqref{eq:P5} are strictly positive, so
$\xi_{ij} = \sign{\bigl(p_i^2 q_i^2 p_j^2 q_j^2\bigr)}$,
$\xi_{lm} = \sign{\bigl(p_l^2 q_l^2 p_m^2 q_m^2\bigr)} \in \{\pm 1\}$,
and
$\varepsilon_1 \varepsilon_2 \varepsilon_3 \varepsilon_4\, S_{ij} S_{lm}
    = 
    \frac{\sqrt{\prod_{\nu = 1}^{4} p_\nu^2 q_\nu^2}}{\prod_{\nu = 1}^{4} p_\nu q_\nu}
    =
    \sign{\bigl(\prod_{\nu = 1}^{4} p_\nu q_\nu\bigr)}
    =  \sign{(z_1 z_3)}$,
that is, $S_{ij} S_{lm} = \chi$, by~\eqref{eq:N12}. In particular, $S_{ij}$, $S_{lm}\neq 0$ and
$S_{lm} = \chi\, \xi_{ij}\, S_{ij}$, because $S_{ij}^2 = \xi_{ij}$ by~\eqref{eq:P5}.

Third, for the rest of the proof, 
we set $v \coloneqq e_i t_i + e_j t_j$ and $w \coloneqq e_l t_l + e_m t_m$. 
If exactly one of $v$, $w$ is a pole, then $v + w \notin \Lambda$,
because $\Lambda \subseteq \Lambda'$, see~\eqref{eq:N9}, and $\mathbf{i}K' + \Lambda'$ is invariant
under translations by $\Lambda'$ and under reflection, that is,
$-(\mathbf{i}K' + \Lambda') = -\mathbf{i}K' + \Lambda' = \mathbf{i}K' + \Lambda'$,
because $\Lambda' = -\Lambda'$ and $2\mathbf{i}K' \in \Lambda'$;
in this case neither (i) nor (ii) holds.

Fourth, assume that both $v$ and $w$ are poles. Then, by
Lemma~\ref{lem:main:3}, $e_i = e_j$, $e_l = e_m$, and $t_i + t_j$,
$t_l + t_m$ are poles, so, by Remark~\ref{rem:lem:main:6:1},
$v = e_i(2m'K + \mathbf{i}K')$ and $w = e_l(2m''K + \mathbf{i}K')$ for some
$m'$, $m'' \in \mathbb{Z}$ with $(-1)^{m'} = S_{ij}$ and $(-1)^{m''} = S_{lm}$.
Hence $(-1)^{m' + m''} = \chi$ because $S_{ij} S_{lm} = \chi$. 
Now
$v + w = \pm\bigl(2(m' + m'')K + 2\mathbf{i}K'\bigr)$ if $e_i = e_l$, and
$v + w = \pm 2(m' - m'')K$ if $e_i = -e_l$. By~\eqref{eq:N9}, for $M < 1$ each
of these points lies in $\Lambda$ if and only if $m' + m''$ is even
(because $m' - m'' \equiv m' + m'' \pmod{2}$), that is, $\chi = 1$;
for $M > 1$, since
$2(m' + m'')K + 2\mathbf{i}K' = 2(m' + m'' - 1)K + (2K + 2\mathbf{i}K')$, 
it follows that the first
point lies in $\Lambda$ if and only if $m' + m''$ is odd, that is, $\chi = -1$,
and the second if and only if $m' + m''$ is even, that is, $\chi = 1$. Thus
$v + w \in \Lambda$ if and only if $\chi = 1$ for $M < 1$, and if and only if
$e_i e_l = -\chi$ for $M > 1$, which is condition (ii); condition (i) does not
hold in this case.

Finally, assume that neither $v$ nor $w$ is a pole. Then, by
Lemma~\ref{lem:main:6},
$(\operatorname{sn}, \operatorname{cn}, \operatorname{dn})(v)
= (\mathbf{i}\, S_{ij} h_{ij1}, \, S_{ij} h_{ij2}, \, h_{ij3})$ and
$(\operatorname{sn}, \operatorname{cn}, \operatorname{dn})(w)
= (\mathbf{i}\, S_{lm} h_{lm1}, \, S_{lm} h_{lm2}, \, h_{lm3})$,
and, by Lemma~\ref{lem:signs}, $v + w \in \Lambda$ if and only if, for some
$n \in \{0, 1\}$,
$S_{ij} h_{ij1} = (-1)^{1 + \varrho} S_{lm} h_{lm1}$,
$S_{ij} h_{ij2} = (-1)^{n + \varrho} S_{lm} h_{lm2}$, and
$h_{ij3} = (-1)^{n} h_{lm3}$. Substituting $S_{lm} = \chi\, \xi_{ij}\, S_{ij}$
and dividing the first two equalities by $S_{ij} \neq 0$, we get the three
equalities of condition (i); condition (ii) does not hold in this case.
Since exactly one of the three cases above holds, it follows that $v + w \in \Lambda$ if and
only if (i) or (ii) holds.
\end{proof}

\begin{remark}\label{rem:prop:signs}
Proposition~\ref{prop:signs} gives the following description of the sign
patterns $(e_i, e_j, e_l, e_m)$ satisfying
$e_i t_i + e_j t_j + e_l t_l + e_m t_m \in \Lambda$, called \emph{witnesses}
below, for fixed indices
$i$, $j$, $l$, $m$ with $\{i, j, l, m\} = \{1, 2, 3, 4\}$. Note first that
the patterns $(e_i, e_j, e_l, e_m)$ and $(-e_i, -e_j, -e_l, -e_m)$ are
witnesses simultaneously, because $\Lambda = -\Lambda$; thus, of the sixteen
patterns, only the eight with $e_i = 1$ need to be considered. By
Lemma~\ref{lem:main:3} and Remark~\ref{rem:lem:main:6:1}, $e_i t_i + e_j t_j$ is
a pole if and only if $e_i = e_j$ and $k^2 = (g_i^2 - 1)(g_j^2 - 1)$, and
analogously for $e_l t_l + e_m t_m$. Three configurations arise depending on whether
the equalities $k^2 = (g_i^2 - 1)(g_j^2 - 1)$ and $k^2 = (g_l^2 - 1)(g_m^2 - 1)$ hold.
If neither holds, then each of
the eight patterns $(1, e_j, e_l, e_m)$ is a witness if and only if it satisfies condition (i).
If exactly one of these equalities holds, say the first, then the four
patterns $(1, 1, e_l, e_m)$ make $e_i t_i + e_j t_j$ a pole and are not
witnesses, and each of the remaining four, $(1, -1, e_l, e_m)$, is a witness if
and only if it satisfies condition (i); analogously if the second equality holds.
If both equalities hold, then, for each of the four patterns
$(1, 1, e_l, -e_l)$ and $(1, -1, e_l, e_l)$ with $e_l \in \{\pm 1\}$, exactly
one of the points $e_i t_i + e_j t_j$ and $e_l t_l + e_m t_m$ is a pole, so
these patterns are not witnesses. For the two patterns $(1, 1, e_l, e_l)$,
condition (ii) applies: for $M < 1$, both are witnesses if $\chi = 1$ and
neither is if $\chi = -1$; for $M > 1$, exactly one is, namely the one with
$e_l = -\chi$. Each of the remaining two patterns $(1, -1, e_l, -e_l)$ is 
a witness if and only if it satisfies condition (i); 
for $M > 1$ they need to be tested only if all
witnesses are sought, because a witness has already been found.
\end{remark}

\section{Flexible \texorpdfstring{$m \times n$}{m x n} nets}
\label{sec:examples}
In this section, we assemble flexible $m \times n$ nets from $3 \times 3$
GQS-nets. Section~\ref{sec:ex:construction} constructs such nets and their
flexion in closed form (Proposition~\ref{prop:ex:flex}), and
Section~\ref{sec:ex:algorithm} turns the construction into an algorithm
and gives examples.

\subsection{Construction}
\label{sec:ex:construction}
The flexibility criterion
of~\cite[Theorem~5]{IzmestievNawratilSamarehRadSharifmoghaddam2026},
recalled in Section~\ref{sec:intro}, presupposes a net whose
$3 \times 3$ subnets are flexible; we construct such a net and its flexion 
without using this criterion: Theorem~\ref{th:main} gives flexions of the subnets
explicitly, and we make them agree on the common faces.
For the GQS-nets of Lemma~\ref{lem:ex:patterns} below, the three requirements
described at the beginning of Section~\ref{sec:notation} can be met explicitly:
their flat angles are determined by those at one of the vertices through
equalities and complements to~$\pi$ in one of four patterns; the faces can
be chosen to fit together; and the agreement of the flexions of
neighbouring subnets can be verified directly.

In Definition~\ref{def:gqs}, we take the pairs of indices $\{2, 3\}$ and
$\{1, 4\}$; the pairs $\{1, 2\}$ and $\{3, 4\}$ can be treated in the
same way, while the pairs $\{1, 3\}$ and $\{2, 4\}$ call for a different
approach, because the equations of Proposition~\ref{prop:eqMrs} link only
neighbouring vertices.
Assume that $\alpha_i$, $\beta_i$, $\gamma_i$, $\delta_i \in (0, \pi)$
satisfy condition~\eqref{eq:N1} and $\overline{\alpha}_i$,
$\overline{\beta}_i$, $\overline{\gamma}_i$,
$\overline{\delta}_i \in (0, \pi)$ for each $i = 1, \ldots, 4$, so that
all sines below are positive. By Proposition~\ref{prop:eqMrs},
conditions~\eqref{eq:N4} and~\eqref{eq:N5} are equivalent to the eight
equations of that proposition; and, under condition~\eqref{eq:N4}, 
in each of conditions~(i)--(iv) of Definition~\ref{def:gqs} with
$(i, j, l, m) = (2, 3, 1, 4)$, the equality for the pair $(2, 3)$ is
equivalent to two of the equations (a)--(d) of Lemma~\ref{lem:gqs}
applied to this pair, and the equality for the pair $(1, 4)$ to the same
two applied to the pair $(1, 4)$. In total, these are twelve
equations between sine ratios in the sixteen variables $\alpha_i$,
$\beta_i$, $\gamma_i$, $\delta_i$; the inequalities of conditions~(i)
and~(iii) are dealt with in Lemma~\ref{lem:ex:patterns} below. Consider
first the equalities $f_2 = f_3$ and $f_1 = f_4$ of condition~(i), which
are equivalent to equations (a) and (b) of Lemma~\ref{lem:gqs} applied to
these two pairs. Taking products and quotients of the equations for the
pairs $(2, 3)$ and $(1, 4)$, we get that the twelve equations are
equivalent to the system
\begin{equation}\label{eq:E1}
\begin{aligned}
    \frac{\sin{\gamma_2}}{\sin{\beta_2}} &= \frac{\sin{\gamma_3}}{\sin{\beta_3}},
    &\qquad
    \frac{\sin{\overline{\gamma}_2}}{\sin{\overline{\beta}_2}} &= \frac{\sin{\overline{\gamma}_3}}{\sin{\overline{\beta}_3}},
    \\
    \frac{\sin{\delta_2}}{\sin{\alpha_2}} &= \frac{\sin{\delta_3}}{\sin{\alpha_3}},
    &\qquad
    \frac{\sin{\overline{\delta}_2}}{\sin{\overline{\alpha}_2}} &= \frac{\sin{\overline{\delta}_3}}{\sin{\overline{\alpha}_3}},
    \\
    \frac{\sin{\gamma_1}}{\sin{\beta_1}} &= \frac{\sin{\gamma_4}}{\sin{\beta_4}},
    &\qquad
    \frac{\sin{\overline{\gamma}_1}}{\sin{\overline{\beta}_1}} &= \frac{\sin{\overline{\gamma}_4}}{\sin{\overline{\beta}_4}},
    \\
    \frac{\sin{\delta_1}}{\sin{\alpha_1}} &= \frac{\sin{\delta_4}}{\sin{\alpha_4}},
    &\qquad
    \frac{\sin{\overline{\delta}_1}}{\sin{\overline{\alpha}_1}} &= \frac{\sin{\overline{\delta}_4}}{\sin{\overline{\alpha}_4}},
    \\
    \frac{\sin{\alpha_1}\sin{\delta_1}}{\sin{\beta_1}\sin{\gamma_1}} &= \frac{\sin{\alpha_2}\sin{\delta_2}}{\sin{\beta_2}\sin{\gamma_2}},
    &\qquad
    \frac{\sin{\overline{\alpha}_1}\sin{\overline{\delta}_1}}{\sin{\overline{\beta}_1}\sin{\overline{\gamma}_1}} &= \frac{\sin{\overline{\alpha}_2}\sin{\overline{\delta}_2}}{\sin{\overline{\beta}_2}\sin{\overline{\gamma}_2}},
    \\
    \frac{\sin{\alpha_3}\sin{\delta_3}}{\sin{\beta_3}\sin{\gamma_3}} &= \frac{\sin{\alpha_4}\sin{\delta_4}}{\sin{\beta_4}\sin{\gamma_4}},
    &\qquad
    \frac{\sin{\overline{\alpha}_3}\sin{\overline{\delta}_3}}{\sin{\overline{\beta}_3}\sin{\overline{\gamma}_3}} &= \frac{\sin{\overline{\alpha}_4}\sin{\overline{\delta}_4}}{\sin{\overline{\beta}_4}\sin{\overline{\gamma}_4}}.
\end{aligned}
\end{equation}
Conditions~(ii), (iii), and~(iv) lead to the same computation with equations
(c) and (d), (a) and (d), and (b) and (c) of Lemma~\ref{lem:gqs} in place of
(a) and (b). Since (c) and (d) are (a) and (b) with inverted right-hand sides,
it follows that the resulting systems are~\eqref{eq:E1} with the following
substitutions in the right-hand sides of its first four rows, for
$i \in \{3, 4\}$: for condition~(ii), $\alpha_i \leftrightarrow \beta_i$,
$\gamma_i \leftrightarrow \delta_i$,
$\overline{\alpha}_i \leftrightarrow \overline{\beta}_i$, and
$\overline{\gamma}_i \leftrightarrow \overline{\delta}_i$; for
condition~(iii), only $\overline{\alpha}_i \leftrightarrow \overline{\beta}_i$
and $\overline{\gamma}_i \leftrightarrow \overline{\delta}_i$; and for
condition~(iv), only $\alpha_i \leftrightarrow \beta_i$ and
$\gamma_i \leftrightarrow \delta_i$.
The last two rows are the same in all
four cases. All the computations are verified in Section~1 of the
supplement (\texttt{ehelper.nb})~\cite{Nurmatov2026Repo}.
The following lemma states the simplest solutions we found
of~\eqref{eq:E1} and of its three analogues: four families of flat
angles satisfying these systems identically.

\begin{lemma}\label{lem:ex:patterns}
Assume that $\alpha_1$, $\beta_1$, $\gamma_1$, $\delta_1 \in (0, \pi)$ satisfy
condition~\eqref{eq:N1}, and $\alpha_i$, $\beta_i$, $\gamma_i$,
$\delta_i$, $i = 1, \ldots, 4$, are given by one of the following systems:
\begin{itemize}
\item[(a)] $\alpha_1 = \alpha_3 = \pi - \alpha_2 = \pi - \alpha_4$,
$\beta_1 = \beta_3 = \pi - \beta_2 = \pi - \beta_4$,
$\gamma_1 = \gamma_3 = \pi - \gamma_2 = \pi - \gamma_4$,
$\delta_1 = \delta_3 = \pi - \delta_2 = \pi - \delta_4$;
\item[(b)] $\alpha_1 = \beta_4 = \pi - \alpha_2 = \pi - \beta_3$,
$\beta_1 = \alpha_4 = \pi - \beta_2 = \pi - \alpha_3$,
$\gamma_1 = \delta_4 = \pi - \gamma_2 = \pi - \delta_3$,
$\delta_1 = \gamma_4 = \pi - \delta_2 = \pi - \gamma_3$;
\item[(c)] $\alpha_1 = \alpha_4 = \pi - \alpha_2 = \pi - \alpha_3$,
$\beta_1 = \beta_4 = \pi - \beta_2 = \pi - \beta_3$,
$\gamma_1 = \gamma_3 = \pi - \gamma_2 = \pi - \gamma_4$,
$\delta_1 = \delta_3 = \pi - \delta_2 = \pi - \delta_4$;
\item[(d)] $\alpha_1 = \beta_4 = \pi - \alpha_2 = \pi - \beta_3$,
$\beta_1 = \alpha_4 = \pi - \beta_2 = \pi - \alpha_3$,
$\gamma_1 = \delta_3 = \pi - \gamma_2 = \pi - \delta_4$,
$\delta_1 = \gamma_3 = \pi - \delta_2 = \pi - \gamma_4$.
\end{itemize}
If $\overline{\alpha}_1$, $\overline{\beta}_1$, $\overline{\gamma}_1$,
$\overline{\delta}_1 \in (0, \pi)$, then $\alpha_i$, $\beta_i$, $\gamma_i$,
$\delta_i \in (0, \pi)$ satisfy conditions~\eqref{eq:N1}, \eqref{eq:N4},
and~\eqref{eq:N5} for each $i = 1, \ldots, 4$,
$\delta_1 + \delta_2 + \delta_3 + \delta_4 = 2\pi$, $M > 0$, and
condition (i), (ii), (iii), (iv) of Definition~\ref{def:gqs} in cases
(a), (b), (c), (d), respectively; in particular, every polyhedron with
these flat angles is a GQS-net.
\end{lemma}

\begin{proof}
First, in each of the cases (a)--(d), for each $i = 1, \ldots, 4$,
$\alpha_i$, $\beta_i$, $\gamma_i$, $\delta_i \in (0, \pi)$ satisfy
condition~\eqref{eq:N1} and $\overline{\alpha}_i$, $\overline{\beta}_i$,
$\overline{\gamma}_i$, $\overline{\delta}_i \in (0, \pi)$, hence
$M_i > 0$ by~\eqref{eq:N3}. Also,
$\delta_1 + \delta_2 + \delta_3 + \delta_4 = 2\pi$. Next, the flat angles
satisfy conditions~\eqref{eq:N4} and~\eqref{eq:N5}, and condition (i),
(ii), (iii), (iv) of Definition~\ref{def:gqs} with
$(i, j, l, m) = (2, 3, 1, 4)$ in cases (a), (b), (c), (d), respectively.
Hence $M > 0$ by~\eqref{eq:N4}, and every polyhedron with these flat
angles is a GQS-net by Definition~\ref{def:gqs}. All
the computations are verified in Section~2 of the supplement
(\texttt{ehelper.nb})~\cite{Nurmatov2026Repo}.
\end{proof}

We now assemble $m \times n$ nets all of whose $3 \times 3$ subnets are
GQS-nets given by Lemma~\ref{lem:ex:patterns}, and make their flexions
agree. We first fix the notation.

\begin{figure}[htbp!]
    \centering
    \input{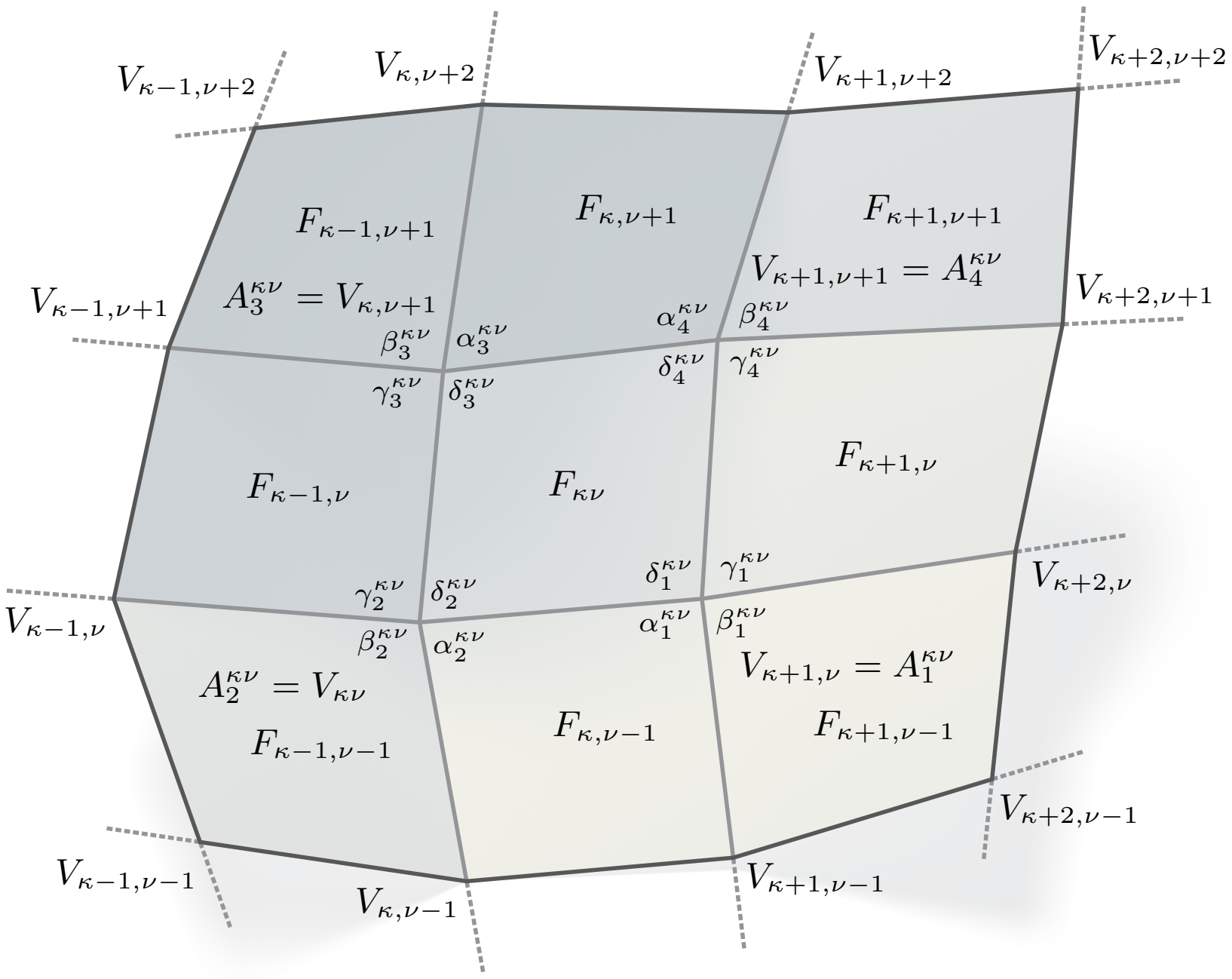}%
    \caption{Notation for vertices, faces, and flat angles of the subnet
with the central face $F_{\kappa\nu}$; dashed segments indicate the rest
of the net.}
    \label{fig:04}
\end{figure}

For $1 \leq \kappa \leq m - 2$ and $1 \leq \nu \leq n - 2$, the
\emph{block}~$(\kappa, \nu)$ is the set of the vertex labels $V_{ij}$ with
$\kappa - 1 \leq i \leq \kappa + 2$ and $\nu - 1 \leq j \leq \nu + 2$,
and a $3 \times 3$ net with these labels is called a \emph{subnet} with the
\emph{central face} $F_{\kappa\nu}$;
in an $m \times n$ net, the points
$V_{ij}$ with these labels form such a subnet. For each
$\nu = 1, \ldots, n - 2$, the blocks $(1, \nu), \ldots, (m - 2, \nu)$ form
the \emph{row}~$\nu$ of blocks.
We transfer the notation of the $3 \times 3$ net of
Figure~\ref{fig:03} to the block~$(\kappa, \nu)$ by shifting the indices:
the vertex label $V_{ij}$ of the block corresponds to the vertex
$V_{i-\kappa+1,j-\nu+1}$ of the $3 \times 3$ net, and we write
$A^{\kappa\nu}_i$ for the vertex labels of the block that correspond to
$A_i$ of the $3 \times 3$ net, and, for a subnet with the central face $F_{\kappa\nu}$,
$\alpha^{\kappa\nu}_i$, $\beta^{\kappa\nu}_i$, $\gamma^{\kappa\nu}_i$,
$\delta^{\kappa\nu}_i$ for its flat angles at its vertex $A^{\kappa\nu}_i$
that correspond to $\alpha_i$, $\beta_i$, $\gamma_i$, $\delta_i$,
respectively, $i = 1, \ldots, 4$. Here
$A^{\kappa\nu}_1 = V_{\kappa+1,\nu}$, $A^{\kappa\nu}_2 = V_{\kappa\nu}$,
$A^{\kappa\nu}_3 = V_{\kappa,\nu+1}$,
$A^{\kappa\nu}_4 = V_{\kappa+1,\nu+1}$; see Figure~\ref{fig:04}.

Two blocks are \emph{neighbours} if they are $(\kappa, \nu)$ and
$(\kappa, \nu + 1)$, or $(\kappa, \nu)$ and $(\kappa + 1, \nu)$. The
neighbours $(\kappa, \nu)$ and $(\kappa, \nu + 1)$ share the vertex labels
$A^{\kappa\nu}_4 = A^{\kappa,\nu+1}_1$ and
$A^{\kappa\nu}_3 = A^{\kappa,\nu+1}_2$, and the neighbours
$(\kappa, \nu)$ and $(\kappa + 1, \nu)$ share
$A^{\kappa\nu}_1 = A^{\kappa+1,\nu}_2$ and
$A^{\kappa\nu}_4 = A^{\kappa+1,\nu}_3$; in an $m \times n$ net, the flat
angles of its subnets with the central faces $F_{\kappa\nu}$ and
$F_{\kappa,\nu+1}$, or $F_{\kappa\nu}$ and $F_{\kappa+1,\nu}$, at these
vertices that lie in the same face coincide, that is,
\begin{equation}\label{eq:E2}
\begin{aligned}
(\alpha^{\kappa,\nu+1}_1, \beta^{\kappa,\nu+1}_1, \gamma^{\kappa,\nu+1}_1,
\delta^{\kappa,\nu+1}_1)
&= (\delta^{\kappa\nu}_4, \gamma^{\kappa\nu}_4, \beta^{\kappa\nu}_4,
\alpha^{\kappa\nu}_4),
\\
(\alpha^{\kappa,\nu+1}_2, \beta^{\kappa,\nu+1}_2, \gamma^{\kappa,\nu+1}_2,
\delta^{\kappa,\nu+1}_2)
&= (\delta^{\kappa\nu}_3, \gamma^{\kappa\nu}_3, \beta^{\kappa\nu}_3,
\alpha^{\kappa\nu}_3),
\\
(\alpha^{\kappa+1,\nu}_2, \beta^{\kappa+1,\nu}_2, \gamma^{\kappa+1,\nu}_2,
\delta^{\kappa+1,\nu}_2)
&= (\beta^{\kappa\nu}_1, \alpha^{\kappa\nu}_1, \delta^{\kappa\nu}_1,
\gamma^{\kappa\nu}_1),
\\
(\alpha^{\kappa+1,\nu}_3, \beta^{\kappa+1,\nu}_3, \gamma^{\kappa+1,\nu}_3,
\delta^{\kappa+1,\nu}_3)
&= (\beta^{\kappa\nu}_4, \alpha^{\kappa\nu}_4, \delta^{\kappa\nu}_4,
\gamma^{\kappa\nu}_4).
\end{aligned}
\end{equation}

Given numbers $\alpha_1$, $\beta_1$, $\gamma_1$, $\delta_1 \in (0, \pi)$,
indices $1 \leq \kappa_0 \leq m - 2$ and $1 \leq \nu_0 \leq n - 2$, and one
of the systems (a)--(d) of Lemma~\ref{lem:ex:patterns} for each
$\nu = 1, \ldots, n - 2$, denoted $\rho_\nu$, the
\emph{$(\rho_1, \ldots, \rho_{n-2})$-assembly} is the collection of
numbers $\alpha^{\kappa\nu}_i$, $\beta^{\kappa\nu}_i$,
$\gamma^{\kappa\nu}_i$, $\delta^{\kappa\nu}_i$, $i = 1, \ldots, 4$,
$1 \leq \kappa \leq m - 2$, $1 \leq \nu \leq n - 2$, such that
$(\alpha^{\kappa_0\nu_0}_1, \beta^{\kappa_0\nu_0}_1,
\gamma^{\kappa_0\nu_0}_1, \delta^{\kappa_0\nu_0}_1)
= (\alpha_1, \beta_1, \gamma_1, \delta_1)$, the numbers with each
superscript $\kappa\nu$ satisfy the system $\rho_\nu$, the first two
equalities of~\eqref{eq:E2} hold for $1 \leq \kappa \leq m - 2$ and
$1 \leq \nu \leq n - 3$, and the last two for $1 \leq \kappa \leq m - 3$
and $1 \leq \nu \leq n - 2$.
For all $\kappa$, $\nu$, and $i$, the numbers
$(\alpha^{\kappa\nu}_i, \beta^{\kappa\nu}_i, \gamma^{\kappa\nu}_i,
\delta^{\kappa\nu}_i)$ form a \emph{quadruple} of the assembly, and the
numbers $(\overline{\alpha}^{\kappa\nu}_i, \overline{\beta}^{\kappa\nu}_i,
\overline{\gamma}^{\kappa\nu}_i, \overline{\delta}^{\kappa\nu}_i)$
obtained from them by~\eqref{eq:N2} form its \emph{barred quadruple}.

\begin{remark}\label{rem:ex:assembly}
The assembly exists and is unique. Indeed, it is unique, because each of
the systems (a)--(d) determines the numbers with a superscript
$\kappa\nu$ from any one of their quadruples, \eqref{eq:E2} determines
two of these quadruples from the numbers with a neighbouring superscript
$\kappa \pm 1, \nu$ or $\kappa, \nu \pm 1$, and every pair
$(\kappa, \nu)$ is reached from $(\kappa_0, \nu_0)$ through neighbouring
pairs. It exists, because the numbers prescribed pair by pair, starting
from $(\kappa_0, \nu_0)$ and each time by~\eqref{eq:E2} and the system
$\rho_\nu$ from an already prescribed neighbouring pair, satisfy the
systems $\rho_\nu$ and~\eqref{eq:E2} for all neighbouring pairs: the two
quadruples that~\eqref{eq:E2} gives agree with $\rho_\nu$, and the numbers
obtained do not depend on the neighbouring pair, as verified in
Section~3 of the supplement (\texttt{ehelper.nb})~\cite{Nurmatov2026Repo}.
\end{remark}

A subnet with the central face $F_{\kappa\nu}$ whose flat angles
$\alpha^{\kappa\nu}_i$, $\beta^{\kappa\nu}_i$, $\gamma^{\kappa\nu}_i$,
$\delta^{\kappa\nu}_i$ are the numbers with these names in the assembly,
$i = 1, \ldots, 4$, is called a \emph{realization} of the
block~$(\kappa, \nu)$. For a subnet with the central face $F_{\kappa\nu}$,
we write $\theta^{\kappa\nu}_i$ for its central dihedral angles that
correspond to $\theta_i$, and attach the superscript $\kappa\nu$ in the
same way to its corner dihedral angles. Subnets with the central faces
$F_{\kappa\nu}$ and $F_{\kappa,\nu+1}$, or $F_{\kappa\nu}$ and
$F_{\kappa+1,\nu}$, have six face labels in common, and these faces meet
along the edges at the two shared vertex labels; we say that the two
subnets \emph{agree on the common faces} if their dihedral angles along
these edges coincide.

\begin{lemma}\label{lem:ex:faces}
Assume that $m$, $n \geq 3$ and, for each $\kappa = 1, \ldots, m - 2$ and
$\nu = 1, \ldots, n - 2$, the numbers
$\delta^{\kappa\nu}_1, \ldots, \delta^{\kappa\nu}_4 \in (0, \pi)$ satisfy
$\delta^{\kappa\nu}_1 + \delta^{\kappa\nu}_2 + \delta^{\kappa\nu}_3
+ \delta^{\kappa\nu}_4 = 2\pi$ and
$\delta^{\kappa\nu}_1 + \delta^{\kappa\nu}_2 \geq \pi$. Then there are
convex quadrilaterals, one for each such pair $(\kappa, \nu)$, with the
consecutive vertex labels $V_{\kappa+1,\nu}$, $V_{\kappa\nu}$,
$V_{\kappa,\nu+1}$, $V_{\kappa+1,\nu+1}$ and the angles
$\delta^{\kappa\nu}_1, \ldots, \delta^{\kappa\nu}_4$ at them, respectively,
such that sides with the same labels have the same length.
\end{lemma}

\begin{proof}
First, for a convex quadrilateral $A_1A_2A_3A_4$ with the angles
$\delta_1, \ldots, \delta_4$ at $A_1, \ldots, A_4$, respectively, 
the components of
$\overrightarrow{A_2A_1} + \overrightarrow{A_1A_4} + \overrightarrow{A_4A_3}
+ \overrightarrow{A_3A_2} = 0$ along the normals of $A_4A_1$ and of
$A_3A_4$ give
$|A_3A_4| = \frac{|A_1A_2|\sin{\delta_1} - |A_2A_3|\sin{(\delta_1 + \delta_2)}}{\sin{\delta_4}}$
and
$|A_4A_1| = \frac{|A_2A_3|\sin{\delta_3} - |A_1A_2|\sin{(\delta_2 + \delta_3)}}{\sin{\delta_4}}$.
Conversely, if $\delta_1, \ldots, \delta_4 \in (0, \pi)$ sum to~$2\pi$,
$|A_1A_2|$, $|A_2A_3| > 0$, and $|A_3A_4|$ and $|A_4A_1|$ given by the
formulas are positive, then the walk from $A_2$ through $A_1$, $A_4$, $A_3$ back to $A_2$ along
sides of these four lengths, turning in the same direction by
$\pi - \delta_1$, $\pi - \delta_4$, $\pi - \delta_3$, $\pi - \delta_2$ at
$A_1$, $A_4$, $A_3$, $A_2$, respectively, closes up: by these formulas, its
displacement has zero components along the normals of $A_4A_1$ and of
$A_3A_4$, which are not parallel because $\delta_4 \in (0, \pi)$. Its
turns, each in $(0, \pi)$, sum to $2\pi$, so it gives a convex quadrilateral with the angles
$\delta_1, \ldots, \delta_4$ at $A_1, \ldots, A_4$, respectively.

Second, for each pair $(\kappa, \nu)$, we apply the formulas for
$|A_3A_4|$ and $|A_4A_1|$ to the quadrilateral for $(\kappa, \nu)$, with
$A_1, \ldots, A_4$ labelled $V_{\kappa+1,\nu}$, $V_{\kappa\nu}$,
$V_{\kappa,\nu+1}$, $V_{\kappa+1,\nu+1}$ and the angles
$\delta^{\kappa\nu}_1, \ldots, \delta^{\kappa\nu}_4$ at them, respectively. 
Its upper side $V_{\kappa,\nu+1}V_{\kappa+1,\nu+1}$ is the lower side of the quadrilateral
for $(\kappa, \nu + 1)$ when $\nu < n - 2$, and its right side
$V_{\kappa+1,\nu}V_{\kappa+1,\nu+1}$ is the left side of the one for
$(\kappa + 1, \nu)$ when $\kappa < m - 2$. We choose the lower sides of the
quadrilaterals for $\nu = 1$ positive and then, for $\nu = 1, \ldots, n - 2$
in turn, the left side $V_{1\nu}V_{1,\nu+1}$ of the quadrilateral for
$(1, \nu)$; the formulas give all other sides, row by row from left to
right. Since $\sin{(\delta^{\kappa\nu}_1 + \delta^{\kappa\nu}_2)} \leq 0$,
it follows that the upper side of each quadrilateral is positive once its
lower and left sides are; and, by induction on~$\kappa$, the right side of
the quadrilateral for $(\kappa, \nu)$ is
$|V_{\kappa+1,\nu}V_{\kappa+1,\nu+1}| = \prod_{i=1}^{\kappa}
    \frac{\sin{\delta^{i\nu}_3}}{\sin{\delta^{i\nu}_4}}
    \Bigl(|V_{1\nu}V_{1,\nu+1}| - \sum_{i=1}^{\kappa} |V_{i\nu}V_{i+1,\nu}|\,
    \frac{\sin{(\delta^{i\nu}_2 + \delta^{i\nu}_3)}}{\sin{\delta^{i\nu}_3}}
    \prod_{j=1}^{i-1} \frac{\sin{\delta^{j\nu}_4}}{\sin{\delta^{j\nu}_3}}\Bigr)$,
where the lower sides $|V_{i\nu}V_{i+1,\nu}|$ are already given. Thus,
choosing each $|V_{1\nu}V_{1,\nu+1}|$ large enough, 
in turn for $\nu = 1, \ldots, n - 2$, all sides are positive, and,
by the converse in the first step, the quadrilaterals are convex with the angles
$\delta^{\kappa\nu}_1, \ldots, \delta^{\kappa\nu}_4$; sides with the same
labels have the same length by construction.
\end{proof}

\begin{lemma}\label{lem:ex:agree}
Assume that the flat angles of two subnets with the central faces
$F_{\kappa\nu}$ and $F_{\kappa,\nu+1}$ satisfy the first two equalities
of~\eqref{eq:E2} and the hypotheses of Theorem~\ref{th:main} with the
signs $e^{\kappa\nu}_1, \ldots, e^{\kappa\nu}_4$ and
$e^{\kappa,\nu+1}_1, \ldots, e^{\kappa,\nu+1}_4$, respectively, where
$e^{\kappa\nu}_3e^{\kappa\nu}_4 = e^{\kappa,\nu+1}_1e^{\kappa,\nu+1}_2$.
Then, for $e_0 = \pm 1$ and all but finitely many $t$ as in
Theorem~\ref{th:main} applied to the flat angles of the second subnet, the
polyhedra of the flexions it gives with the parameter~$t$ and the
sign~$e_0$ agree on the common faces with those of the flexions given by
Theorem~\ref{th:main} applied to the flat angles of the first subnet after
the cyclic permutation of the indices by~$2$, with the parameter~$t$ and
the sign $-e_0e^{\kappa\nu}_3e^{\kappa,\nu+1}_2$.
\end{lemma}

\begin{proof}
After the cyclic permutation of the indices by~$2$ of the first subnet,
its vertices $A_1$ and $A_2$ are the vertices $A_2$ and $A_1$ of the
second, respectively, so the common edge of the central faces is $A_1A_2$
in both, and the signs of the first are $e^{\kappa\nu}_3$,
$e^{\kappa\nu}_4$, $e^{\kappa\nu}_1$, $e^{\kappa\nu}_2$ by
Remark~\ref{rem:th:main:2}. If $e^{\kappa,\nu+1}_2 = -e^{\kappa\nu}_3$, then
$e^{\kappa,\nu+1}_1 = -e^{\kappa\nu}_4$, because
$e^{\kappa\nu}_3e^{\kappa\nu}_4 = e^{\kappa,\nu+1}_1e^{\kappa,\nu+1}_2$,
and the sign of the first is
$-e_0e^{\kappa\nu}_3e^{\kappa,\nu+1}_2 = e_0$. Otherwise,
$e^{\kappa,\nu+1}_1 = e^{\kappa\nu}_4$,
$e^{\kappa,\nu+1}_2 = e^{\kappa\nu}_3$, and the sign of the first is
$-e_0$. In the latter case, we replace
$(e_0, e^{\kappa,\nu+1}_1, \ldots, e^{\kappa,\nu+1}_4)$ by
$(-e_0, -e^{\kappa,\nu+1}_1, \ldots, -e^{\kappa,\nu+1}_4)$. This preserves
the products $e_0e_1$, $e_0e_2$, $e_2e_3$, through which alone the signs
enter the formulas of Theorem~\ref{th:main}, and condition~\eqref{eq:N6},
because $\Lambda = -\Lambda$; hence it preserves the polyhedra of the
flexions for the flat angles of the second subnet. After the replacement,
$e^{\kappa,\nu+1}_1 = -e^{\kappa\nu}_4$,
$e^{\kappa,\nu+1}_2 = -e^{\kappa\nu}_3$, and these polyhedra of the second
have the sign $-e_0$, the sign of the first. Thus, in both cases,
$e^{\kappa,\nu+1}_1 = -e^{\kappa\nu}_4$,
$e^{\kappa,\nu+1}_2 = -e^{\kappa\nu}_3$, and the polyhedra of the two
have the same parameter and sign. By~\eqref{eq:E2},
\eqref{eq:N5}, Theorem~\ref{th:main}, and Remark~\ref{rem:th:main:1}, in
these numberings the two subnets have the same $u$ and $x_1$, so their
parameters range over the same admissible set $I(x_1)$, and, for all but
finitely many~$t$ in it, such polyhedra agree on the common faces, as
verified in Section~4 of the supplement
(\texttt{ehelper.nb})~\cite{Nurmatov2026Repo}.
\end{proof}

\begin{remark}\label{rem:ex:agree}
Lemma~\ref{lem:ex:agree} also holds for the subnets with the central faces
$F_{\kappa\nu}$ and $F_{\kappa+1,\nu}$, with the last two equalities
of~\eqref{eq:E2} in place of the first two, the condition
$e^{\kappa\nu}_4e^{\kappa\nu}_1 = e^{\kappa+1,\nu}_2e^{\kappa+1,\nu}_3$,
Theorem~\ref{th:main} applied to the flat angles of the second subnet after
the cyclic permutation of the indices by~$1$ and to those of the first
after the cyclic permutation by~$3$, and the signs $e_0$ for the second and
$-e_0e^{\kappa\nu}_4e^{\kappa+1,\nu}_3$ for the first. Indeed, the proof of
Lemma~\ref{lem:ex:agree} applies with these permutations, as verified for
them in Section~4 of the supplement
(\texttt{ehelper.nb})~\cite{Nurmatov2026Repo}.
\end{remark}

\begin{proposition}\label{prop:ex:flex}
Assume that $\alpha_1$, $\beta_1$, $\gamma_1$,
$\delta_1 \in (0, \pi)$ satisfy condition~\eqref{eq:N1},
$\overline{\alpha}_1$, $\overline{\beta}_1$, $\overline{\gamma}_1$,
$\overline{\delta}_1 \in (0, \pi)$, $m$, $n \geq 3$, $1 \leq \kappa_0 \leq m - 2$,
$1 \leq \nu_0 \leq n - 2$, and each of
$\rho_1, \ldots, \rho_{n-2}$ is one of the systems (a)--(d) of
Lemma~\ref{lem:ex:patterns}.
Then there is a flexible $m \times n$ net with the flat angles prescribed
by the $(\rho_1, \ldots, \rho_{n-2})$-assembly, and a flexion of it is
given on each of its subnets by the formulas of Theorem~\ref{th:main} for the
signs $(e_1, e_2, e_3, e_4) = (1, -1, 1, -1)$.
\end{proposition}

\begin{proof}
First, for each block~$(\kappa, \nu)$, there is a polyhedron with the flat
angles $\alpha^{\kappa\nu}_i$, $\beta^{\kappa\nu}_i$,
$\gamma^{\kappa\nu}_i$, $\delta^{\kappa\nu}_i$, $i = 1, \ldots, 4$, given
by the assembly and with a flexion given by Theorem~\ref{th:main} for the
signs $(1, -1, 1, -1)$, that is, a flexible realization of the
block~$(\kappa, \nu)$. Indeed, every quadruple of the assembly is obtained
from $(\alpha_1, \beta_1, \gamma_1, \delta_1)$ by the permutations of the
entries in~\eqref{eq:E2} and the permutations and complements to~$\pi$ in
the systems (a)--(d), which preserve condition~\eqref{eq:N1} and keep the
quadruple and its barred quadruple in $(0, \pi)^4$, see Section~2 of the
supplement (\texttt{ehelper.nb})~\cite{Nurmatov2026Repo}. 
Hence, by Lemma~\ref{lem:ex:patterns}, these flat angles satisfy
conditions~\eqref{eq:N1}, \eqref{eq:N4}, and~\eqref{eq:N5} for each
$i = 1, \ldots, 4$, $\delta_1 + \delta_2 + \delta_3 + \delta_4 = 2\pi$,
and $M > 0$. Then, by Lemmas~2 and~4 of~\cite{NURMATOV2026104102}, $r_i$,
$s_i$, $f_i$, $M \neq 0, 1$ are well-defined positive reals.
In each of the systems (a)--(d), the quadruples at
$A_2$ and $A_3$ are the complements to~$\pi$ of those at $A_1$ and $A_4$,
respectively, and taking complements preserves $r_i$, $s_i$, $f_i$ and
replaces $\sigma_i$ by $2\pi - \sigma_i$, see~\eqref{eq:N2}
and~\eqref{eq:N3}. Thus $t_1 = t_2 \pm 2K$ and $t_3 = t_4 \pm 2K$
by~\eqref{eq:N10}, \eqref{eq:N11}, and Table~\ref{tab:tpqsigma}. Hence
$t_1 - t_2 + t_3 - t_4 \in \{0, \pm 4K\} \subset \Lambda$
by~\eqref{eq:N9}, which is condition~\eqref{eq:N6} with
$(e_1, e_2, e_3, e_4) = (1, -1, 1, -1)$, and Theorem~\ref{th:main}
applies.

Second, we fix the faces. In each of the systems (a)--(d),
$\delta^{\kappa\nu}_1 + \delta^{\kappa\nu}_2
= \delta^{\kappa\nu}_3 + \delta^{\kappa\nu}_4 = \pi$, and hence
Lemma~\ref{lem:ex:faces} applies to the numbers $\delta^{\kappa\nu}_i$ of
the assembly. We fix the quadrilaterals it gives as the faces
$F_{\kappa\nu}$ and, for the other labels $F_{ij}$, quadrilaterals with the
angles given by the assembly at their inner vertex labels, that is, the labels $V_{ij}$ with
$0 < i < m$ and $0 < j < n$, and short enough edges 
from the inner vertex labels to the others, completed by the
vertices labelled $V_{00}$, $V_{m0}$, $V_{0n}$, $V_{mn}$, such that sides
with the same labels have the same length; they can be chosen convex,
because the angles given by the assembly lie in $(0, \pi)$. Every 
polyhedron of a flexion given by Theorem~\ref{th:main} in the first step
can be rebuilt from these quadrilaterals with the same
flat and dihedral angles: the fixed faces along the edges of the fixed
central face are attached to it at the dihedral angles of the polyhedron
along these edges, which is always possible,
and Bricard's system~\eqref{eq:N13}, which holds by
 Lemma~\ref{lem:bricard}, gives, in the notation of the $3 \times 3$ net,
$\angle B_iA_iC_i = \beta_i$ for $i = 1, \ldots, 4$,
so the fixed corner faces fit, as verified in
Section~5 of the supplement (\texttt{ehelper.nb})~\cite{Nurmatov2026Repo}.

Third, we realize the blocks: $(1, 1)$ at a parameter~$t$ and a
sign~$e_0$, then $(\kappa + 1, 1)$ from $(\kappa, 1)$ for
$\kappa = 1, \ldots, m - 3$, and then $(\kappa, \nu + 1)$ from
$(\kappa, \nu)$ for $\kappa = 1, \ldots, m - 2$ and
$\nu = 1, \ldots, n - 3$, each with the parameter and sign given by
Remark~\ref{rem:ex:agree} or Lemma~\ref{lem:ex:agree}, respectively, whose
hypotheses hold by the first step and~\eqref{eq:E2}, both sides of their
sign conditions being~$-1$, and all of them rebuilt from the fixed
quadrilaterals of the second step; here, by Remark~\ref{rem:th:main:2},
every realization is a polyhedron of the flexions of Theorem~\ref{th:main}
in each numbering of the remark and the lemma, with some parameter and sign. 
Then the realizations of a block and
of the one realized from it agree on
the common faces, and those of all neighbours agree once those of
$(\kappa, \nu)$ and $(\kappa + 1, \nu)$ do for $\kappa = 1, \ldots, m - 3$
and $\nu = 1, \ldots, n - 2$. This holds for $\nu = 1$, one of the two
being realized from the other, and, if it holds for some
$\nu \leq n - 3$, then the realization of $(\kappa + 1, \nu + 1)$ from
$(\kappa + 1, \nu)$ coincides with its realization from
$(\kappa, \nu + 1)$, because, by these agreements, 
both have the same $\theta_1$ and $\theta_4$ at
the vertex label $V_{\kappa+1,\nu+1}$ in the numbering of
Remark~\ref{rem:ex:agree}, hence the same parameter and sign, as verified in
Section~6 of the supplement (\texttt{ehelper.nb})~\cite{Nurmatov2026Repo},
using Theorem~\ref{th:main} and the strict positivity of its radicands, and
thus it holds for $\nu + 1$. Each realization and each step of this induction excludes
finitely many values of a dihedral angle along an edge of a central face.
Since $c_{11i} = -\sin{\alpha_i}\sin{\gamma_i} \neq 0$, and
$c_{00i} \neq 0$ by~\eqref{eq:N1}, it follows that each equation
of~\eqref{eq:N13} has at most two solutions in either variable for a fixed
value of the other; hence, back along the order of realization, each
excluded value comes from only finitely many values of~$t$, and we
restrict $t$ to an interval free of them. The realizations of the block~$(1, 1)$
are pairwise non-congruent, being polyhedra of one flexion of
Theorem~\ref{th:main}, and, by induction along the order of realization,
all realizations depend continuously on~$t$,
because the formulas of Theorem~\ref{th:main} do and a realization cannot
switch with~$t$ between the two signs of Theorem~\ref{th:main}, which give
different polyhedra for the same parameter, as verified in Section~6 of
the supplement (\texttt{ehelper.nb})~\cite{Nurmatov2026Repo}.

Finally, every edge between two faces lies at a vertex of a central face,
and all realizations containing both faces give it the same dihedral
angle, because their blocks form a rectangle in the grid of blocks and the
realizations of neighbours agree. Placing the fixed quadrilaterals of the 
second step one by one along the
sides with the same labels, at these dihedral angles, gives an
$m \times n$ net whose subnets are realizations of their blocks,
independently of the order, because the faces around an inner vertex close
up, as in a realization containing them, and the grid of faces is a disc.
These nets depend continuously on~$t$ and are pairwise non-congruent,
because the realizations do and those of the block~$(1, 1)$ are.
Hence they form a flexion, which is given on each of their subnets by the
formulas of Theorem~\ref{th:main} for the signs $(1, -1, 1, -1)$.
\end{proof}

\begin{remark}\label{rem:ex:flex:1}
Proposition~\ref{prop:ex:flex} also holds with the signs
$(1, -1, -1, 1)$ in place of the signs $(1, -1, 1, -1)$. Indeed, its proof applies:
in the first step, $t_1 - t_2 - t_3 + t_4 \in \{0, \pm 4K\} \subset \Lambda$,
which is condition~\eqref{eq:N6} for these signs, and both sides of the
sign conditions of Lemma~\ref{lem:ex:agree} and Remark~\ref{rem:ex:agree}
are $-1$ and~$1$, respectively.
\end{remark}

\begin{remark}\label{rem:ex:flex:2}
The central faces of the nets of Proposition~\ref{prop:ex:flex} are
trapezoids, because $\delta_1 + \delta_2 = \delta_3 + \delta_4 = \pi$ in each
of the systems (a)--(d). Nevertheless, their subnets are not T-nets of
Sauer and Graf~\cite{SauerGraf1931} in general. Indeed, the flat angles of
a $3 \times 3$ T-net satisfy
$\cos{\alpha_i}\cos{\gamma_i} = \cos{\beta_i}\cos{\delta_i}$,
$i = 1, \ldots, 4$, see~\cite[Section~3.1.1]{izmestiev-2017}, but for
$(\alpha_1, \beta_1, \gamma_1, \delta_1) = (\pi/3, \pi/3, \pi/2, 2\pi/3)$,
which satisfies the hypotheses of Proposition~\ref{prop:ex:flex}, we have
$|\cos{\alpha^{\kappa\nu}_i}\cos{\gamma^{\kappa\nu}_i}
- \cos{\beta^{\kappa\nu}_i}\cos{\delta^{\kappa\nu}_i}| = 1/4$ for all
$\kappa$, $\nu$, and~$i$, because every quadruple of the assembly is
obtained from $(\alpha_1, \beta_1, \gamma_1, \delta_1)$ by permutations and
complements to~$\pi$. Hence no subnet is a T-net.
\end{remark}

\begin{figure}[htbp!]
    \centering
    \includegraphics[width=0.9\textwidth]{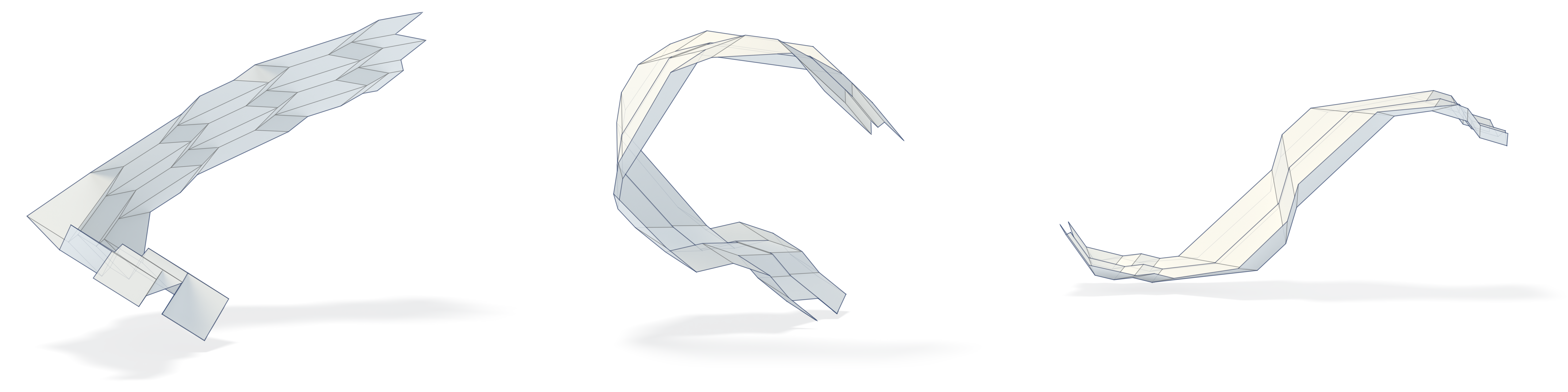}\\[-0.2ex]
    {\footnotesize%
    \makebox[0.3\textwidth]{(i)~$t = 1.12$}%
    \makebox[0.3\textwidth]{(ii)~$t = 0.38$}%
    \makebox[0.3\textwidth]{(iii)~$t = 0.44$}}
    \\[1.5ex]
    \includegraphics[width=0.9\textwidth]{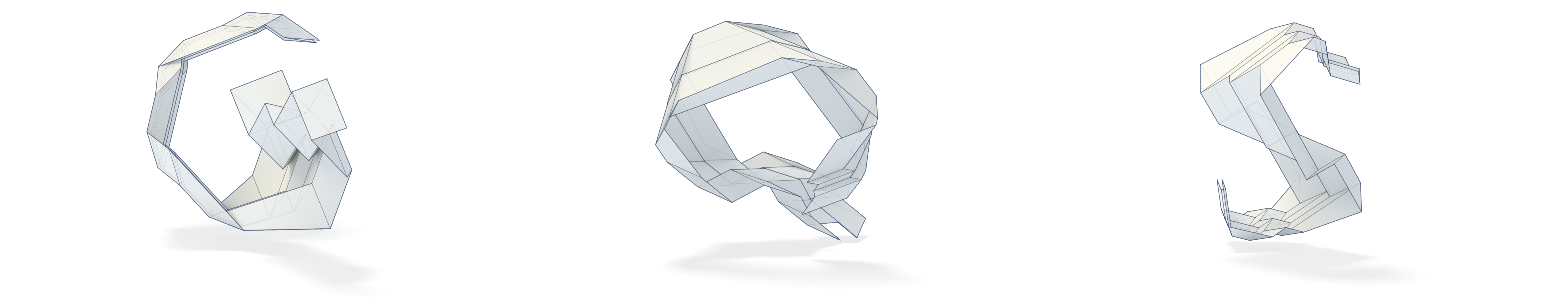}\\[-0.2ex]
    {\footnotesize%
    \makebox[0.3\textwidth]{(i)~$t = 6.25$}%
    \makebox[0.3\textwidth]{(ii)~$t = 0.64$}%
    \makebox[0.3\textwidth]{(iii)~$t = 0.78$}}
    \caption{Three flexible $5 \times 12$ nets (i)--(iii) assembled from
    GQS-nets, each at two values of the parameter~$t$ of the
    block~$(\kappa_0, \nu_0)$, from the same viewing direction and at the
    same scale. Bottom: the nets read GQS; top: another position of the same
    nets. Table~\ref{tab:ex:nets:GQS} gives the data of the nets.}
    \label{fig:ex:nets:GQS}
\end{figure}

\subsection{Algorithm and examples}
\label{sec:ex:algorithm}
To build concrete examples, we use the construction in the proof of
Proposition~\ref{prop:ex:flex}. We choose numbers $\alpha_1$, $\beta_1$,
$\gamma_1$, $\delta_1 \in (0, \pi)$ that satisfy condition~\eqref{eq:N1}
and $\overline{\alpha}_1$, $\overline{\beta}_1$, $\overline{\gamma}_1$,
$\overline{\delta}_1 \in (0, \pi)$, integers $m$, $n \geq 3$, indices
$1 \leq \kappa_0 \leq m - 2$ and $1 \leq \nu_0 \leq n - 2$, systems
$\rho_1, \ldots, \rho_{n-2}$ among (a)--(d) of
Lemma~\ref{lem:ex:patterns}, the signs $(1, -1, 1, -1)$ or, by
Remark~\ref{rem:ex:flex:1}, $(1, -1, -1, 1)$, and a sign~$e_0$. Then we
compute the numbers of the $(\rho_1, \ldots, \rho_{n-2})$-assembly and fix
the faces as in the second step of the proof, with the sides of the
central faces chosen as in the proof of Lemma~\ref{lem:ex:faces}. As in the
third step, but starting from the block~$(\kappa_0, \nu_0)$, we realize
this block at a parameter~$t$ and the sign~$e_0$, then the other blocks of
its row and then those of each column, each from a realized neighbour, in
either direction, with the parameter and sign related to those of the
neighbour by Remark~\ref{rem:ex:agree} or Lemma~\ref{lem:ex:agree},
respectively. Finally, we place the fixed faces at the dihedral angles of
these realizations, as in the last step. 
Because the relations of Remark~\ref{rem:ex:agree} and
Lemma~\ref{lem:ex:agree} determine each of the two realizations of
neighbours from the other, the nets obtained in this way are those of the
proof of Proposition~\ref{prop:ex:flex}, with the block~$(1, 1)$ realized
at the parameter and sign so determined. 
Hence, by that proof, the realizations of all neighbours agree on the
common faces, the faces of the net are planar convex quadrilaterals with
the flat angles of the assembly at the inner vertices, and its dihedral
angles along the edges of the central faces are those given by
Theorem~\ref{th:main}. The supplement (\texttt{gqs.py})~\cite{Nurmatov2026Repo}
implements this algorithm. At every position drawn, it verifies these
properties numerically, with the flat and dihedral angles recomputed from
the vertices of the net, and checks that the faces keep their shape; it also
detects self-intersections in every frame of an exported motion.
Figure~\ref{fig:ex:nets:GQS} shows three of the nets obtained with it, for
$(\alpha_1, \beta_1, \gamma_1, \delta_1) = (45^\circ, 45^\circ, 90^\circ, 150^\circ)$,
$(142.5^\circ, 127.5^\circ, 15^\circ, 60^\circ)$, and
$(75^\circ, 22.5^\circ, 112.5^\circ, 120^\circ)$
($M = 4/3$, $\sqrt{3}/2$, and~$\sqrt{3}/2$), respectively, each at two
positions of its motion; two of them have the signs $(1, -1, -1, 1)$
of Remark~\ref{rem:ex:flex:1}. Section~\ref{sm:ex} gives their data and
further examples: four nets that form the digits ``2026'' and two
$6 \times 60$ nets whose systems~$\rho_\nu$, chosen by a beam search, make
them approximately follow the side silhouettes of the Stanford
bunny~\cite{StanfordRepo} and of Blub the fish~\cite{CraneRepo}.

\section{Application}
\label{sec:application}
\begin{figure}[htbp!]
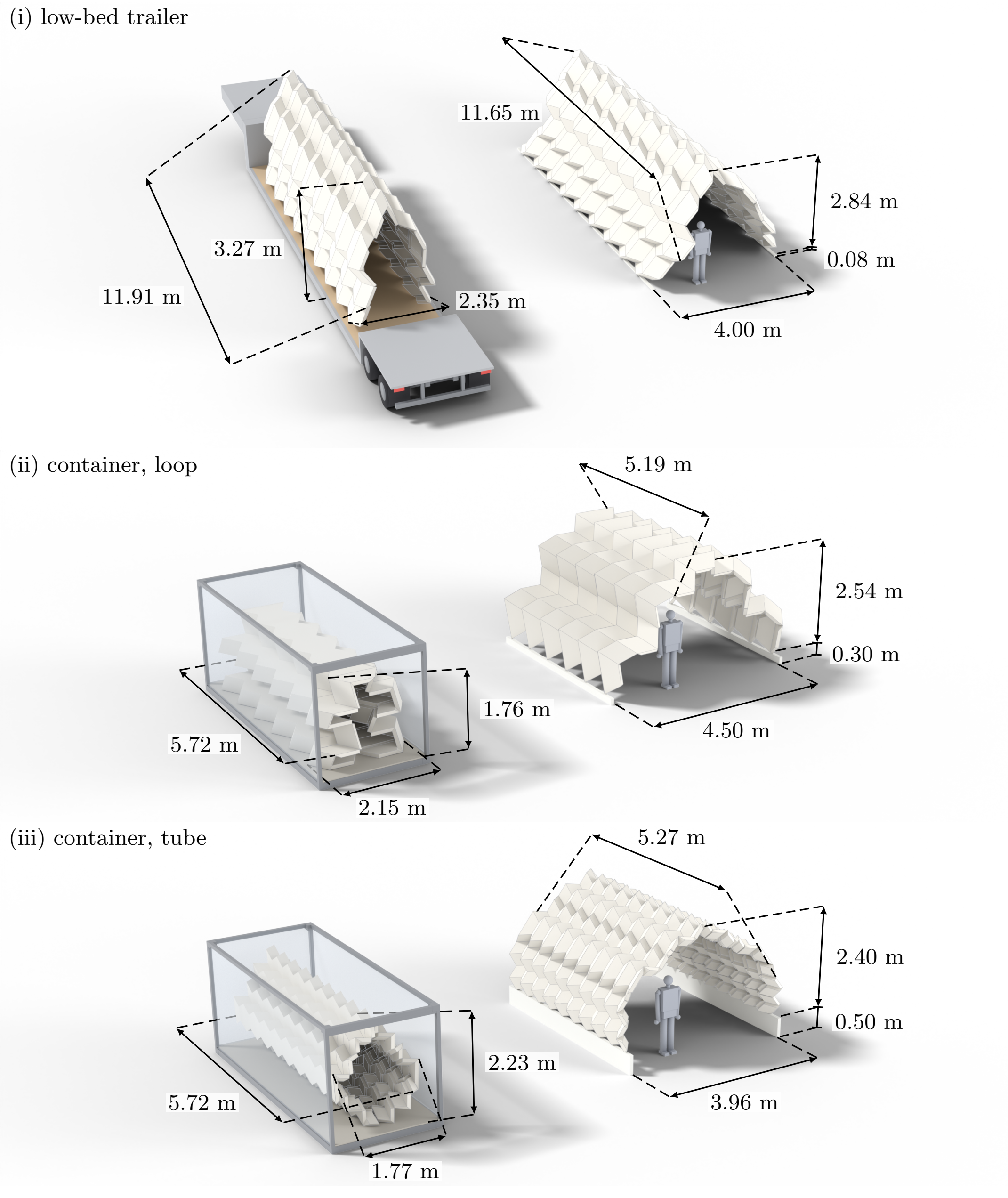

    \centering
    \input{fig/fig_app_render_H3s.tex}\\[0pt]
    \input{fig/fig_app_render_P1.tex}\\[0pt]
    \input{fig/fig_app_render_B.tex}
    \caption{The halls (i)--(iii) in the transport (left)
and deployed (right) positions, with a person $1.80$~m tall. Values,
left: width, height, and length; right: span, rise to the crown of the
midline of the profile, length, and the height of the feet above the floor.
The hall~(i) stands on plinths and (ii) and~(iii) on knee walls. Schematic
trailer~\cite{KSLL2} and 1CC container~\cite{ISO668};
Section~\ref{sm:app:render} gives the further details of
the rendering.}
    \label{fig:app:render}
\end{figure}

We now apply the nets of Proposition~\ref{prop:ex:flex} to a design
problem with prescribed dimensions and loads: a transportable arch hall
(Figure~\ref{fig:app:render}), shipped narrow and expanded on site.
In its transport position, the hall must pass through the door of a
20-ft container of type 1CC, at least $2.286$~m wide and $2.261$~m
high~\cite{ISO668}, or travel on a semi-trailer within the width and height
limits of $2.55$ and $4.00$~m of~\cite[Annex~I, 1.2(a) and 1.3]{EU96-53}.
Deployed, it should span $3.5$ to $4.5$~m and give, inside the panels,
a clear height of at least $2.5$~m and a floor strip at least $2.0$~m wide
with at least $2.0$~m of clear height. On their inner side, the faces carry $80$~mm sandwich
panels with steel face sheets and a polyurethane core, hinged along the
edges of the net. Section~\ref{sm:app} gives the details of this section,
and the supplement (\texttt{Application})~\cite{Nurmatov2026Repo} contains
its data and checks.

The halls (i)--(iii) consist of faces of the $20 \times 20$,
$13 \times 19$, and $20 \times 27$ nets of Proposition~\ref{prop:ex:flex}
with $(\alpha_1, \beta_1, \gamma_1, \delta_1) = (36^\circ, 144^\circ, 72^\circ, 72^\circ)$,
$(120^\circ, 60^\circ, 60^\circ, 60^\circ)$, and
$(60^\circ, 60^\circ, 120^\circ, 60^\circ)$
($M = 5/4$, $9/8$, and~$9/8$), respectively, the signs
$(e_1, e_2, e_3, e_4) = (1, -1, 1, -1)$, and $e_0 = 1$; Table~\ref{tab:app} gives their
blocks~$(\kappa_0, \nu_0)$, systems~$\rho_\nu$, faces, and positions. We omit the boundary
faces of each net, whose edges to the boundary vertices the proof of
Proposition~\ref{prop:ex:flex} leaves free, so that each hall ends on
regular rows of faces; for~(i), we also omit the faces with $j = 1$, so
that its two feet are mirror images. Omitting such faces preserves the flexion
of the net, which moves the remaining faces rigidly and keeps their common
edges.

\begin{figure}[htbp!]
    \centering
    \includegraphics[width=0.44\textwidth]{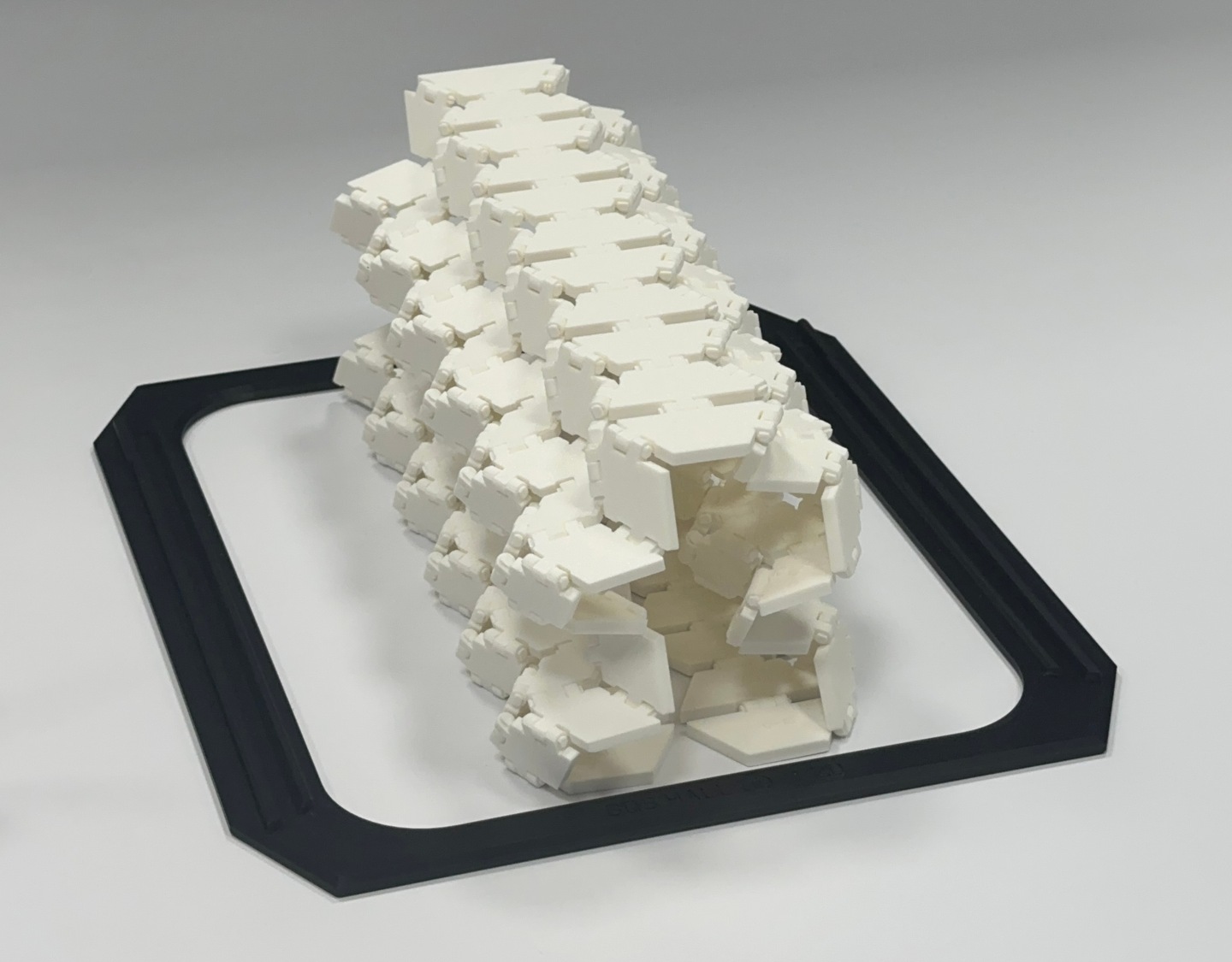}\hspace{0.02\textwidth}%
    \includegraphics[width=0.44\textwidth]{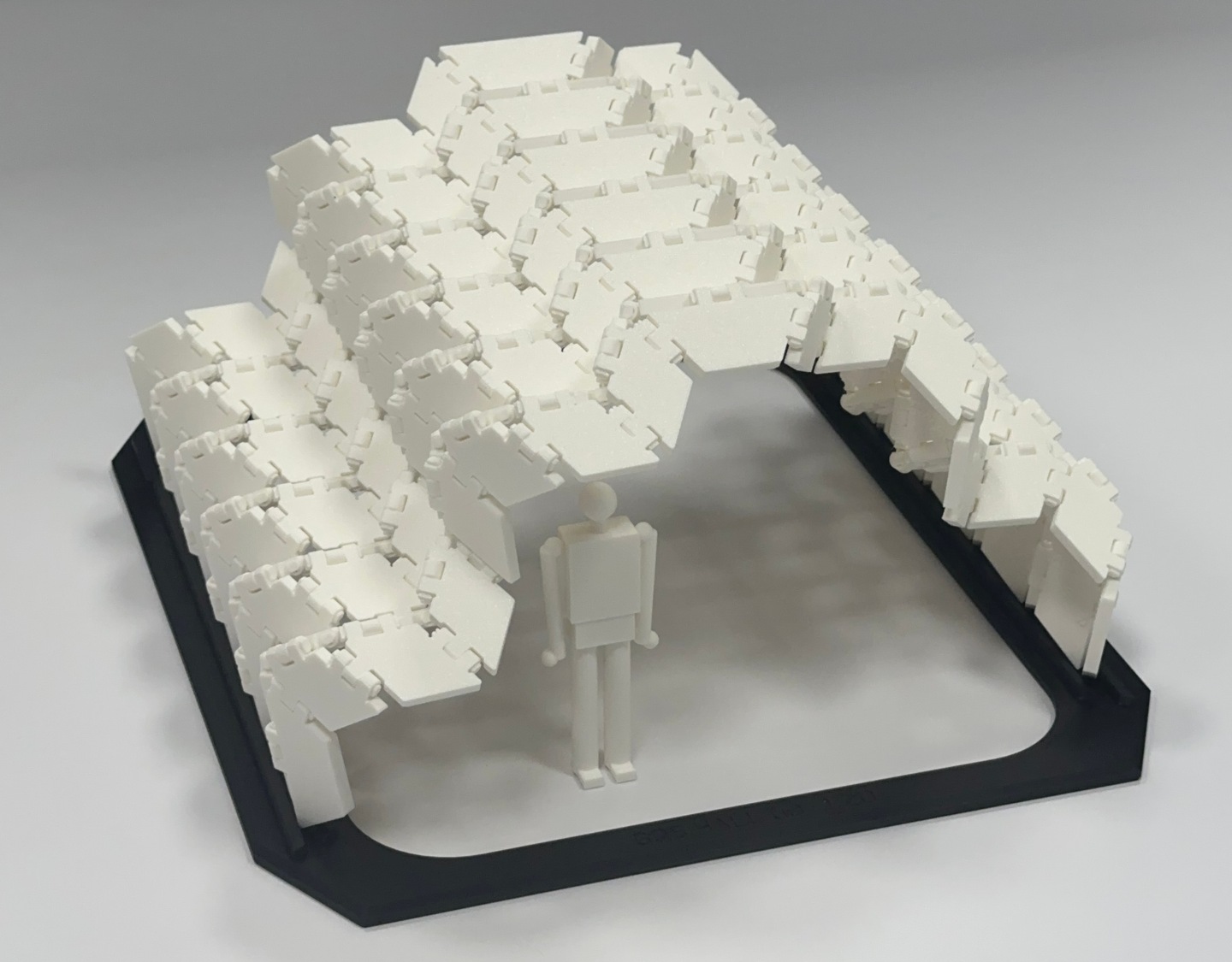}%
    \caption{Photographs of the $1{:}20$ physical model of the hall~(ii),
    folded (left) and deployed on its stand (right), with the person of
    Figure~\ref{fig:app:render} at $1{:}20$.}
    \label{fig:app:model}
\end{figure}

Seen along the mean direction of its rows of blocks, a long net looks
like a planar curve, which we call its \emph{profile}. We measure its
length on its \emph{midline}, the mean of the projected polylines of the
vertices~$V_{ij}$ with even and with odd~$i$. This length changes by only
about $3\%$ between the deployed and the transport positions of the halls
(Table~\ref{tab:app}). Hence a hall cannot pack by shortening its profile.
Instead, the transport position folds the profile: with its two legs side
by side for~(i), which needs a low-bed semi-trailer, and into a nearly
closed loop or a tube for (ii) and~(iii), which pass through the door, see
Figure~\ref{fig:app:render} and the profiles in Figure~\ref{fig:app:profiles}.
A $1{:}20$ physical model of~(ii) folds by hand into such a loop and opens
again (Figure~\ref{fig:app:model} and Section~\ref{sm:app:model}).

We found the systems~$\rho_\nu$ by searches: at one value of the
parameter~$t$ of the block~$(\kappa_0, \nu_0)$, the
profile must follow a parabola or a circular arc of the given span, and at
a second, the net must fit the door or the trailer.
The motion between the two values of~$t$ is
checked numerically, at finitely many positions, for self-intersections and
collisions of the panels, see Section~\ref{sm:app:nets}. The motion
of~(iii) passes through $t = 1/\sqrt{2}$, which Theorem~\ref{th:main}
excludes (Section~\ref{sm:app:drive}).

In the transport position, (ii) and~(iii) pass through the door and fit the
floor of the container, and (i) stays within the width and height limits
in the well of an extendable low-bed semi-trailer~\cite{KSLL2}, reaching
$3.57$~m above the road; extended for this hall, 
the trailer exceeds the length limits and may circulate on a special
permit or under a national derogation~\cite{EU96-53}, see
Section~\ref{sm:app:transport}. Deployed, on plinths for~(i) and on knee
walls for (ii) and~(iii), the halls have the largest clear heights $2.74$,
$2.71$, and $2.79$~m and floor strips $2.00$, $2.29$, and $2.30$~m wide
with at least $2.0$~m of clear height.

We check the deployed halls with a linear shell model of the sandwich
under the self-weight and the undrifted and drifted snow on a cylindrical
roof of~\cite[Section~5.3.5]{EN1991-1-3} for a ground snow load
$s_{\mathrm{k}} = 1$~kN/m$^2$. The capacities, the smallest ground snow
loads at which a limit on displacement or stress is reached, with the
stresses read outside small zones around the point connections, are $0.60$,
$0.64$, and $0.68$~kN/m$^2$ with both feet anchored, about two thirds
of~$s_{\mathrm{k}}$, set by the sway of the open ends of~(i) and by the
wrinkling of the face sheets of (ii) and~(iii). With one foot on rails, a
spindle between the feet, and one locked hinge, as in the drive of
Section~\ref{sm:app:drive}, they are $0.43$, $0.04$, and $0.04$~kN/m$^2$,
set by wrinkling along a hinge line of~(i) and by the spreading of the foot
on the rails of (ii) and~(iii). In both cases, the buckling factors of the
global modes exceed~$1$. Each limit points to a measure to be examined: a
beam or carriage along the foot on the rails, frames at the open ends, and
hinge hardware supporting the face sheets. Section~\ref{sm:app:snow} gives
the shell model, its limits and assumptions, and the mesh study.

\section{Discussion}
\label{sec:discussion}
We now return to the three tasks of Section~\ref{sec:intro} and to
Izmestiev's second and third questions recalled there, and then discuss directions for
further work and applications.

First, the conditions of Propositions~\ref{prop:eqMrs}, \ref{prop:shifts},
and~\ref{prop:signs} can be checked directly on given flat angles, for example on
those of the subnets of the numerically computed mechanisms
of~\cite[Section~3]{Yorov2026}; for such angles, they can be confirmed
only up to numerical error.

Second, Theorem~\ref{th:main} answers Izmestiev's second question for
this type, on the real part of the
configuration space: all numbers that satisfy its hypotheses are the
flat angles of a flexible polyhedron, and hence, for such numbers, this
real part is nonempty. The
third question, on self-intersections, is addressed here only
numerically: the supplement (\texttt{gqs.py})~\cite{Nurmatov2026Repo}
tests the nets of Sections~\ref{sec:examples} and~\ref{sm:ex} for self-intersections at
finitely many positions of their motions, and
Section~\ref{sec:application} checks the halls, with their panels, in
the same way. For rigid-foldable
quadrilateral meshes, no efficient way to predict self-intersection
other than simulation is known~\cite[Section~8(c)]{HeGuest2020}. Here,
however, the flexions are in closed form: by the proof of
Proposition~\ref{prop:ex:flex} and the construction of the examples in
Section~\ref{sec:examples}, the
central dihedral angles of all subnets of a net are explicit functions
of the parameter~$t$ of the block~$(\kappa_0, \nu_0)$, and, together
with the faces, they determine the net up to a rigid motion. 
Thus self-intersections along an interval of~$t$
could be excluded exactly, by locating the values of~$t$ at which two
faces touch outside their common edge or vertex; such a method remains to be developed.

Third, choosing neighbouring $3 \times 3$ subnets is settled here only for
the simplest GQS-nets. The assembly of Section~\ref{sec:examples} uses
the GQS-nets of Lemma~\ref{lem:ex:patterns}, with the pairs of indices
$\{2, 3\}$ and $\{1, 4\}$ of Definition~\ref{def:gqs}, one of the
systems (a)--(d) for each row of blocks, and the same signs,
$(1, -1, 1, -1)$ or $(1, -1, -1, 1)$, on all subnets
(Proposition~\ref{prop:ex:flex} and Remark~\ref{rem:ex:flex:1}). How to
build larger meshes from $3 \times 3$ blocks is open in
general~\cite[Section~5]{Jiang2024}.

Each of these restrictions leaves room for further work: the pairs of
indices $\{1, 3\}$ and $\{2, 4\}$ call for a different approach, see
Section~\ref{sec:examples}; the systems (a)--(d) of
Lemma~\ref{lem:ex:patterns} are only the simplest solutions we found
of~\eqref{eq:E1} and of its three analogues; and assemblies with signs
that change between subnets remain to be studied.
Beyond GQS-nets, the subnets
may be any polyhedra of equimodular elliptic type. For such subnets, the
first two tasks are settled, and the third remains, with its three
requirements. Sharing flat angles means solving~\eqref{eq:E2} together
with the hypotheses of Theorem~\ref{th:main} in every block.
Lemma~\ref{lem:ex:faces} gives central faces that fit together if
$\delta_1 + \delta_2 \geq \pi$ in every block, and
Lemma~\ref{lem:ex:agree} and Remark~\ref{rem:ex:agree} show that two
neighbouring subnets move compatibly under a condition on their signs.
If this condition holds for all neighbours, then the induction in the
third step of the proof of Proposition~\ref{prop:ex:flex}, which does not
use the systems (a)--(d), shows that the subnets also move compatibly
around every $2 \times 2$ square of blocks. It remains to meet all these
conditions together.

For shape design, the bunny and fish nets of Section~\ref{sm:ex:animals} and the halls
of Section~\ref{sec:application} follow a prescribed curve only
approximately and only at one value of~$t$. Two goals go further: nets with
subnets of equimodular elliptic type that approximate whole
three-dimensional shapes, and nets whose flexion passes from one prescribed
shape to another, for instance from the bunny to the fish. Both could first
be sought among the nets of Proposition~\ref{prop:ex:flex}, for the second
by choosing the systems~$\rho_\nu$ so that the net matches two silhouettes
at two values of~$t$; the halls solve a simpler version, with only a bound
on the size of the net at the second value. A surface, however, also needs freedom along the rows of
blocks, where the systems~$\rho_\nu$ do not change. This freedom could
come from the side lengths of the faces, left partly free by the proofs
of Lemma~\ref{lem:ex:faces} and Proposition~\ref{prop:ex:flex}, and from
other subnets of equimodular elliptic type. The approximation of given
shapes by flexible quadrilateral meshes is also open in
general~\cite[Section~5]{Jiang2024}.

Finally, unlike the facades and roofs that pass from a flat to a curved
state~\cite{Tachi2010freeform}, the nets of
Proposition~\ref{prop:ex:flex} have no flat position.
Indeed, every inner vertex is a vertex of the central face of a subnet,
whose flat angles satisfy condition~\eqref{eq:N1} by the first step of
the proof of Proposition~\ref{prop:ex:flex}; in a flat position, a
signed sum of the flat angles of the four faces at this vertex would
vanish modulo~$2\pi$, which is ruled out by condition~\eqref{eq:N1}.
Since the developed and flat-folded states
of~\cite[Section~2.3]{Tachi2010freeform} are flat positions, these nets
are neither developable, as the metamaterials and curved shapes
of~\cite{Silverberg2014, Dudte2016, Callens2018, Dieleman2020} are, nor
flat-foldable for compact packaging~\cite[Section~1]{Tachi2010freeform}.
The argument needs only the four faces at one inner vertex; hence the
arch halls of Section~\ref{sec:application}, formed by faces of such
nets, have no flat position either. The nets thus suit uses that need a
motion of rigid faces but no flat position, as these halls do, for
instance envelopes, partitions, facades, and shades in transformable
architecture with thick rigid panels and
hinges~\cite[Sections~1 and~5]{Tachi2010freeform}. For such uses, they
offer flexions in closed form for all $m, n \geq 3$, flat angles
determined by those at one vertex and by the systems~$\rho_\nu$, and
partly free side lengths.

\section{Conclusion}
\label{sec:concl}
We have given a constructive approach to the equimodular elliptic type
and to flexible $m \times n$ nets built from its polyhedra. Elementary
algebra now decides whether given flat angles have the type and for
which signs its third condition holds
(Propositions~\ref{prop:eqMrs}, \ref{prop:shifts}, and~\ref{prop:signs}). 
The algebraic form
of this condition leads to the GQS-nets, which have the type and
include the elliptic QS-nets (Definition~\ref{def:gqs} and
Theorem~\ref{th:gqs}). Theorem~\ref{th:main} realizes the flat angles of
every polyhedron of the type by a flexible polyhedron, with its flexion
in closed form. From the simplest GQS-nets,
Proposition~\ref{prop:ex:flex} assembles flexible $m \times n$ nets for 
all $m, n \geq 3$, with flexions in closed form;
the letters of Section~\ref{sec:examples} and the digits and animal
silhouettes of Section~\ref{sm:ex} show how different their
shapes can be. Section~\ref{sec:application} designs three further
nets of Proposition~\ref{prop:ex:flex} as transportable arch halls of
$80$~mm panels and checks them numerically for transport, clear height,
and snow load; a $1{:}20$ physical model of one of them shows its motion
(Figure~\ref{fig:app:model}).
The supplement~\cite{Nurmatov2026Repo}
verifies the symbolic computations, constructs and checks the nets, and
contains the print files and the assembly guide of the physical model. Left
for future work are a method that 
detects self-intersections exactly, not only numerically;
assemblies with other GQS-nets or with subnets of general
equimodular elliptic type; and nets that approximate whole
three-dimensional shapes or pass from one prescribed shape to another.

\section*{Acknowledgments}
This research has been supported by the baseline funding of the KAUST Computational Sciences Group.

\section*{Declarations}
During the preparation of this work, the authors used the Claude models
Opus~5, Opus~5.5, Fable~5, and Fable~5.1 (Anthropic) to search the literature;
double-check citations against the cited sources; improve the language and
readability of the text; double-check the authors' proofs for errors;
check, revise, and help write the code of the supplement; help with
rendering the schematic Figures~\ref{fig:01ab}--\ref{fig:04} and the nets
of Figures~\ref{fig:ex:nets:GQS}, \ref{fig:app:render}, \ref{fig:ex:nets:2026},
and~\ref{fig:ex:nets:animals} in Blender; help with designing the
physical model of Section~\ref{sm:app:model}, including its print files
and assembly guide; and help with producing the video of the supplement,
whose narration was generated with Gemini~3.8 Flash TTS (Google). The
authors reviewed, edited, and verified all content, including all code,
computations, figures, and references. The authors assume responsibility for all content.

\bibliographystyle{abbrvurl}
\bibliography{010_references}

% ---------------------------------------------------------------------------------------------------------
% The supplementary materials, in the same document (arXiv version): SM numbers for the pages, sections, figures, tables,
% equations, and theorems, the title "Supplementary Materials: ...", and a reference list of their own ([SM1]--[SM12]), so
% that all cross-references between the paper and the supplementary materials are links.
\clearpage
\makeatletter
\renewcommand{\thepage}{SM\arabic{page}}\setcounter{page}{1}
\setcounter{section}{0}\setcounter{figure}{0}\setcounter{table}{0}\setcounter{equation}{0}\setcounter{theorem}{0}\setcounter{footnote}{0}
\renewcommand{\thesection}{SM\arabic{section}}
\renewcommand{\thefigure}{SM\arabic{figure}}
\renewcommand{\thetable}{SM\arabic{table}}
\begin{center}
{\LARGE Supplementary Materials: \TheTitle\par}
\vskip 1.5em
{\large\TheAuthors\par}
\end{center}
\vskip 1.5em
% Citations: the labels [SMn] and links to the reference list below; nothing is written to the .aux file.
\def\@extra@b@citeb{@SM}
\def\@citex[#1]#2{%
  \leavevmode
  \let\@citea\@empty
  \@cite{%
    \@for\@citeb:=#2\do{%
      \@citea
      \def\@citea{,\penalty\@m\ }%
      \edef\@citeb{\expandafter\@firstofone\@citeb\@empty}%
      \@ifundefined{b@\@citeb\@extra@b@citeb}{%
        \hbox{\reset@font\bfseries ?}%
        \G@refundefinedtrue
        \@latex@warning{Citation `\@citeb' of the supplementary materials undefined}%
      }{%
        \hyper@@link[cite]{}{cite.\@citeb\@extra@b@citeb}{\csname b@\@citeb\@extra@b@citeb\endcsname}%
      }%
    }%
  }{#1}%
}
% Citation labels of the supplementary materials ([SM1]--[SM12]), from a separate BibTeX run (abbrvurl) over their citations.
\expandafter\gdef\csname b@EU96-53@SM\endcsname{SM1}
\expandafter\gdef\csname b@CraneRepo@SM\endcsname{SM2}
\expandafter\gdef\csname b@Dhondt2004@SM\endcsname{SM3}
\expandafter\gdef\csname b@EN1991-1-3@SM\endcsname{SM4}
\expandafter\gdef\csname b@ISO668@SM\endcsname{SM5}
\expandafter\gdef\csname b@KSLL2test@SM\endcsname{SM6}
\expandafter\gdef\csname b@KSLL2testDE@SM\endcsname{SM7}
\expandafter\gdef\csname b@KSLL2@SM\endcsname{SM8}
\expandafter\gdef\csname b@MichelinXTA2@SM\endcsname{SM9}
\expandafter\gdef\csname b@Nurmatov2026Repo@SM\endcsname{SM10}
\expandafter\gdef\csname b@StanfordRepo@SM\endcsname{SM11}
\expandafter\gdef\csname b@Zenkert@SM\endcsname{SM12}

\makeatother

% The text of the supplementary materials: input by supplementary.tex (the separate PDF) and by main.tex (the arXiv version, after the paper).
These supplementary materials give the data of the nets of
Section~\ref{sec:ex:algorithm}, further examples, and the details of the
arch halls of Section~\ref{sec:application}. The notebooks, the code,
the data, and the video are in the repository~\cite{Nurmatov2026Repo}, which the paper
calls its supplement. All renderings, in
Figures~\ref{fig:01ab}--\ref{fig:app:render}, \ref{fig:ex:nets:2026},
and~\ref{fig:ex:nets:animals}, were made with
Blender~5.2.1.

\section{Examples of flexible \texorpdfstring{$m \times n$}{m x n} nets}
\label{sm:ex}
All nets below are constructed as in Section~\ref{sec:ex:algorithm} and
checked numerically with \texttt{gqs.py} as described there; this tool
also lets one set all the data of the construction and shows the motion of
a net interactively. Its settings for the pictures of these nets and of
Figure~\ref{fig:ex:nets:GQS}, including the views, are listed in the
repository (\texttt{parameters.txt})~\cite{Nurmatov2026Repo}.
Table~\ref{tab:ex:nets:GQS} gives the data of the nets of
Figure~\ref{fig:ex:nets:GQS}.

\begin{table}[htbp!]
\centering
\footnotesize
\caption{The nets of Figure~\ref{fig:ex:nets:GQS}, with the signs
$(e_1, e_2, e_3, e_4) = (1, -1, -1, 1)$ and $e_0 = 1$ for~(i),
$(1, -1, -1, 1)$ and $e_0 = -1$ for~(ii), and $(1, -1, 1, -1)$ and
$e_0 = 1$ for~(iii): $(\alpha_1, \beta_1, \gamma_1, \delta_1)$ are the flat
angles at~$A_1$ of the block~$(\kappa_0, \nu_0)$, and
$\rho_1, \ldots, \rho_{10}$ are the systems (a)--(d) of
Lemma~\ref{lem:ex:patterns} in the assembly.}
\label{tab:ex:nets:GQS}
\setlength{\tabcolsep}{3pt}
\renewcommand{\arraystretch}{1.2}
\begin{tabular}{c|c|c|c|c|c}
Fig. & $M$ & $(\alpha_1, \beta_1, \gamma_1, \delta_1)$ & $m \times n$ & $(\kappa_0, \nu_0)$ & $(\rho_1, \ldots, \rho_{10})$ \\
\hline
(i)   & $4/3$        & $(45^\circ, 45^\circ, 90^\circ, 150^\circ)$      & $5 \times 12$ & $(1, 1)$ & abababbaaa \\
(ii)  & $\sqrt{3}/2$ & $(142.5^\circ, 127.5^\circ, 15^\circ, 60^\circ)$ & $5 \times 12$ & $(1, 7)$ & bbababaaaa \\
(iii) & $\sqrt{3}/2$ & $(75^\circ, 22.5^\circ, 112.5^\circ, 120^\circ)$ & $5 \times 12$ & $(1, 3)$ & baabbabbba
\end{tabular}
\end{table}

\subsection{The digits 2026}
\label{sm:ex:2026}
\begin{figure}[htbp!]
    \centering
    \includegraphics[width=0.9\textwidth]{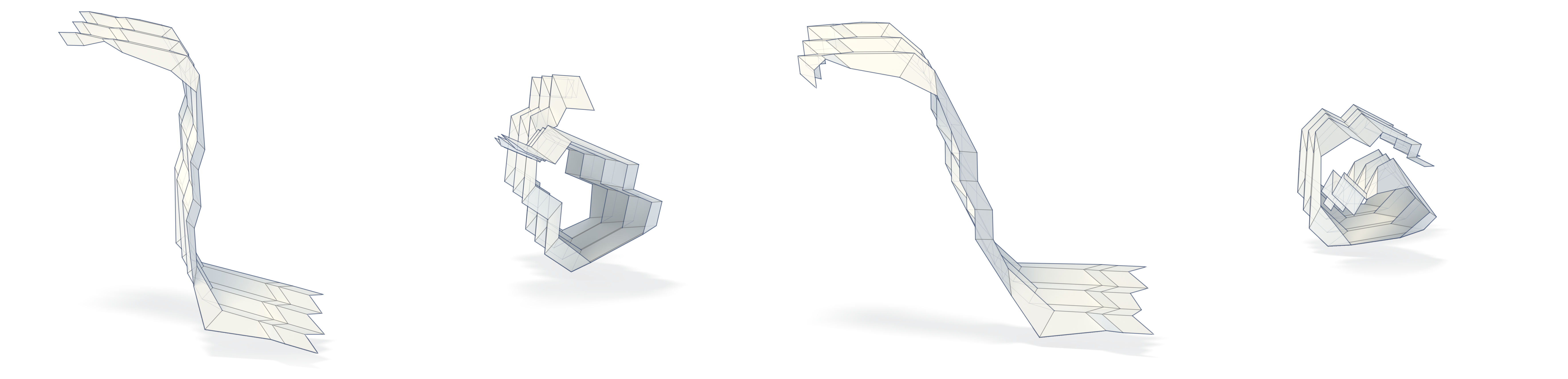}\\[-0.2ex]
    {\footnotesize%
    \makebox[0.225\textwidth]{(i)~$t = 0.06$}%
    \makebox[0.225\textwidth]{(ii)~$t = 0.45$}%
    \makebox[0.225\textwidth]{(iii)~$t = 1.65$}%
    \makebox[0.225\textwidth]{(iv)~$t = 0.5$}}
    \\[1.5ex]
    \includegraphics[width=0.9\textwidth]{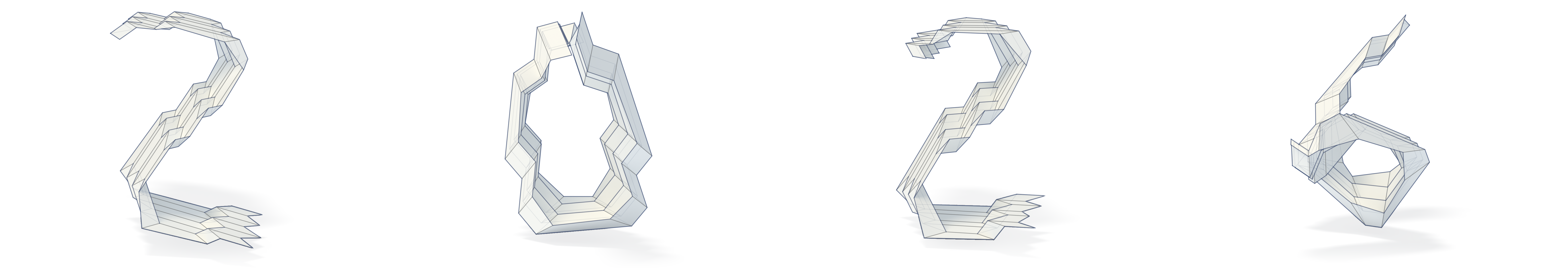}\\[-0.2ex]
    {\footnotesize%
    \makebox[0.225\textwidth]{(i)~$t = 0.24$}%
    \makebox[0.225\textwidth]{(ii)~$t = 0.16$}%
    \makebox[0.225\textwidth]{(iii)~$t = 2.93$}%
    \makebox[0.225\textwidth]{(iv)~$t = 0.1$}}
    \caption{Four flexible nets (i)--(iv) assembled from GQS-nets, each shown
    at two values of the parameter~$t$ of the block~$(\kappa_0, \nu_0)$, from
    the same viewing direction and at the same scale. Bottom: the nets read 2026; top:
    the same nets at another position of their motion, which is checked
    numerically to pass from one to the other without self-intersections. 
    The data of the nets are given in
    Table~\ref{tab:ex:nets:2026}.}
    \label{fig:ex:nets:2026}
\end{figure}

Figure~\ref{fig:ex:nets:2026} shows four nets that form the digits ``2026'', and
Table~\ref{tab:ex:nets:2026} gives their data.

\begin{table}[htbp!]
\centering
\footnotesize
\caption{The nets of Figure~\ref{fig:ex:nets:2026}, all with the signs
$(e_1, e_2, e_3, e_4) = (1, -1, 1, -1)$ and $e_0 = 1$:
$(\alpha_1, \beta_1, \gamma_1, \delta_1)$ are the flat angles at~$A_1$ of
the block~$(\kappa_0, \nu_0)$, and $\rho_1, \ldots, \rho_{12}$ are the
systems (a)--(d) of Lemma~\ref{lem:ex:patterns} in the assembly.}
\label{tab:ex:nets:2026}
\setlength{\tabcolsep}{3pt}
\renewcommand{\arraystretch}{1.2}
\begin{tabular}{c|c|c|c|c|c}
Fig. & $M$ & $(\alpha_1, \beta_1, \gamma_1, \delta_1)$ & $m \times n$ & $(\kappa_0, \nu_0)$ & $(\rho_1, \ldots, \rho_{12})$ \\
\hline
(i)   & $\sqrt{2}$   & $(135^\circ, 22.5^\circ, 142.5^\circ, 82.5^\circ)$ & $6 \times 14$ & $(1, 1)$ & aabbbaaabbba \\
(ii)  & $\sqrt{3/2}$ & $(127.5^\circ, 67.5^\circ, 120^\circ, 135^\circ)$  & $6 \times 14$ & $(1, 1)$ & cabaabbbaaba \\
(iii) & $5/4$        & $(72^\circ, 144^\circ, 144^\circ, 108^\circ)$      & $6 \times 14$ & $(1, 1)$ & cabbbaaabbba \\
(iv)  & $4/3$        & $(45^\circ, 90^\circ, 45^\circ, 30^\circ)$         & $6 \times 14$ & $(1, 1)$ & cbbbbbbbaaba
\end{tabular}
\end{table}

\begin{figure}[htbp!]
    \centering
    \includegraphics[width=0.9\textwidth]{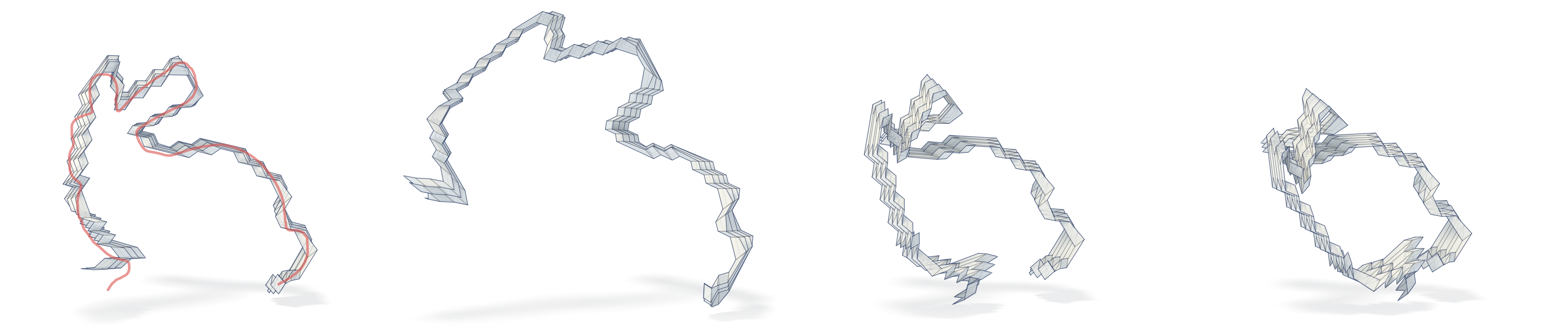}\\[-0.2ex]
    {\footnotesize%
    \makebox[0.225\textwidth]{(i)~$t = 0.5$}%
    \makebox[0.225\textwidth]{$t = 0.093$}%
    \makebox[0.225\textwidth]{$t = 0.71$}%
    \makebox[0.225\textwidth]{$t = 0.83$}}
    \\[1.5ex]
    \includegraphics[width=0.9\textwidth]{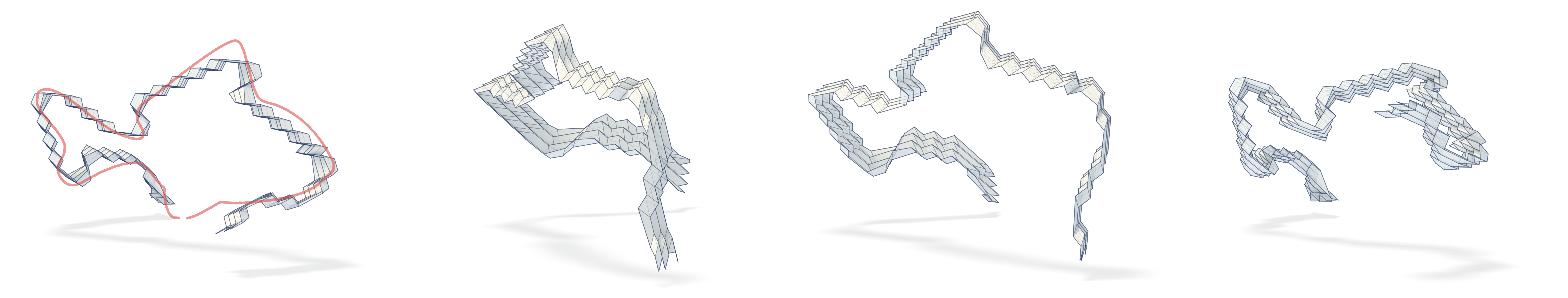}\\[-0.2ex]
    {\footnotesize%
    \makebox[0.225\textwidth]{(ii)~$t = 0.96$}%
    \makebox[0.225\textwidth]{$t = 0.15$}%
    \makebox[0.225\textwidth]{$t = 0.49$}%
    \makebox[0.225\textwidth]{$t = 1.21$}}
    \caption{Two flexible $6 \times 60$ nets (i), (ii) assembled from GQS-nets, each shown
    at four values of the parameter~$t$ of the block~$(\kappa_0, \nu_0)$, from
    the same viewing direction and at the same scale. Left: the nets approximately follow the side
    silhouettes (red) of the Stanford bunny~\cite{StanfordRepo} and of Blub the
    fish~\cite{CraneRepo}; right: the same nets at three other positions of their
    motion, which is checked numerically to pass through all four positions without self-intersections.}
    \label{fig:ex:nets:animals}
\end{figure}

\subsection{Side silhouettes}
\label{sm:ex:animals}
The systems $\rho_1, \ldots, \rho_{n-2}$ of an assembly can be
chosen freely, and this freedom can be used to make a net approximately follow a prescribed shape: the profile of a long net (Section~\ref{sec:application}) depends
on the choice of the systems~$\rho_\nu$. As the target
we take the side silhouette of a mesh: the outer contour of its side view,
cut open at the bottom. For a fixed~$t$, we then choose
$\rho_1, \rho_2, \ldots$ one after another by a beam search: every partial
net kept so far is extended by one more row of blocks, with each of the
systems (a)--(d) of Lemma~\ref{lem:ex:patterns}, its profile is compared, up to a similarity, with the corresponding part of the silhouette, and only
the closest few partial nets are kept. Repeating this
for several values of~$t$ and of the angles at~$A_1$, and keeping only nets
whose motion between the positions shown is checked numerically to be free of self-intersections, we get the two
$6 \times 60$ nets of Figure~\ref{fig:ex:nets:animals}, which approximately
follow the side silhouettes (red) of the Stanford
bunny~\cite{StanfordRepo} and of Blub the fish~\cite{CraneRepo}; each is
shown in this shape and at three other positions of its motion. Both nets
have the same block~$(\kappa_0, \nu_0) = (1, 1)$, the same angles
$(\alpha_1, \beta_1, \gamma_1, \delta_1) = (45^\circ, 90^\circ, 45^\circ, 30^\circ)$
at~$A_1$ of this block (hence the same $M = 4/3$), and the same signs
$e_1, \ldots, e_4$, $e_0$ as the net~(iv) of
Figure~\ref{fig:ex:nets:2026}; only $n$ and the systems~$\rho_\nu$ differ,
and $\rho_1, \ldots, \rho_{58}$ are listed in the repository
(\texttt{parameters.txt})~\cite{Nurmatov2026Repo}.

\section{The arch halls}
\label{sm:app}
This section gives the details behind Section~\ref{sec:application}. The
scripts that produce or check each value of
Sections~\ref{sm:app:nets}--\ref{sm:app:snow} and~\ref{sm:app:render} are
listed in the repository (\texttt{Application}, Sections~4 and~5 of
\texttt{README.md})~\cite{Nurmatov2026Repo}.

\begin{table}[htbp!]
\centering
\footnotesize
\caption{The halls (i)--(iii), shown in
Figure~\ref{fig:app:render} and Figures~\ref{fig:app:profiles} and~\ref{fig:app:snow}, all with
$(e_1, e_2, e_3, e_4) = (1, -1, 1, -1)$ and $e_0 = 1$; $t$ is the
parameter of the block~$(\kappa_0, \nu_0)$; $\ell$: the length of the
profile, that is, of the midline of Figure~\ref{fig:app:profiles};
$L$: the overall length along the hall;
$H_{\mathrm{c}}$: the largest clear height above the floor, with the feet at the
given height; $W_{20}$: the width of the widest floor strip with at least
$2.0$~m of clear height; $s_{\mathrm{k},\max}$ (kN/m$^2$): the capacity, the
smallest ground snow load at which a limit on displacement or stress is
reached; $\lambda_1$: the buckling factor of the first global mode at
$s_{\mathrm{k}} = 1$~kN/m$^2$.}
\label{tab:app}
\setlength{\tabcolsep}{3pt}
\renewcommand{\arraystretch}{1.2}
\begin{tabular}{l|c|c|c}
 & (i) & (ii) & (iii) \\
\hline
profile; transport position & parabola-like; low-bed & parabola-like; loop & circle-like; tube \\
$(\alpha_1, \beta_1, \gamma_1, \delta_1)$ & $(36^\circ, 144^\circ, 72^\circ, 72^\circ)$ & $(120^\circ, 60^\circ, 60^\circ, 60^\circ)$ & $(60^\circ, 60^\circ, 120^\circ, 60^\circ)$ \\
$M$ & $5/4$ & $9/8$ & $9/8$ \\
$m \times n$, $(\kappa_0, \nu_0)$ & $20 \times 20$, $(2, 1)$ & $13 \times 19$, $(1, 1)$ & $20 \times 27$, $(1, 1)$ \\
$\rho_1, \ldots, \rho_{n-2}$ & baaacabaaccabacaaa & bccccccabaccccccb & aaacabacacacccacacabacaaa \\
faces $F_{ij}$ of the hall & $1 \leq i \leq 18$, $2 \leq j \leq 18$ & $1 \leq i \leq 11$, $1 \leq j \leq 17$ & $1 \leq i \leq 18$, $1 \leq j \leq 25$ \\
panels; hinges & $306$; $577$ & $187$; $346$ & $450$; $857$ \\
typical panel edge (m) & $0.61$ & $0.53$ & $0.33$ \\
$t$: deployed; transport & $0.08$; $0.32$ & $0.7265$; $1.4019$ & $0.15$; $0.75$ \\
$\ell$ (m): deployed; transport & $7.81$; $8.06$ & $8.79$; $9.07$ & $7.48$; $7.70$ \\
deployed span $\times$ rise $\times$ $L$ (m) & $4.00 \times 2.84 \times 11.65$ & $4.50 \times 2.54 \times 5.19$ & $3.96 \times 2.40 \times 5.27$ \\
transport $W \times H \times L$ (m) & $2.35 \times 3.27 \times 11.91$ & $2.15 \times 1.76 \times 5.72$ & $1.77 \times 2.23 \times 5.72$ \\
feet above the floor (m) & $0.08$ & $0.30$ & $0.50$ \\
$H_{\mathrm{c}}$; $W_{20}$ (m) & $2.74$; $2.00$ & $2.71$; $2.29$ & $2.79$; $2.30$ \\
$s_{\mathrm{k},\max}$; $\lambda_1$, feet anchored & $0.60$; $1.26$ & $0.64$; $3.12$ & $0.68$; $2.66$ \\
$s_{\mathrm{k},\max}$; $\lambda_1$, spindle and lock & $0.43$; $1.51$ & $0.04$; $2.07$ & $0.04$; $2.03$
\end{tabular}
\end{table}

\begin{figure}[htbp!]
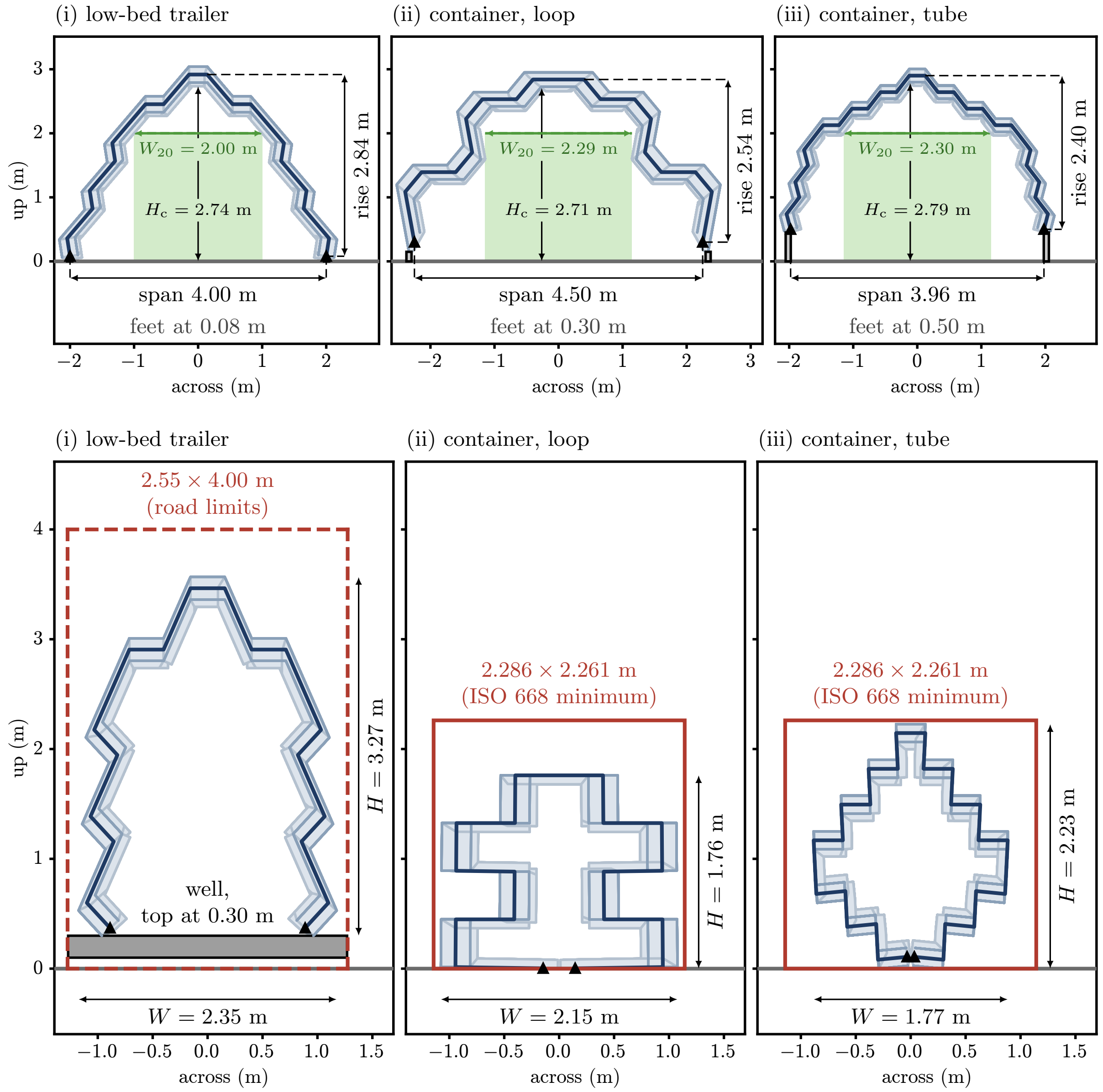

    \centering
    \input{fig/fig_app_deployed.tex}\\[6pt]
    \input{fig/fig_app_transport.tex}
    \caption{Profiles of the halls (i)--(iii), seen along the rows of
blocks. Top: deployed positions, with the $80$~mm panels (light); the
polylines of the vertices~$V_{ij}$ for each~$i$ (thin; they coincide for
all even and for all odd~$i$); the midline (dark), the mean of the even
and odd polylines, with its ends on the foot lines (triangles); and the
knee walls (grey). Green: the widest floor strip with at least $2.0$~m of
clear height; $H_{\mathrm{c}}$ and~$W_{20}$ as in Table~\ref{tab:app}. Bottom:
transport positions, drawn the same way, (i) in the well (grey) of the
semi-trailer within the road limits of~\cite[Annex~I, 1.2(a) and
1.3]{EU96-53} (red, dashed), and (ii) and~(iii) within the minimum door
opening of~\cite{ISO668} (red). $W$ and~$H$: width and height with the
panels.}
    \label{fig:app:profiles}
\end{figure}

\subsection{Nets and panels}
\label{sm:app:nets}
Table~\ref{tab:app} gives the data of the halls. The net of~(i) was found by a
beam search that fixes $\rho_1, \rho_2, \ldots$ one after another, followed
by a local search, and those of (ii) and~(iii) by an enumeration of the
sequences with $\rho_\nu = \rho_{n-1-\nu}$. For~(iii), the systems (a)
and~(c) were then exchanged in~$\rho_\nu$ for $\nu = 4, 8, 10, 12, 14, 16, 18, 22$. This leaves the profile unchanged in both positions and reduces
the skew in plan, the spread over~$j$ of the mean positions along the hall
of the vertices $V_{1j}, \ldots, V_{19j}$ of the hall, from $2.73$ to $0.63$~m (deployed)
and from $2.45$ to $0.56$~m (transport).
The panels are hinged on their
outer face sheets, along the edges of the net.
At each hinge, both panels are bevelled by the bisector
plane at the position of smallest inner angle; the corners of the panels
of~(i) are cut at $90$ vertices. The faces of the bare nets are free of
self-intersections at $25$ positions between the two values of~$t$. With
the panels, at $365$, $181$, and $259$ positions of the motions of (i),
(ii), and~(iii), no two panels overlap; the smallest distance between
panels without a common vertex is $188$~mm for~(i), between $19.8$
and~$22.8$~mm for~(ii), and $35$~mm for~(iii), each at the transport
position.

\subsection{Transport and clear height}
\label{sm:app:transport}
The hall~(i) travels in the well of an extendable low-bed
semi-trailer~\cite{KSLL2} and reaches $3.57$~m above the road. The trailer is $13.21$~m long when closed and is extended by $5.55$~m to
$18.76$~m~\cite{KSLL2} for this hall; the articulated vehicle then exceeds
the limits of $16.50$~m of~\cite[Annex~I, 1.1]{EU96-53} and of $12.00$~m
from the king pin to the rear~\cite[Annex~I, 1.6]{EU96-53}, and may
circulate on a special permit for indivisible
loads~\cite[Article~4(3)]{EU96-53} or under a national
derogation~\cite[Article~4(4)]{EU96-53}; whether a folded hall is an
indivisible load is not examined. The halls (ii) and~(iii) fit the floor
of the container, at least $5.867$~m long~\cite{ISO668}.

Inside, the halls meet the targets of height and width
(Table~\ref{tab:app} and Figure~\ref{fig:app:profiles}, top), on plinths
for~(i) and on knee walls for (ii) and~(iii). The knee walls are plain
walls and ship flat; the plinths are plain strips $0.03$~m high. For~(ii),
the clear height $H_{\mathrm{c}} = 2.71$~m is reached in the notches between
the terraces of its crown; under the flat parts of the crown, it is
$2.66$~m.

\subsection{Mobility and drive}
\label{sm:app:drive}
In the drive we propose, one spindle between the feet drives the
motion, which has one degree of freedom: the compatibility matrix of the
rigid panels has a null space of dimension~$7$ (the six rigid motions and
one more) at $9$ equally spaced values of~$t$ from the deployed to the
transport position, both included, in each hall, and of dimension~$6$ with
one hinge locked, for each of five hinges at the deployed, the middle, and
the transport position. One foot keeps its place across the hall, and the other
runs on rails across it. Both feet also run on carriages along the hall,
because the halls lengthen during folding (the period-two translation
$V_{i+2,j} - V_{ij}$ along the rows of blocks grows by $3$, $14$, and $12\%$,
and the length~$L$ of Table~\ref{tab:app}, measured over the stepped open
ends, by $2$, $10$, and $9\%$). Hence anchoring a foot
along its length blocks the motion; in use, the first foot is anchored. The
forces of this drive are not computed here.

The special position $t = 1/\sqrt{2}$ of~(iii) (Section~\ref{sec:application})
lies between its deployed and transport positions, $t = 0.15$ and~$0.75$,
and Theorem~\ref{th:main} excludes it: for its block~$(\kappa_0, \nu_0) = (1, 1)$,
\eqref{eq:N3} gives $x_1 = -4$ and $z_1 = z_2 = 2$, so that the denominators
$y_2 u (z_2 + x_1 t^2)$ and $y_1 u (z_1 + x_1 t^2)$ of the formulas of
Theorem~\ref{th:main} for $\cot(\theta_2/2)$ and $\cot(\theta_4/2)$ vanish at
this~$t$. Through this~$t$, the dihedral
angles at $204$ of the $857$ hinges change sign, and the two faces at each
of them are coplanar there; at all other hinges, the face normals make at
least $70.52^\circ$. At this~$t$, the null space above has dimension~$18$.
None of the $9$ samples lies within $0.001$ of this~$t$; whether another
motion branches off there and whether the drive passes through it are not
examined.

\begin{figure}[htbp!]
    \centering
    \input{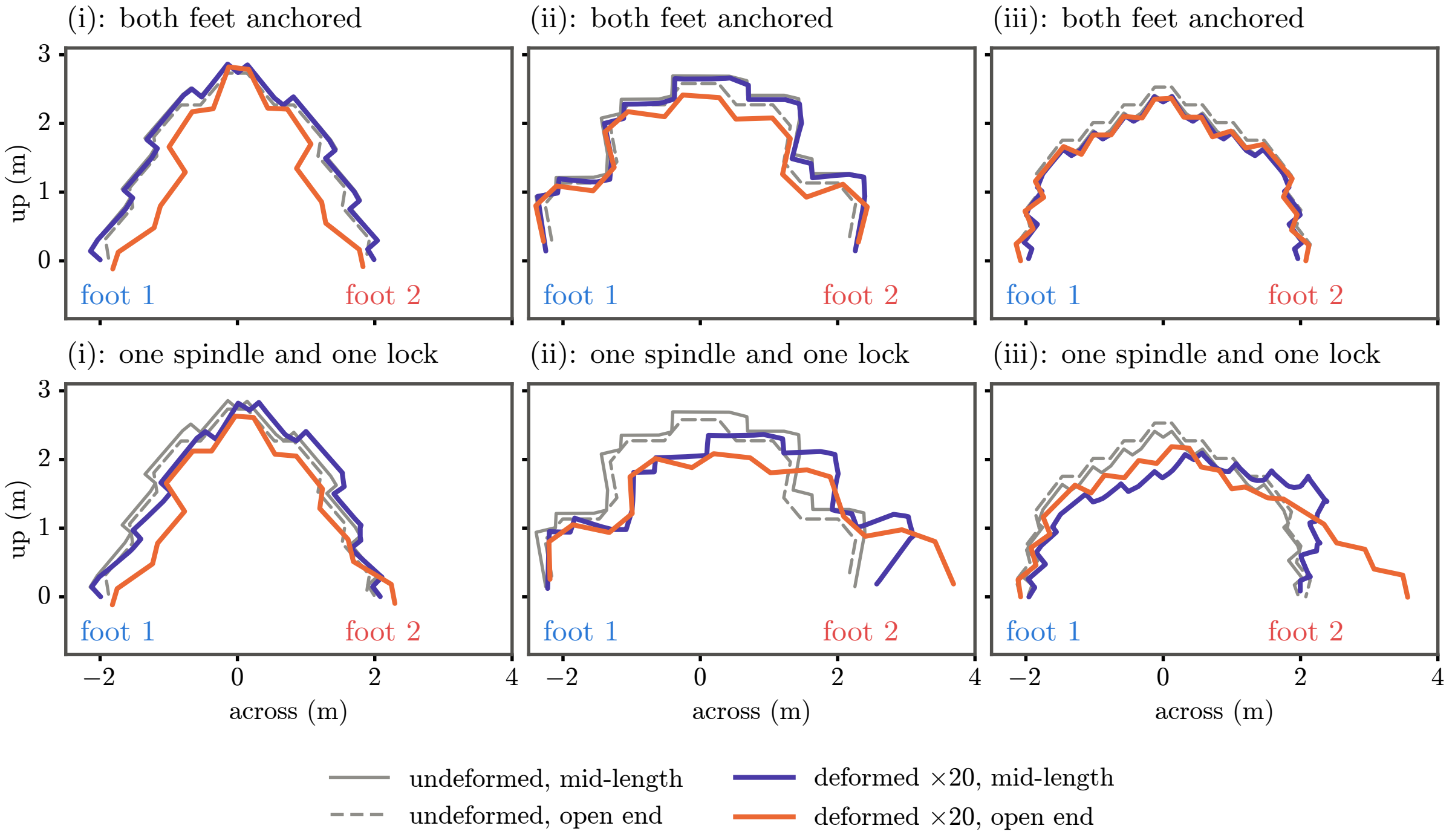}
    \caption{Profiles of the halls (i)--(iii) under self-weight and
undrifted snow with $s_{\mathrm{k}} = 1$~kN/m$^2$, at mid-length and at the free
edge of the weaker open end, with the displacements scaled by~$20$.
Top: both feet anchored; bottom: foot~2 on rails with one spindle and
one locked hinge.}
    \label{fig:app:snow}
\end{figure}

\subsection{Snow check}
\label{sm:app:snow}
We check the deployed halls under snow with a linear shell model of the
sandwich, including the shear of the core. The loads are the self-weight
and the undrifted and drifted snow on a cylindrical roof
of~\cite[Section~5.3.5]{EN1991-1-3} for a ground snow load
$s_{\mathrm{k}} = 1$~kN/m$^2$. Under these loads, we limit the displacement to
span$/200$. Under $1.35$ times the self-weight plus $1.5$ times the
snow, we limit the stresses in the face sheets and the core, limit the
compression of the face sheets also to their wrinkling
stress~\cite[Eq.~(6.14)]{Zenkert}, and
require a buckling factor $\lambda_1 \geq 1$ for the global modes. Local
modes at the point supports, the panel corners, the end of the
spindle, and the latch of a locked hinge depend on the mesh and on the
detailing of these connections, belong to their design, and are not
counted, for example the buckling factors $0.76$ at a support of~(ii) under
drifted snow and $0.57$ at the end of the spindle of~(i). We computed with
CalculiX~2.23~\cite{Dhondt2004} and quadratic shell elements, $12$ per
panel edge for (i) and~(ii) and $8$ for the larger model of~(iii). Around
the panel corners and the supports, where the point connections
concentrate the stresses, the stresses are read outside a radius of
$0.2$~m for (i) and~(ii) and $0.15$~m for~(iii) around every net vertex,
and, for~(i) with one locked hinge, also outside $0.15$~m around its latch.

The panel has steel face sheets of $0.6$~mm ($E_{\mathrm{f}} = 210$~GPa,
yield stress $280$~MPa) and a polyurethane core of $78.8$~mm
($E_{\mathrm{c}} = 4$~MPa, $G_{\mathrm{c}} = 3.5$~MPa, shear strength
$0.10$~MPa); these are assumed values, typical of such panels. The shape
coefficient is $0.8$ for undrifted snow and two triangles with the peaks
$1.0$ and $2.0$ for drifted snow; the wrinkling stress is
$0.5\,(E_{\mathrm{f}}E_{\mathrm{c}}G_{\mathrm{c}})^{1/3} = 71.6$~MPa, and
the stress limit of the face sheets is $280/1.1 = 254.5$~MPa, with an
assumed material factor $1.1$.

Table~\ref{tab:app} gives the capacities and buckling factors with both
feet anchored and with one foot on rails, one spindle, and one locked
hinge; Figure~\ref{fig:app:snow} shows the deformed profiles. All
buckling factors exceed~$1$. With both feet anchored, the halls carry
$0.60$ to $0.68$~kN/m$^2$, about two thirds of the $s_{\mathrm{k}}$ of the check,
set by the sway of the open ends of~(i) under undrifted snow and by the
wrinkling of the face sheets of (ii) and~(iii) under drifted snow. 
A coarser mesh, $8$ instead of $12$ elements per panel edge, changes the
capacities of (i) and~(ii) by at most $10\%$; $6$ instead of $8$ changes
the vertical displacements of~(iii) by at most $2.5\%$. With one spindle, the foot on the rails of (ii)
and~(iii), a hinged zigzag held only at the end of the spindle, spreads
by up to $20$~mm under self-weight and $86$~mm under self-weight and snow, which sets
their capacity. For~(i), wrinkling along a hinge line under drifted snow
sets it, outside $0.15$~m around the latch of the locked hinge, where the
factored self-weight alone exceeds the wrinkling limit. For~(i) with both feet anchored, the wrinkling limit alone gives the
capacities $0.82$, $0.69$, and $0.58$~kN/m$^2$ at the radii $0.2$, $0.15$,
and $0.1$~m.

The shell model assumes ideal hinges without play,
point supports, and rails that hold the feet down; it is geometrically
linear, without wind or imperfections. The wrinkling stress and the
shear strength of the core have no material factor, and the wrinkling
stress assumes full core support, also next to the hinges.

\subsection{Physical model}
\label{sm:app:model}
We built a $1{:}20$ physical model of the hall~(ii) of Table~\ref{tab:app}
(Figure~\ref{fig:app:model}). Its $187$ faces, those of the hall~(ii) at
$1{:}20$, are $4$~mm thick, that is, $80$~mm at full size, and extend to
both sides of the net, whereas the panels of Section~\ref{sec:application}
lie on its inner side; in both, the hinge axes lie on the edges of the net.
We printed the faces as separate parts in polylactic acid and joined them
with $640$ printed pins at the $346$ hinges. The hinges turn freely and have
little play: by design, pins of $2.8$~mm diameter turn in bores $0.05$~mm
wider, with axial gaps of $0.05$ to $0.06$~mm. A stand holds the feet at
their deployed distance, and the model stands on it by its own weight; off
the stand, it folds by hand into a nearly closed loop like that of the
transport position and opens again. The model shows the motion of the net
qualitatively, without the knee walls and the drive of the hall. Along the
computed motion from the deployed to the transport position, the printed
parts are checked numerically at $181$ positions: they keep their design
clearances, and faces not joined by a hinge stay at least $1.6$~mm apart.
The print files, the assembly guide, and the result of this check are in
the repository (\texttt{Prototype})~\cite{Nurmatov2026Repo}.

\subsection{Renderings}
\label{sm:app:render}
The trailer (K\"assbohrer K.SLL 2/2H-18/20, extended) has the width
$2.54$~m and the well top at $0.30$~m from~\cite{KSLL2}, the well thickness
and the tyre size from~\cite{KSLL2test}, the deck height over the axles
from~\cite{KSLL2testDE}, and the tyre diameter from~\cite{MichelinXTA2};
the container has the minimum internal size, the minimum door opening, and
the nominal external size of~\cite{ISO668}. In the values of Figure~\ref{fig:app:render}, the rise is measured to the
crown of the midline; the outer top of the panels lies $0.12$, $0.11$,
and $0.10$~m above it. The deployed length is measured over the stepped
open ends. The plinths and knee walls are $80$~mm thick and centred under
the lowest panel corners. The packed tube of~(iii) leaves $35$~mm between
its feet. The arrows and dashed lines of the values are hidden behind opaque
parts and dimmed behind the translucent walls and roof of the container.
Under the feet of (i)--(iii) are triangular gaps (seats, not designed),
and the ground shows through the open ends. Sub-pixel offsets of $0.1$ to
$2$~mm keep touching faces apart.

\begingroup\makeatletter\@fileswfalse\renewcommand\@biblabel[1]{[SM#1]}
% The reference list of the supplementary materials (abbrvurl), from a separate BibTeX run over their citations.

\endgroup

\end{document}